\documentclass[12pt,a4paper]{article}
\usepackage[active]{srcltx}
\usepackage{amsmath,amsfonts,amssymb,epsf}
\usepackage[T1]{fontenc}
\usepackage{verbatim} 
\usepackage[colorlinks=true,linkcolor=blue,urlcolor=blue]{hyperref}

\usepackage{pstricks-add}

\newcommand{\cP}{\mathcal{P}}

\newcommand{\cL}{\mathcal{L}}

\def\R{{\mathbb R}}

\def\N{{\mathbb N}}

\def\K{{\mathbb K}}
\def\Z{{\mathbb Z}}

\def\inv{{}^{-1}}

\def\Ga{{\rm GA}}
\def\msk{\medskip}

\def\ssk{\smallskip}

\def\giantbreak{\par \ifdim\lastskip<2\bigskipamount \removelastskip

         \penalty-400 \giantskip\fi}

\def\nin{\noindent}

\def\pagebreak{\vskip 0pt plus 0.0001fil\break}
\def\linebreak{\break}

\def\Gl{{\rm GL}}

\def\1{{\bf 1}}

\def\theta{{\vartheta}}

\def\cL{{\cal L}}

\makeatletter
\def\blfootnote{\xdef\@thefnmark{}\@footnotetext}
\makeatother

\begin{document}

\title{Commutative and Non-commutative Parallelogram Geometry:  an Experimental Approach}

\author{Wolfgang Bertram}

\date{  }

    \maketitle
    
 \begin{abstract}   
 By ``parallelogram geometry'' we mean the elementary, ``commutative'', geometry corresponding to
 vector addition, and by  ``trapezoid geometry'' a certain  ``non-commutative deformation'' of
 the former.
This text presents an elementary
 approach via exercises using dynamical software, hopefully accessible to a wide mathematical
audience, from undergraduate students and high school teachers to researchers,
 proceeding in three steps:
 
 \ssk
(1) {\em experimental geometry}, 

(2) {\em algebra (linear algebra and elementary group theory)},
and

(3)  {\em axiomatic geometry}. 
\end{abstract}

\blfootnote{

keywords: projective geometry, group, torsor, associativity, Desargues theorem, lattice, geogebra

2010 AMS classification: 
06C05,  
08A02, 
51A05,  
97G70,  	
97U70  	
}

\bigskip
{\bf Introduction.}
Sometimes, fundamental research leads to elementary results that can 
easily be explained to a wide audience of non-specialists, and which already bear the germ of
the much more sophisticated mathematics lying at the background.
In the present work I will try to show that our
joint paper with M.\ Kinyon  \cite{BeKi12}  is an example of this situation, and
 that even in such classical domains 
as plane geometry of points and lines, something new can be said and, maybe, taught. 
In a nutshell, everything is contained in the following figure:
%
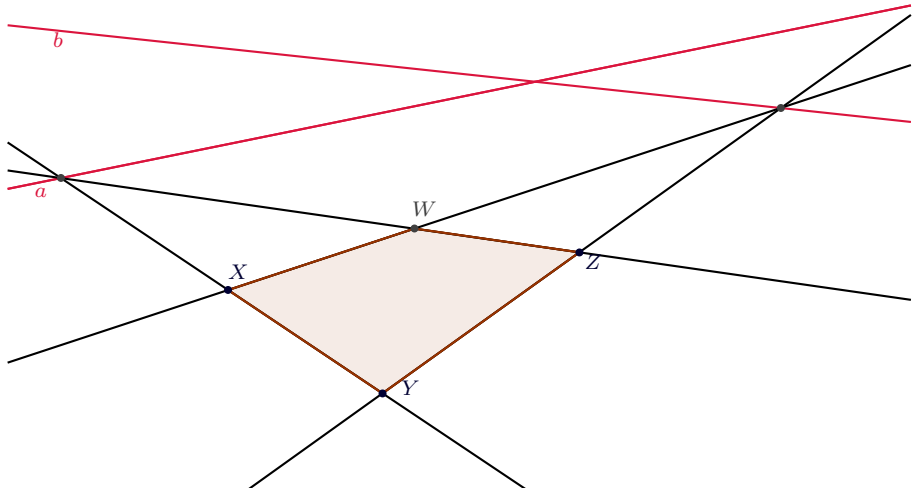
\begin{center}
\newrgbcolor{dcrutc}{0.86 0.08 0.24}
\newrgbcolor{qqqqtt}{0 0 0.2}
\newrgbcolor{zzttqq}{0.6 0.2 0}
\psset{xunit=0.5cm,yunit=0.5cm,algebraic=true,dotstyle=o,dotsize=3pt 0,linewidth=0.8pt,arrowsize=3pt 2,arrowinset=0.25}
\begin{pspicture*}(-4.3,-7.32)(19.46,6.3)
\pspolygon[linecolor=zzttqq,fillcolor=zzttqq,fillstyle=solid,opacity=0.1](1.5,-2.08)(6.4,-0.45)(10.74,-1.08)(5.56,-4.82)
\psplot[linecolor=dcrutc]{-4.3}{19.46}{(--22.3--3.08*x)/15.04}
\psplot[linecolor=dcrutc]{-4.3}{19.46}{(--84.64-2.04*x)/18.9}
\psplot[linecolor=dcrutc]{-4.3}{19.46}{(-22.3-3.08*x)/-15.04}
\psplot{-4.3}{19.46}{(--4.33--2.74*x)/-4.06}
\psplot{-4.3}{19.46}{(-45.76--3.74*x)/5.18}
\psplot{-4.3}{19.46}{(--6.42-1.97*x)/13.64}
\psplot{-4.3}{19.46}{(-37.48--4.83*x)/14.54}
\psline[linecolor=zzttqq](1.5,-2.08)(6.4,-0.45)
\psline[linecolor=zzttqq](6.4,-0.45)(10.74,-1.08)
\psline[linecolor=zzttqq](10.74,-1.08)(5.56,-4.82)
\psline[linecolor=zzttqq](5.56,-4.82)(1.5,-2.08)
\begin{scriptsize}
\rput[bl](-3.58,0.38){\dcrutc{$a$}}
\rput[bl](-3.1,4.36){\dcrutc{$b$}}
\psdots[dotstyle=*,linecolor=qqqqtt](5.56,-4.82)
\rput[bl](6.06,-4.86){\qqqqtt{$Y$}}
\psdots[dotstyle=*,linecolor=qqqqtt](1.5,-2.08)
\rput[bl](1.5,-1.78){\qqqqtt{$X$}}
\psdots[dotstyle=*,linecolor=qqqqtt](10.74,-1.08)
\rput[bl](10.88,-1.52){\qqqqtt{$Z$}}
\psdots[dotstyle=*,linecolor=darkgray](-2.9,0.89)
\psdots[dotstyle=*,linecolor=darkgray](16.04,2.75)
\psdots[dotstyle=*,linecolor=darkgray](6.4,-0.45)
\rput[bl](6.34,-0.12){\darkgray{$W$}}
\end{scriptsize}
\end{pspicture*}
\end{center}

\nin 
The quadrangle $XYZW$ is not a parallelogram, but its construction has something in common with the one of a
parallelogram: we may call it a ``non-commutative parallelogram''. 
This terminology is  a wink towards ``non-commutative geometry'':  indeed, the picture illustrates the fundamental
process of passing from a commutative, associative law (vector addition, corresponding  to usual
parallelograms) to a non-commutative, but still associative law: the key observation is 
 that, when the lines $a,b$ and the point $Y$ are kept fixed,
  the law given by $(X,Z) \mapsto W$
is associative and non-commutative; moreover it depends ``nicely'' on the parameters $Y,a,b$.
Dynamical software, such as {\tt geogebra}, permits to vizualize such dependence in a very convincing 
way: we construct $W$ as a {\em function} of $X,Y,Z,a,b$, and then directly {\em observe} its behaviour
as the arguments move.  In our opinion, this way of vizualizing
could give students a good intuitive idea of some
basic notions from analysis in a geometric  context, for instance, {\em continuity} and {\em singularity}.
On the other hand, algebraic notions can also visualized: {\em normal subgroups}, 
{\em semidirect products}, as well as their combination with notions from analysis: {\em deformation},
{\em contraction}. 

\msk
This ``experimental approach'' (given in Chapter \ref{sec:1}) automatically leads to formalization and to mathematical
argumentation: to simplify the situation, we suppress the line $b$ from our description by choosing
it as ``line at infinity''; then  the function $(X,Y,Z,a) \mapsto W$ may be studied from the point of 
view of affine analytic geometry -- for vectors $x,y,z,w \in \R^2$ and a linear
form $\alpha:\R^2 \to \R$ having the line $a$ as kernel, this function is given by the formula
(Chapter \ref{sec:Compute}) 
\[
w  = \alpha(z)\alpha(y)\inv (x-y) + z \, . \leqno{\rm (A)} 
\]
The same function can be studied from the point of view of ``synthetic'' geometry (Chapter \ref{sec:Geometry}), 
 in terms   of meets of lines ($\land$) and joins of points ($\lor$)
\[
w =  \Bigl( \bigl( (x \lor y) \land a \bigr) \lor z \Bigr)
\land
\Bigl( \bigl( ( z \lor y) \land b \bigr) \lor x \Bigr) \, , \leqno{\rm (B)}
\]    
which places it into the context of {\em lattice theory}.
In Chapter \ref{sec:Compute}
 we present, again in exercise form, the elementary linear algebra related to Formula (A),
and in Chapter \ref{sec:Geometry}  some of the axiomatic aspects related to Formula (B). 
Since most of these exercises are non-standard, but rather easy, we hope that this presentation
could be useful for teachers both at university or high-school looking for fresh and new material
related to otherwise standard topics of the curriculum.

\msk
Finally, both formulas (A) and (B) lead to important generalizations and research topics:
 (A) generalizes to matrix or operator spaces (see section \ref{subsec:GeneralizationII})
and is closely related to the general
framework of associative algebras (including $C^*$-algebras, and thus the framework of
non-commutative geometry mentioned above);  and formula (B)
can be written for any lattice, and thus raises the question of the r\^ole that our function might play
in the general theory of lattices. 
The case best understood so far is the case of lattices from {\em projective geometries}:
for the specialist it will not be surprising that, axiomatically, {\em Desargues' theorem} is at the root of associativity
of the law $(X,Z) \mapsto W$ (Chapter \ref{sec:Geometry}); 
but this law preserves many of its interesting properties (though not
associativity) for certain {\em non-Desarguesian projective planes}  (see \cite{BeKi12}), and thus leads to the
realm of geometries related to {\em non-associative algebras}, such as {\em alternative} or {\em Jordan
algebras} (since {\em symmetric spaces} play an important r\^ole in this context, we take, in the
experimental part,  a closer look on {\em inversions} and {\em powers}).  
Some more remarks on these and other  topics are given in the text, 
and the reader might agree that elementary questions can lead
quickly to deep and non-trivial mathematical problems.

          
\newpage          
          
           \setcounter{tocdepth}{3} 
           
    \tableofcontents
    

\newpage

    \section{Geometry: experimental approach}\label{sec:1}
    
The following series of exercises is designed to give an ``experimental approach to
incidence geometry'': we draw figures only using points, lines joining points, 
intersections of lines, and parallels; we discover properties of these figures
just by looking at dynamical images, but at this stage we do not aim at giving mathematical proofs.   
The first step of these exercices is to install on your computer some interactive geometry
software -- see  \url{http://en.wikipedia.org/wiki/Interactive_geometry_software} for
an overview concerning existing software. 
For our purposes, we found that {\tt geogebra},  freely 
available under the address 
    \url{http://www.geogebra.org/cms/en}, is very convenient (and it has been used to create the images given in the sequel); for this reason we will
use the word  ``geogebra'' instead of the longer ``interactive geometry software''. 
 
 %
    
 \subsection{Incidence geometry: preliminaries}\label{subsec:Incprelim}
 
 Informally,
 {\em incidence geometry} may be defined as the geometry using
 only the following three kinds of the operations available on geogebra:
 
 \begin{itemize}
\item[1.]  {\em line passing through two points}:
   for two distinct points $x,y$ in the plane (or in space), we denote by $x \lor y$
   the affine line spanned by $x$ and $y$,
\item[2.] {\em intersection of two lines}:
   for two distinct and non-parallel lines $u,v$ in the plane, we denote by $u \land v$
   their intersection point,
\item[3.] {\em  parallel line}: for a given line $\ell$ and a point $x$, there is exactly one line
$k$ through $x$ that is parallel to $\ell$.  We will use the notation
$k=x \lor (\ell \land i)$ to denote this unique parallel of $\ell$ through $x$.
\end{itemize}  

The notation used in the preceding point is simply a convention.
The reason to use this convention is that we imagine a point ``infinitely far'' on the line $\ell$,
given by intersecting $a$ with an ``infinitely far line'' $i$; then the parallel to $x$ is the line joining
this infinitely far point $\ell \land i$ with $x$  --  the two horizontal lines in the following
figure will become ``true parallels'': 

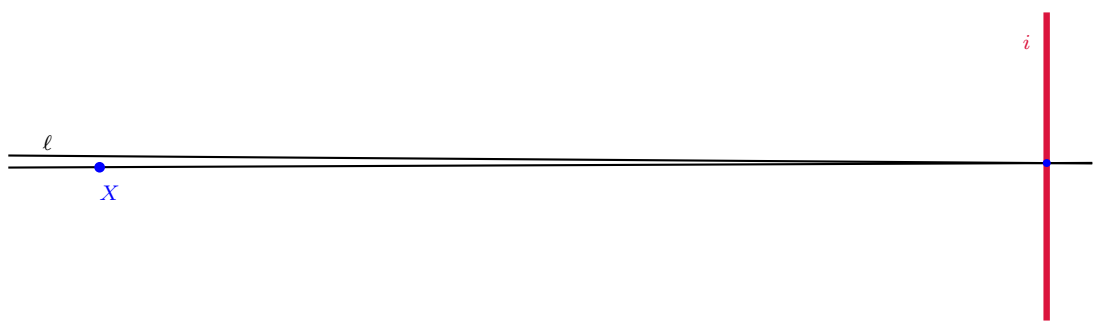
\begin{center}
\newrgbcolor{dcrutc}{0.86 0.08 0.24}
\psset{xunit=0.6cm,yunit=0.3cm,algebraic=true,dotstyle=o,dotsize=3pt 0,linewidth=0.8pt,arrowsize=3pt 2,arrowinset=0.25}
\begin{pspicture*}(-4.3,-7.32)(19.46,6.3)
\psplot{-4.3}{19.46}{(--10.8-0.18*x)/-20.76}
\psplot{-4.3}{19.46}{(--1.79--0.3*x)/-20.36}
\psline[linewidth=2.4pt,linecolor=dcrutc](18.46,-7.32)(18.46,6.3)
\begin{scriptsize}
\psdots[dotstyle=*,linecolor=blue](18.46,-0.36)
\psdots[dotsize=4pt 0,dotstyle=*,linecolor=blue](-2.3,-0.54)
\rput[bl](-2.3,-2.0){\blue{$X$}}
\rput[bl](-3.54,0.2){$\ell$}
\rput[bl](17.94,4.7){\dcrutc{$i$}}
\end{scriptsize}
\end{pspicture*}
\end{center}

We don't bother, for the moment, whether such a line $i$  ``really exists'': it is just a symbol
in a mathematical formula which is executed by geogebra in the way just explained.

 \subsection{Parallelogram geometry}\label{subsec:ParGeo}
 
 \subsubsection{Usual view}\label{sec:ParUsual}
 
 Given three non-collinear points $x,y,z$, construct, using geogebra, a fourth point
 \begin{equation}
w :=  \Bigl( \bigl( (x \lor y) \land i \bigr) \lor z \Bigr)
\land
\Bigl( \bigl( ( z \lor y) \land i \bigr) \lor x \Bigr) \, ,  
\label{eqn:ParUsual}
\end{equation}
that is, according to our convention, the intersection of the parallel
of $x \lor y$ through $z$ with the one of $z \lor y$ through $x$. In other words, construct
the last vertex $w$  in a parallelogram having vertices $x,y,z,w$, like this:

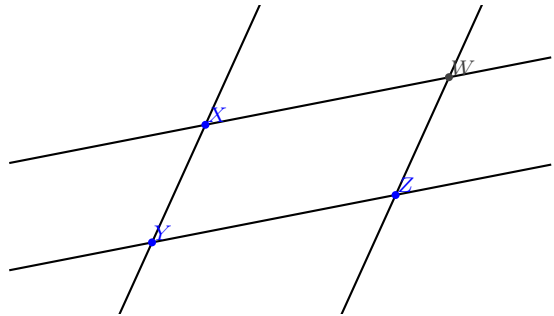
\begin{center}
\psset{xunit=0.3cm,yunit=0.3cm,algebraic=true,dotstyle=o,dotsize=3pt 0,linewidth=0.8pt,arrowsize=3pt 2,arrowinset=0.25}
\begin{pspicture*}(-4.3,-7.32)(19.46,6.3)
\psplot{-4.3}{19.46}{(-19.97--5.2*x)/2.34}
\psplot{-4.3}{19.46}{(-48.76--2.1*x)/10.68}
\psplot{-4.3}{19.46}{(--1.86--2.1*x)/10.68}
\psplot{-4.3}{19.46}{(-70.6--5.2*x)/2.34}
\begin{scriptsize}
\psdots[dotstyle=*,linecolor=blue](1.96,-4.18)
\rput[bl](2.04,-4.06){\blue{$Y$}}
\psdots[dotstyle=*,linecolor=blue](4.3,1.02)
\rput[bl](4.38,1.14){\blue{$X$}}
\psdots[dotstyle=*,linecolor=blue](12.64,-2.08)
\rput[bl](12.72,-1.96){\blue{$Z$}}
\psdots[dotstyle=*,linecolor=darkgray](14.98,3.12)
\rput[bl](15.06,3.24){\darkgray{$W$}}
\end{scriptsize}
\end{pspicture*}
\end{center}

 
Now move one (or several) of the points $x,y,z$, and observe how the position of $w$ 
changes when varying $x,y$ or $z$. 
 What happens to $w$ if you drag the point $x$ onto the line $y \lor z$: does it tend 
to some well-defined point? Note that the initial construction
 works  well if we assume that $x,y,z$ are not collinear, but it  is not defined
  if $x,y,z$ are on a common line, say $\ell$.
Nevertheless, even if $x,y,z$ are collinear, geogebra seems to attribute a well-defined position to
the point $w$ (in fact, a position on the line $\ell$)!  
In the language of analysis, this means that the map associating to the triple
$(x,y,z)$ the fourth point $w$ is {\em continuous}, or, more precisely, that
it admits a {\em continuous extension} from its initial domain of definition
(non-collinear triples) to the bigger set of {\em all} triples.
For geometry, the collinear case is of major importance --
we will investigate below (see \ref{sec:Coll2}) in more detail what is going on in this case.

\subsubsection{Perspective view}\label{sec:ParPersp}
 
Let us have another view on a parallelogram: choose a line
$a$  in the plane (which you draw horizontally)
and three non-collinear points $x,y,z$ in the plane such that the lines
$x \lor y$ and $z \lor y$ are not parallel to $a$. 
Using geogebra, construct the point 
\begin{equation}
w :=  \Bigl( \bigl( (x \lor y) \land a \bigr) \lor z \Bigr)
\land
\Bigl( \bigl( ( z \lor y) \land a \bigr) \lor x \Bigr) \ .
\label{eqn:ParPersp}
\end{equation}
Your image may look like this: 

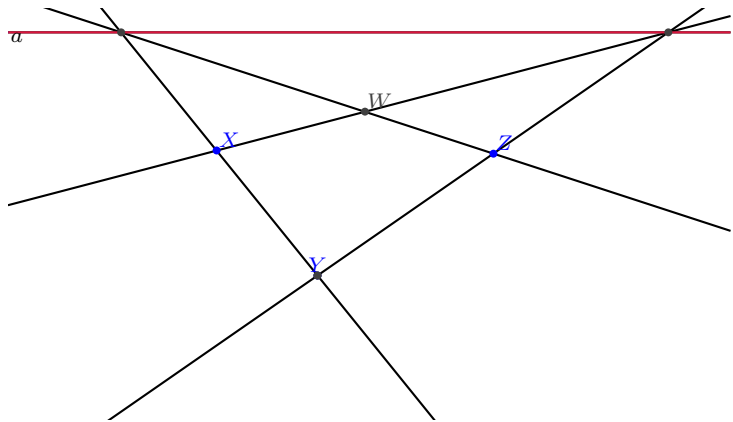
\begin{center}
\newrgbcolor{dcrutc}{0.86 0.08 0.24}
\psset{xunit=0.4cm,yunit=0.4cm,algebraic=true,dotstyle=o,dotsize=3pt 0,linewidth=0.8pt,arrowsize=3pt 2,arrowinset=0.25}
\begin{pspicture*}(-4.3,-7.32)(19.46,6.3)
\psplot{-4.3}{19.46}{(--128.84-0*x)/23.34}
\psplot{-4.3}{19.46}{(-38.44--4.04*x)/5.78}
\psplot{-4.3}{19.46}{(-15.91--4.14*x)/-3.32}
\psplot{-4.3}{19.46}{(-13.73-3.92*x)/-14.85}
\psplot{-4.3}{19.46}{(-65.24--4.02*x)/-12.24}
\psplot[linecolor=dcrutc]{-4.3}{19.46}{(-128.84-0*x)/-23.34}
\begin{scriptsize}
\psdots[dotstyle=*,linecolor=blue](-4.42,5.52)
\rput[bl](-4.2,5.2){$a$}
\psdots[dotstyle=*,linecolor=blue](5.88,-2.54)
\rput[bl](5.58,-2.42){\blue{$Y$}}
\psdots[dotstyle=*,linecolor=blue](11.66,1.5)
\rput[bl](11.74,1.62){\blue{$Z$}}
\psdots[dotstyle=*,linecolor=blue](2.56,1.6)
\rput[bl](2.64,1.72){\blue{$X$}}
\psdots[dotstyle=*,linecolor=darkgray](17.41,5.52)
\psdots[dotstyle=*,linecolor=darkgray](-0.58,5.52)
\psdots[dotstyle=*,linecolor=darkgray](7.44,2.89)
\rput[bl](7.52,3){\darkgray{$W$}}
\psdots[dotstyle=*,linecolor=darkgray](5.88,-2.54)
\psdots[dotstyle=*,linecolor=darkgray](5.88,-2.54)
\end{scriptsize}
\end{pspicture*}
 \end{center}

You may imagine  this drawing being  a perspective view onto a plane in 3-dimensional space,
where the red line $a$ represents the ``horizon''.
Drag the line $a$ in your drawing further and further up, away from $x,y,z$.
Observe that the view becomes more and more a view ``from above'', and that the
 figure formed by $x,y,z,w$ looks more and more like a ``usual'' parallelogram, which arises
 when $a=i$ is the line at infinity. 
 You may now move the point $x$ (or $z$, or $y$)
 and see what happens, in particular,  in the collinear case, when  $x$ approaches the line $y \lor z$.
 
 \ssk
 There is also another,  rather special, kind of position for $y$,
 having the effect that the lines
$ \bigr( (x \lor y) \land a \bigr) \lor z$ and $ \bigl( ( z \lor y) \land a \bigr) \lor x$ are parallel,
and hence the intersection point $w$ does not exist in our drawing plane --
find out this position! 
For simplicity, to avoid such complications,  in the following subsections on parallelograms, 
we will stick  to ``usual'' parallelograms, as in subsection \ref{sec:ParUsual} above (for 
perspective parallelograms, most of what we are going to say remains true, but the operations
need not be  defined for {\em all} triples $(x,y,z)$.)

\subsubsection{Commutativity, associativity}\label{sec:ParAss}

The fourth vertex $w$ in a parallelogram depends on $x,y,z$, therefore we introduce the notation
$w:= x +_y z$, and sometimes we write $w = (xyz)$.
In other words, we consider $w$ as a function of $(x,y,z)$.
If $y$ is fixed as ``origin'', we often write $o$ instead of $y$, and let $x+z:= x+_o z$. 

\ssk
The law $(x,z) \mapsto x + z$ is commutative: indeed, this  is obvious since
 the operations $\lor$ and $\land$ are symmetric
in both arguments, and formula (\ref{eqn:ParUsual}) is therefore symmetric in $x$ and $z$,
whence $x + z=z+x$. 

\ssk
Exericise: show that   $+$ is associative:
\begin{equation}
(x + u) + v = x + (u+v) \, .
\end{equation}
In order to visualize this property,
 construct the points $(x+u)+v$ and $x+(u+v)$, like in the following illustration;
move the points $x,u,v$ around and notice that the ``last'' vertex of the blue and of
the yellow parallelogram always coincide. Convince yourself that
this remains true even if some of these points are collinear!

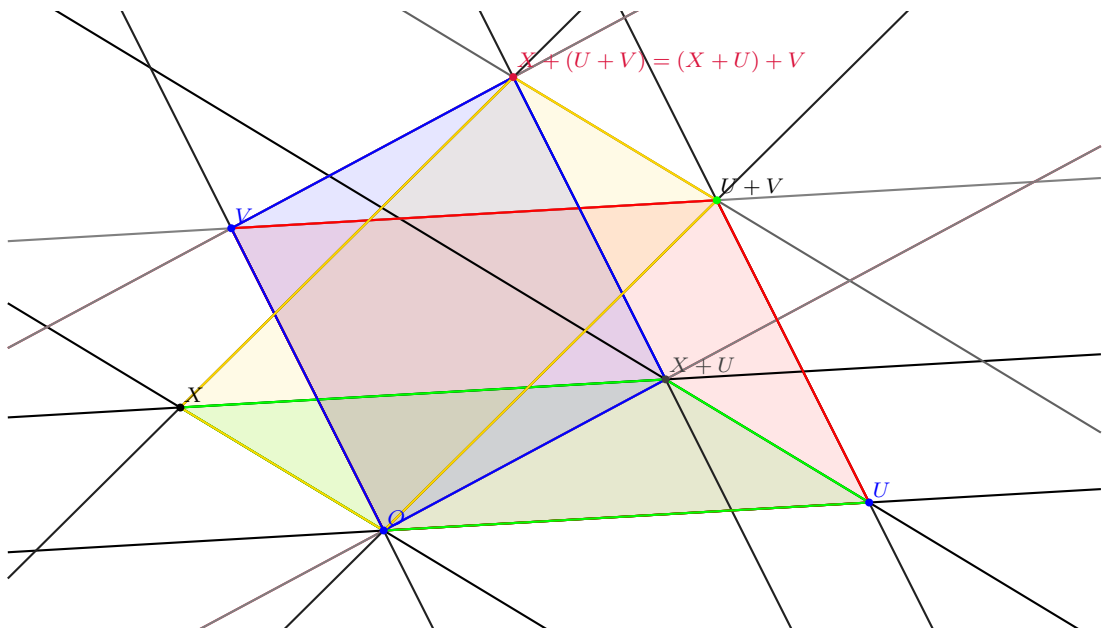
\begin{center}
\newrgbcolor{sqsqsq}{0.13 0.13 0.13}
\newrgbcolor{dcrutc}{0.86 0.08 0.24}
\newrgbcolor{wqwqwq}{0.38 0.38 0.38}
\newrgbcolor{ffdxqq}{1 0.84 0}
\psset{xunit=0.6cm,yunit=0.6cm,algebraic=true,dotstyle=o,dotsize=3pt 0,linewidth=0.8pt,arrowsize=3pt 2,arrowinset=0.25}
\begin{pspicture*}(-4.3,-7.32)(19.66,6.3)
\pspolygon[linecolor=red,fillcolor=red,fillstyle=solid,opacity=0.1](0.6,1.5)(3.94,-5.18)(14.58,-4.56)(11.24,2.12)
\pspolygon[linecolor=green,fillcolor=green,fillstyle=solid,opacity=0.1](3.94,-5.18)(-0.52,-2.46)(10.12,-1.84)(14.58,-4.56)
\pspolygon[linecolor=blue,fillcolor=blue,fillstyle=solid,opacity=0.1](3.94,-5.18)(0.6,1.5)(6.78,4.84)(10.12,-1.84)
\pspolygon[linecolor=ffdxqq,fillcolor=ffdxqq,fillstyle=solid,opacity=0.1](-0.52,-2.46)(3.94,-5.18)(11.24,2.12)(6.78,4.84)
\psplot{-4.3}{19.66}{(--12.39--2.72*x)/-4.46}
\psplot{-4.3}{19.66}{(-57.56--0.62*x)/10.64}
\psplot{-4.3}{19.66}{(-25.85--0.62*x)/10.64}
\psplot{-4.3}{19.66}{(-19.32--2.72*x)/-4.46}
\psplot[linecolor=sqsqsq]{-4.3}{19.66}{(--9.02-6.68*x)/3.34}
\psplot[linecolor=sqsqsq]{-4.3}{19.66}{(--82.16-6.68*x)/3.34}
\psplot[linecolor=gray]{-4.3}{19.66}{(--15.59--0.62*x)/10.64}
\psplot[linecolor=dcrutc]{-4.3}{19.66}{(--45.17-3.34*x)/-6.18}
\psplot[linecolor=dcrutc]{-4.3}{19.66}{(-7.27-3.34*x)/-6.18}
\psplot[linecolor=wqwqwq]{-4.3}{19.66}{(-2927.89--198.96*x)/-326.23}
\psplot[linecolor=gray]{-4.3}{19.66}{(--7.27--3.34*x)/6.18}
\psplot[linecolor=gray]{-4.3}{19.66}{(-45.17--3.34*x)/6.18}
\psplot[linecolor=sqsqsq]{-4.3}{19.66}{(--1948.52-211.8*x)/105.9}
\psline[linecolor=red](0.6,1.5)(3.94,-5.18)
\psline[linecolor=red](3.94,-5.18)(14.58,-4.56)
\psline[linecolor=red](14.58,-4.56)(11.24,2.12)
\psline[linecolor=red](11.24,2.12)(0.6,1.5)
\psline[linecolor=green](3.94,-5.18)(-0.52,-2.46)
\psline[linecolor=green](-0.52,-2.46)(10.12,-1.84)
\psline[linecolor=green](10.12,-1.84)(14.58,-4.56)
\psline[linecolor=green](14.58,-4.56)(3.94,-5.18)
\psline[linecolor=blue](3.94,-5.18)(0.6,1.5)
\psline[linecolor=blue](0.6,1.5)(6.78,4.84)
\psline[linecolor=blue](6.78,4.84)(10.12,-1.84)
\psline[linecolor=blue](10.12,-1.84)(3.94,-5.18)
\psplot[linecolor=sqsqsq]{-4.3}{19.66}{(-66.58--7.3*x)/7.3}
\psplot[linecolor=sqsqsq]{-4.3}{19.66}{(-14.16--7.3*x)/7.3}
\psline[linecolor=ffdxqq](-0.52,-2.46)(3.94,-5.18)
\psline[linecolor=ffdxqq](3.94,-5.18)(11.24,2.12)
\psline[linecolor=ffdxqq](11.24,2.12)(6.78,4.84)
\psline[linecolor=ffdxqq](6.78,4.84)(-0.52,-2.46)
\begin{scriptsize}
\psdots[dotstyle=*,linecolor=blue](3.94,-5.18)
\rput[bl](4.02,-5.06){\blue{$O$}}
\psdots[dotstyle=*](-0.52,-2.46)
\rput[bl](-0.44,-2.34){$X$}
\psdots[dotstyle=*,linecolor=blue](14.58,-4.56)
\rput[bl](14.66,-4.44){\blue{$U$}}
\psdots[dotstyle=*,linecolor=darkgray](10.12,-1.84)
\rput[bl](10.2,-1.72){\darkgray{$X+U$}}
\psdots[dotstyle=*,linecolor=blue](0.6,1.5)
\rput[bl](0.68,1.62){\blue{$V$}}
\psdots[dotstyle=*,linecolor=green](11.24,2.12)
\rput[bl](11.32,2.24){\black{$U+V$}}
\psdots[dotstyle=*,linecolor=dcrutc](6.78,4.84)
\rput[bl](6.86,4.96){\dcrutc{$X+(U+V)=(X+U)+V$}}
\end{scriptsize}
\end{pspicture*}
\end{center}

Make another copy of the same file and add the point $v+x$ to your image;
mark parallelograms as in the following illustration such that a
parallelepiped appears visually, having two of the parallelograms from the preceding 
illustration as ``diagonal parallelograms'':

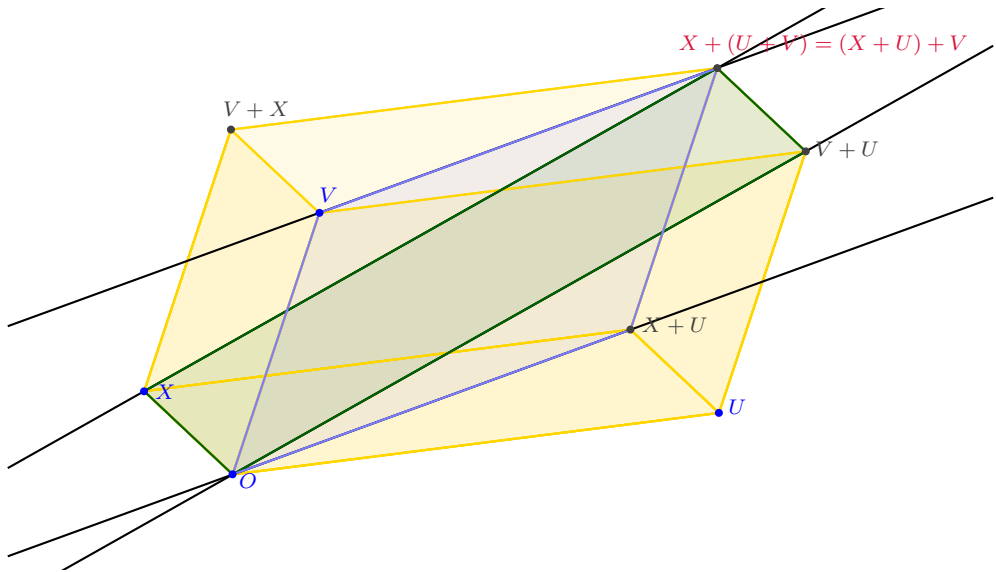
\begin{center}
\newrgbcolor{ffdxqq}{1 0.84 0}
\newrgbcolor{qqwuqq}{0 0.39 0}
\newrgbcolor{xdxdff}{0.49 0.49 1}
\psset{xunit=0.6cm,yunit=0.6cm,algebraic=true,dotstyle=o,dotsize=3pt 0,linewidth=0.8pt,arrowsize=3pt 2,arrowinset=0.25}
\begin{pspicture*}(-3.71,-6.97)(17.89,5.41)
\pspolygon[linecolor=ffdxqq,fillcolor=ffdxqq,fillstyle=solid,opacity=0.1](3.12,0.9)(1.22,-4.88)(-0.72,-3.04)(1.18,2.74)
\pspolygon[linecolor=ffdxqq,fillcolor=ffdxqq,fillstyle=solid,opacity=0.1](1.22,-4.88)(11.88,-3.52)(9.94,-1.68)(-0.72,-3.04)
\pspolygon[linecolor=ffdxqq,fillcolor=ffdxqq,fillstyle=solid,opacity=0.1](13.78,2.26)(11.88,-3.52)(9.94,-1.68)(11.84,4.1)
\pspolygon[linecolor=ffdxqq,fillcolor=ffdxqq,fillstyle=solid,opacity=0.1](1.18,2.74)(11.84,4.1)(9.94,-1.68)(-0.72,-3.04)
\pspolygon[linecolor=ffdxqq,fillcolor=ffdxqq,fillstyle=solid,opacity=0.1](3.12,0.9)(13.78,2.26)(11.88,-3.52)(1.22,-4.88)
\pspolygon[linecolor=qqwuqq,fillcolor=qqwuqq,fillstyle=solid,opacity=0.1](-0.72,-3.04)(1.22,-4.88)(13.78,2.26)(11.84,4.1)
\pspolygon[linecolor=xdxdff,fillcolor=xdxdff,fillstyle=solid,opacity=0.1](3.12,0.9)(1.22,-4.88)(9.94,-1.68)(11.84,4.1)
\psline[linecolor=ffdxqq](3.12,0.9)(1.22,-4.88)
\psline[linecolor=ffdxqq](1.22,-4.88)(-0.72,-3.04)
\psline[linecolor=ffdxqq](-0.72,-3.04)(1.18,2.74)
\psline[linecolor=ffdxqq](1.18,2.74)(3.12,0.9)
\psline[linecolor=ffdxqq](1.22,-4.88)(11.88,-3.52)
\psline[linecolor=ffdxqq](11.88,-3.52)(9.94,-1.68)
\psline[linecolor=ffdxqq](9.94,-1.68)(-0.72,-3.04)
\psline[linecolor=ffdxqq](-0.72,-3.04)(1.22,-4.88)
\psline[linecolor=ffdxqq](13.78,2.26)(11.88,-3.52)
\psline[linecolor=ffdxqq](11.88,-3.52)(9.94,-1.68)
\psline[linecolor=ffdxqq](9.94,-1.68)(11.84,4.1)
\psline[linecolor=ffdxqq](11.84,4.1)(13.78,2.26)
\psline[linecolor=ffdxqq](1.18,2.74)(11.84,4.1)
\psline[linecolor=ffdxqq](11.84,4.1)(9.94,-1.68)
\psline[linecolor=ffdxqq](9.94,-1.68)(-0.72,-3.04)
\psline[linecolor=ffdxqq](-0.72,-3.04)(1.18,2.74)
\psline[linecolor=ffdxqq](3.12,0.9)(13.78,2.26)
\psline[linecolor=ffdxqq](13.78,2.26)(11.88,-3.52)
\psline[linecolor=ffdxqq](11.88,-3.52)(1.22,-4.88)
\psline[linecolor=ffdxqq](1.22,-4.88)(3.12,0.9)
\psplot{-3.71}{17.89}{(-46.46--3.2*x)/8.72}
\psplot{-3.71}{17.89}{(-2.11--3.2*x)/8.72}
\psplot{-3.71}{17.89}{(--70.02-7.14*x)/-12.56}
\psplot{-3.71}{17.89}{(-33.04--7.14*x)/12.56}
\psline[linecolor=qqwuqq](-0.72,-3.04)(1.22,-4.88)
\psline[linecolor=qqwuqq](1.22,-4.88)(13.78,2.26)
\psline[linecolor=qqwuqq](13.78,2.26)(11.84,4.1)
\psline[linecolor=qqwuqq](11.84,4.1)(-0.72,-3.04)
\psline[linecolor=xdxdff](3.12,0.9)(1.22,-4.88)
\psline[linecolor=xdxdff](1.22,-4.88)(9.94,-1.68)
\psline[linecolor=xdxdff](9.94,-1.68)(11.84,4.1)
\psline[linecolor=xdxdff](11.84,4.1)(3.12,0.9)
\begin{scriptsize}
\psdots[dotstyle=*,linecolor=blue](1.22,-4.88)
\rput[bl](1.36,-5.17){\blue{$O$}}
\psdots[dotstyle=*,linecolor=blue](-0.72,-3.04)
\rput[bl](-0.48,-3.22){\blue{$X$}}
\psdots[dotstyle=*,linecolor=blue](11.88,-3.52)
\rput[bl](12.09,-3.55){\blue{$U$}}
\psdots[dotstyle=*,linecolor=blue](3.12,0.9)
\rput[bl](3.12,1.14){\blue{$V$}}
\psdots[dotstyle=*,linecolor=darkgray](9.94,-1.68)
\rput[bl](10.19,-1.77){\darkgray{$X+U$}}
\psdots[dotstyle=*,linecolor=darkgray](1.18,2.74)
\rput[bl](1.01,2.99){\darkgray{$V+X$}}
\psdots[dotstyle=*,linecolor=darkgray](11.84,4.1)
\rput[bl](11.00,4.39){\dcrutc{$X+(U+V)=(X+U)+V$}}
\psdots[dotstyle=*,linecolor=darkgray](13.78,2.26)
\rput[bl](13.99,2.1){\darkgray{$V+U$}}
\end{scriptsize}
\end{pspicture*}
\end{center}

Using three-dimensional geometric intuition, you may find this image more 
convincing than the preceding one: it somehow ``explains'' why the associative law should be
true. However, do not forget that our image really is two-dimensional, and three-dimensional
imagination is just an artifice for understanding it. Note also that the point $v+x$ is not really needed
for the construction.

\subsubsection{Associativity: more general version}\label{sec:ParAss2}

We have seen above that in all constructions the choice of the origin, $o$ or $y$,
was completely arbitrary: any point of the plane could play this r\^ole.
Therefore it will be useful to be able to change the origin freely.
This is facilitated by the following more general version of the associative law
(called the {\em para-associative law}): 
\begin{equation}
x +_o (u +_p v) = (x+_o u) +_p v 
\label{eqn:torsor1}
\end{equation}
where $o$ and $p$ may be different points. Exercice: construct the points
$x +_o (u +_p v)$ and $(x+_o u) +_p v $, in a similar way as above:

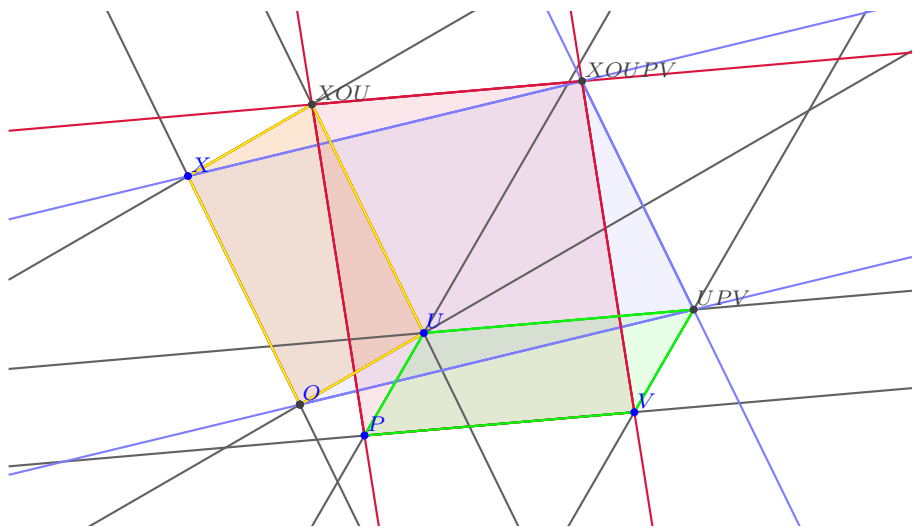
\begin{center} 
\newrgbcolor{wqwqwq}{0.38 0.38 0.38}
\newrgbcolor{dcrutc}{0.86 0.08 0.24}
\newrgbcolor{xdxdff}{0.49 0.49 1}
\psset{xunit=0.5cm,yunit=0.5cm,algebraic=true,dotstyle=o,dotsize=3pt 0,linewidth=0.8pt,arrowsize=3pt 2,arrowinset=0.25}
\begin{pspicture*}(-4.3,-7.32)(19.46,6.3)
\pspolygon[linecolor=xdxdff,fillcolor=xdxdff,fillstyle=solid,opacity=0.1](0.42,1.94)(3.36,-4.12)(13.72,-1.6)(10.78,4.46)
\pspolygon[linecolor=dcrutc,fillcolor=dcrutc,fillstyle=solid,opacity=0.1](5.06,-4.94)(12.16,-4.32)(10.78,4.46)(3.68,3.84)
\pspolygon[linecolor=dcrutc,fillcolor=dcrutc,fillstyle=solid,opacity=0.1](0.42,1.94)(3.68,3.84)(6.62,-2.22)(3.36,-4.12)
\pspolygon[linecolor=green,fillcolor=green,fillstyle=solid,opacity=0.1](5.06,-4.94)(12.16,-4.32)(13.72,-1.6)(6.62,-2.22)
\pspolygon[linecolor=yellow,fillcolor=yellow,fillstyle=solid,opacity=0.1](0.42,1.94)(3.68,3.84)(6.62,-2.22)(3.36,-4.12)
\psplot[linecolor=wqwqwq]{-4.3}{19.46}{(--8.25-6.06*x)/2.94}
\psplot[linecolor=wqwqwq]{-4.3}{19.46}{(-19.82--1.9*x)/3.26}
\psplot[linecolor=wqwqwq]{-4.3}{19.46}{(--5.53--1.9*x)/3.26}
\psplot[linecolor=wqwqwq]{-4.3}{19.46}{(--33.59-6.06*x)/2.94}
\psplot[linecolor=wqwqwq]{-4.3}{19.46}{(-38.21--0.62*x)/7.1}
\psplot[linecolor=wqwqwq]{-4.3}{19.46}{(--21.47-2.72*x)/-1.56}
\psplot[linecolor=wqwqwq]{-4.3}{19.46}{(-19.87--0.62*x)/7.1}
\psplot[linecolor=wqwqwq]{-4.3}{19.46}{(--39.81-2.72*x)/-1.56}
\psplot[linecolor=dcrutc]{-4.3}{19.46}{(--37.61-8.78*x)/1.38}
\psplot[linecolor=dcrutc]{-4.3}{19.46}{(--100.8-8.78*x)/1.38}
\psplot[linecolor=dcrutc]{-4.3}{19.46}{(--633.09--15.71*x)/179.93}
\psplot[linecolor=xdxdff]{-4.3}{19.46}{(-51.15--2.52*x)/10.36}
\psplot[linecolor=xdxdff]{-4.3}{19.46}{(--19.04--2.52*x)/10.36}
\psplot[linecolor=xdxdff]{-4.3}{19.46}{(--1438.95-111.17*x)/53.93}
\psline[linecolor=xdxdff](0.42,1.94)(3.36,-4.12)
\psline[linecolor=xdxdff](3.36,-4.12)(13.72,-1.6)
\psline[linecolor=xdxdff](13.72,-1.6)(10.78,4.46)
\psline[linecolor=xdxdff](10.78,4.46)(0.42,1.94)
\psline[linecolor=dcrutc](5.06,-4.94)(12.16,-4.32)
\psline[linecolor=dcrutc](12.16,-4.32)(10.78,4.46)
\psline[linecolor=dcrutc](10.78,4.46)(3.68,3.84)
\psline[linecolor=dcrutc](3.68,3.84)(5.06,-4.94)
\psline[linecolor=dcrutc](0.42,1.94)(3.68,3.84)
\psline[linecolor=dcrutc](3.68,3.84)(6.62,-2.22)
\psline[linecolor=dcrutc](6.62,-2.22)(3.36,-4.12)
\psline[linecolor=dcrutc](3.36,-4.12)(0.42,1.94)
\psline[linecolor=green](5.06,-4.94)(12.16,-4.32)
\psline[linecolor=green](12.16,-4.32)(13.72,-1.6)
\psline[linecolor=green](13.72,-1.6)(6.62,-2.22)
\psline[linecolor=green](6.62,-2.22)(5.06,-4.94)
\psline[linecolor=yellow](0.42,1.94)(3.68,3.84)
\psline[linecolor=yellow](3.68,3.84)(6.62,-2.22)
\psline[linecolor=yellow](6.62,-2.22)(3.36,-4.12)
\psline[linecolor=yellow](3.36,-4.12)(0.42,1.94)
\begin{scriptsize}
\psdots[dotstyle=*,linecolor=blue](3.36,-4.12)
\rput[bl](3.44,-4){\blue{$O$}}
\psdots[dotstyle=*,linecolor=blue](0.42,1.94)
\rput[bl](0.5,2.06){\blue{$X$}}
\psdots[dotstyle=*,linecolor=blue](6.62,-2.22)
\rput[bl](6.7,-2.1){\blue{$U$}}
\psdots[dotstyle=*,linecolor=blue](5.06,-4.94)
\rput[bl](5.14,-4.82){\blue{$P$}}
\psdots[dotstyle=*,linecolor=blue](12.16,-4.32)
\rput[bl](12.24,-4.2){\blue{$V$}}
\psdots[dotstyle=*,linecolor=darkgray](13.72,-1.6)
\rput[bl](13.8,-1.48){\darkgray{$UPV$}}
\psdots[dotstyle=*,linecolor=darkgray](3.68,3.84)
\rput[bl](3.76,3.96){\darkgray{$XOU$}}
\psdots[dotstyle=*,linecolor=darkgray](3.36,-4.12)
\psdots[dotstyle=*,linecolor=darkgray](3.36,-4.12)
\psdots[dotstyle=*,linecolor=darkgray](10.78,4.46)
\rput[bl](10.86,4.58){\darkgray{$XOUPV$}}
\end{scriptsize}
\end{pspicture*}
\end{center}

(In the drawing we have used the notation $UPV$, etc., instead of $U+_P V$.)
You will note that, wherever you move the points $x,o,u,p,v$, the ``last'' vertex of the blue and
of the red parallelogram always coincide -- this illustrates the para-associative law.
If you drag $o$ to $p$, this gives again the associative law from the preceding illustration. 

\subsubsection{Group law, and the collinear case revisited}\label{sec:Coll2}

We want to show that (for any fixed origin $o$) the plane $P$ with 
$x+z:=x+_o z$ is a {\em commutative group (with neutral element $o$}.
All that remains to be shown is:

\ssk
(a) $o$ is neutral: $x+o=x=o+x$, that is $x+_y y = x = y +_y x$,

(b) for each $x$ there is an element $z$ (called {\em inverse element}) such that $x+z=o$.

\ssk
Using  your dynamical image from subsection \ref{sec:ParUsual}, convince yourself experimentally
that (a) and (b) hold. 
As to (a), this seems to be rather obvious; as to (b), the problem arises that $x,o,z$ will be 
collinear points, and hence our initial definition of $x +_o z$ does not apply. 
Thus, for a better understanding of (b), we have to give, finally,
a good definition of $(xyz)$ for the ``singular'' case when $x,y,z$ are collinear. 
With the help of the para-associative law, this can be done as follows:
first of all,  give an illustration of the special case $u=p$ of the para-associative law:
\begin{equation}
(x+_o p) +_p v = x +_o (p +_p v)  = x +_o v .
\label{eqn:torsor2}
\end{equation}
(We have used for the last equality  that $p+_p v = p$.)  Hence
 the vertex $xoppv$ of the red triangle equals $xov$ and thus does not depend on the choice of $p$:

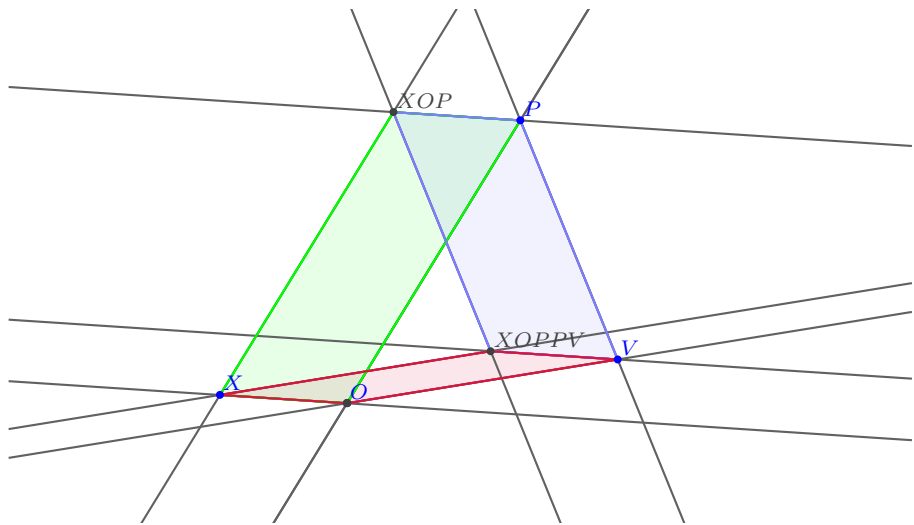
\begin{center}
\newrgbcolor{wqwqwq}{0.38 0.38 0.38}
\newrgbcolor{xdxdff}{0.49 0.49 1}
\newrgbcolor{dcrutc}{0.86 0.08 0.24}
\psset{xunit=0.5cm,yunit=0.5cm,algebraic=true,dotstyle=o,dotsize=3pt 0,linewidth=0.8pt,arrowsize=3pt 2,arrowinset=0.25}
\begin{pspicture*}(-4.3,-7.32)(19.46,6.3)
\pspolygon[linecolor=green,fillcolor=green,fillstyle=solid,opacity=0.1](1.26,-3.92)(5.82,3.58)(9.16,3.36)(4.6,-4.14)
\pspolygon[linecolor=xdxdff,fillcolor=xdxdff,fillstyle=solid,opacity=0.1](5.82,3.58)(9.16,3.36)(11.72,-2.98)(8.38,-2.76)
\pspolygon[linecolor=dcrutc,fillcolor=dcrutc,fillstyle=solid,opacity=0.1](1.26,-3.92)(4.6,-4.14)(11.72,-2.98)(8.38,-2.76)
\psplot[linecolor=wqwqwq]{-4.3}{19.46}{(-12.82-0.22*x)/3.34}
\psplot[linecolor=wqwqwq]{-4.3}{19.46}{(-34.81--1.16*x)/7.12}
\psplot[linecolor=wqwqwq]{-4.3}{19.46}{(-53.38--7.5*x)/4.56}
\psplot[linecolor=wqwqwq]{-4.3}{19.46}{(-27.33--7.5*x)/4.56}
\psplot[linecolor=wqwqwq]{-4.3}{19.46}{(--13.24-0.22*x)/3.34}
\psplot[linecolor=wqwqwq]{-4.3}{19.46}{(--66.68-6.34*x)/2.56}
\psplot[linecolor=wqwqwq]{-4.3}{19.46}{(-53.38--7.5*x)/4.56}
\psplot[linecolor=wqwqwq]{-4.3}{19.46}{(--1200.1-165.18*x)/66.7}
\psplot[linecolor=wqwqwq]{-4.3}{19.46}{(-7.37-0.22*x)/3.34}
\psplot[linecolor=wqwqwq]{-4.3}{19.46}{(-29.37--1.16*x)/7.12}
\psline[linecolor=green](1.26,-3.92)(5.82,3.58)
\psline[linecolor=green](5.82,3.58)(9.16,3.36)
\psline[linecolor=green](9.16,3.36)(4.6,-4.14)
\psline[linecolor=green](4.6,-4.14)(1.26,-3.92)
\psline[linecolor=xdxdff](5.82,3.58)(9.16,3.36)
\psline[linecolor=xdxdff](9.16,3.36)(11.72,-2.98)
\psline[linecolor=xdxdff](11.72,-2.98)(8.38,-2.76)
\psline[linecolor=xdxdff](8.38,-2.76)(5.82,3.58)
\psline[linecolor=dcrutc](1.26,-3.92)(4.6,-4.14)
\psline[linecolor=dcrutc](4.6,-4.14)(11.72,-2.98)
\psline[linecolor=dcrutc](11.72,-2.98)(8.38,-2.76)
\psline[linecolor=dcrutc](8.38,-2.76)(1.26,-3.92)
\begin{scriptsize}
\psdots[dotstyle=*,linecolor=blue](1.26,-3.92)
\rput[bl](1.34,-3.8){\blue{$X$}}
\psdots[dotstyle=*,linecolor=blue](4.6,-4.14)
\rput[bl](4.68,-4.02){\blue{$O$}}
\psdots[dotstyle=*,linecolor=blue](11.72,-2.98)
\rput[bl](11.8,-2.86){\blue{$V$}}
\psdots[dotstyle=*,linecolor=blue](9.16,3.36)
\rput[bl](9.24,3.48){\blue{$P$}}
\psdots[dotstyle=*,linecolor=darkgray](5.82,3.58)
\rput[bl](5.9,3.7){\darkgray{$XOP$}}
\psdots[dotstyle=*,linecolor=darkgray](8.38,-2.76)
\rput[bl](8.46,-2.64){\darkgray{$XOPPV$}}
\psdots[dotstyle=*,linecolor=darkgray](4.6,-4.14)
\psdots[dotstyle=*,linecolor=darkgray](4.6,-4.14)
\psdots[dotstyle=*,linecolor=darkgray](4.6,-4.14)
\psdots[dotstyle=*,linecolor=darkgray](4.6,-4.14)
\end{scriptsize}
\end{pspicture*}
\end{center}

Now flatten the red parallelogram more and more, until $x,o,v$ become collinear.
The fourth point $xov$ is on the same line as $x,o,v$, and it still does not depend on $p$.
Therefore we have found an incidence geometric construction of $xov$ in the collinear case, namely:
{\em choose any point $p$ not on the same line as $x,o,v$, and then construct
first $xop$ and then $(xop)pv$; the result is $xov$:}

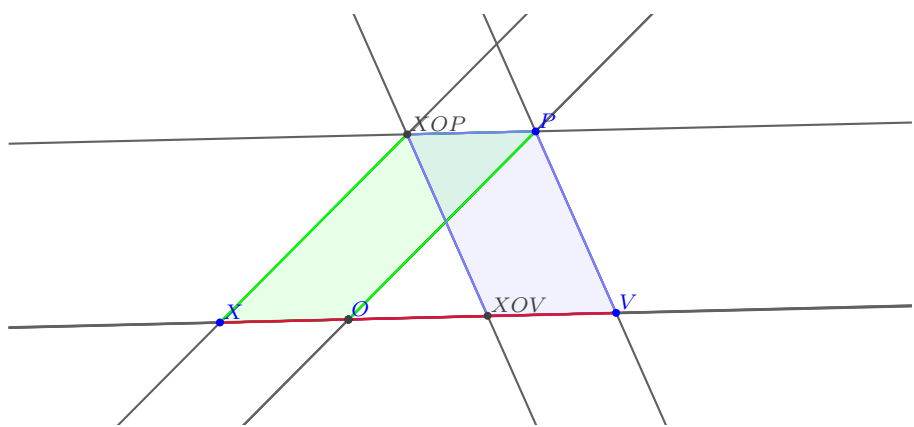
\begin{center}
\newrgbcolor{wqwqwq}{0.38 0.38 0.38}
\newrgbcolor{xdxdff}{0.49 0.49 1}
\newrgbcolor{dcrutc}{0.86 0.08 0.24}
\psset{xunit=0.5cm,yunit=0.4cm,algebraic=true,dotstyle=o,dotsize=3pt 0,linewidth=0.8pt,arrowsize=3pt 2,arrowinset=0.25}
\begin{pspicture*}(-4.3,-7.32)(19.46,6.3)
\pspolygon[linecolor=green,fillcolor=green,fillstyle=solid,opacity=0.1](1.26,-3.92)(6.18,2.32)(9.56,2.42)(4.64,-3.82)
\pspolygon[linecolor=xdxdff,fillcolor=xdxdff,fillstyle=solid,opacity=0.1](6.18,2.32)(9.56,2.42)(11.68,-3.6)(8.3,-3.7)
\pspolygon[linecolor=dcrutc,fillcolor=dcrutc,fillstyle=solid,opacity=0.1](1.26,-3.92)(4.64,-3.82)(11.68,-3.6)(8.3,-3.7)
\psplot[linecolor=wqwqwq]{-4.3}{19.46}{(-13.38--0.1*x)/3.38}
\psplot[linecolor=wqwqwq]{-4.3}{19.46}{(-27.91--0.22*x)/7.04}
\psplot[linecolor=wqwqwq]{-4.3}{19.46}{(-47.75--6.24*x)/4.92}
\psplot[linecolor=wqwqwq]{-4.3}{19.46}{(-27.15--6.24*x)/4.92}
\psplot[linecolor=wqwqwq]{-4.3}{19.46}{(--7.22--0.1*x)/3.38}
\psplot[linecolor=wqwqwq]{-4.3}{19.46}{(--62.68-6.02*x)/2.12}
\psplot[linecolor=wqwqwq]{-4.3}{19.46}{(-47.75--6.24*x)/4.92}
\psplot[linecolor=wqwqwq]{-4.3}{19.46}{(--867.68-124.01*x)/43.67}
\psplot[linecolor=wqwqwq]{-4.3}{19.46}{(-13.34--0.1*x)/3.38}
\psplot[linecolor=wqwqwq]{-4.3}{19.46}{(-27.87--0.22*x)/7.04}
\psline[linecolor=green](1.26,-3.92)(6.18,2.32)
\psline[linecolor=green](6.18,2.32)(9.56,2.42)
\psline[linecolor=green](9.56,2.42)(4.64,-3.82)
\psline[linecolor=green](4.64,-3.82)(1.26,-3.92)
\psline[linecolor=xdxdff](6.18,2.32)(9.56,2.42)
\psline[linecolor=xdxdff](9.56,2.42)(11.68,-3.6)
\psline[linecolor=xdxdff](11.68,-3.6)(8.3,-3.7)
\psline[linecolor=xdxdff](8.3,-3.7)(6.18,2.32)
\psline[linecolor=dcrutc](1.26,-3.92)(4.64,-3.82)
\psline[linecolor=dcrutc](4.64,-3.82)(11.68,-3.6)
\psline[linecolor=dcrutc](11.68,-3.6)(8.3,-3.7)
\psline[linecolor=dcrutc](8.3,-3.7)(1.26,-3.92)
\begin{scriptsize}
\psdots[dotstyle=*,linecolor=blue](1.26,-3.92)
\rput[bl](1.34,-3.8){\blue{$X$}}
\psdots[dotstyle=*,linecolor=blue](4.64,-3.82)
\rput[bl](4.72,-3.7){\blue{$O$}}
\psdots[dotstyle=*,linecolor=blue](11.68,-3.6)
\rput[bl](11.76,-3.48){\blue{$V$}}
\psdots[dotstyle=*,linecolor=blue](9.56,2.42)
\rput[bl](9.64,2.54){\blue{$P$}}
\psdots[dotstyle=*,linecolor=darkgray](6.18,2.32)
\rput[bl](6.26,2.44){\darkgray{$XOP$}}
\psdots[dotstyle=*,linecolor=darkgray](8.3,-3.7)
\rput[bl](8.38,-3.58){\darkgray{$XOV$}}
\psdots[dotstyle=*,linecolor=darkgray](4.64,-3.82)
\psdots[dotstyle=*,linecolor=darkgray](4.64,-3.82)
\psdots[dotstyle=*,linecolor=darkgray](4.64,-3.82)
\psdots[dotstyle=*,linecolor=darkgray](4.64,-3.82)
\end{scriptsize}
\end{pspicture*}
\end{center}
Using this interpretation, show that $xoo=x$ (so $x+o=x$), and that
$z= (oxo)$ is the inverse of $x$:
\begin{center}
\newrgbcolor{wqwqwq}{0.38 0.38 0.38}
\newrgbcolor{xdxdff}{0.49 0.49 1}
\newrgbcolor{dcrutc}{0.86 0.08 0.24}
\psset{xunit=0.5cm,yunit=0.4cm,algebraic=true,dotstyle=o,dotsize=3pt 0,linewidth=0.8pt,arrowsize=3pt 2,arrowinset=0.25}
\begin{pspicture*}(-4.3,-7.32)(19.46,6.3)
\pspolygon[linecolor=green,fillcolor=green,fillstyle=solid,opacity=0.1](1.24,-3.92)(8.3,2.86)(11.7,2.96)(4.64,-3.82)
\psplot[linecolor=wqwqwq]{-4.3}{19.46}{(-13.45--0.1*x)/3.4}
\psplot[linecolor=wqwqwq]{-4.3}{19.46}{(-13.45--0.1*x)/3.4}
\psplot[linecolor=wqwqwq]{-4.3}{19.46}{(-58.43--6.78*x)/7.06}
\psplot[linecolor=wqwqwq]{-4.3}{19.46}{(-36.08--6.78*x)/7.06}
\psplot[linecolor=wqwqwq]{-4.3}{19.46}{(--8.89--0.1*x)/3.4}
\psplot[linecolor=wqwqwq]{-4.3}{19.46}{(--67.32-6.68*x)/-3.66}
\psplot[linecolor=wqwqwq]{-4.3}{19.46}{(-58.43--6.78*x)/7.06}
\psplot[linecolor=wqwqwq]{-4.3}{19.46}{(--1005.04-149.27*x)/-81.79}
\psplot[linecolor=wqwqwq]{-4.3}{19.46}{(-13.45--0.1*x)/3.4}
\psplot[linecolor=wqwqwq]{-4.3}{19.46}{(-13.45--0.1*x)/3.4}
\psline[linecolor=green](1.24,-3.92)(8.3,2.86)
\psline[linecolor=green](8.3,2.86)(11.7,2.96)
\psline[linecolor=green](11.7,2.96)(4.64,-3.82)
\psline[linecolor=green](4.64,-3.82)(1.24,-3.92)
\psline[linecolor=xdxdff](8.3,2.86)(11.7,2.96)
\psline[linecolor=xdxdff](11.7,2.96)(8.04,-3.72)
\psline[linecolor=dcrutc](1.24,-3.92)(4.64,-3.82)
\psline[linecolor=dcrutc](4.64,-3.82)(8.04,-3.72)
\begin{scriptsize}
\psdots[dotstyle=*,linecolor=blue](1.24,-3.92)
\rput[bl](1.22,-3.66){\blue{$X$}}
\psdots[dotstyle=*,linecolor=blue](4.64,-3.82)
\rput[bl](5.04,-3.7){\blue{$O$}}
\psdots[dotstyle=*,linecolor=blue](8.04,-3.72)
\rput[bl](8.32,-3.6){\blue{$Z=OXO$}}
\psdots[dotstyle=*,linecolor=blue](11.7,2.96)
\rput[bl](11.84,2.62){\blue{$P$}}
\psdots[dotstyle=*,linecolor=darkgray](8.3,2.86)
\rput[bl](8.48,2.56){\darkgray{$XOP$}}
\end{scriptsize}
\end{pspicture*}
\end{center}

    \subsection{Trapezoid  geometry}\label{sec:Trapezoid}

Now we are ready to introduce a kind of ``twist'' in the preceding constructions which
will have the effect that we get  {\em non-commutative} group laws. 
 Instead of parallelograms we use trapezoids.
The constructions will depend on the choice of some line $a$ in the plane, to be 
considered as fixed for this section
(drawn horizontally in the images).
The underlying set of our constructions will be the set
 $G:= P \setminus a$  of all points of the plane $P$  not on $a$.

\subsubsection{The generic construction: non-collinear case}\label{subsec:GenericConstruction}

We   recall our convention on the ``line at infinity'' $i$.
For the moment, let us fix also a point $y$ not on $a$, and two other points
$x,z$ such that the line $y \land z$ is not parallel to $a$. 
Construct the point
\begin{equation}
\boxed{
w :=  \Bigl( \bigl( (x \lor y) \land i \bigr) \lor z \Bigr)
\land
\Bigl( \bigl( ( z \lor y) \land a \bigr) \lor x \Bigr) } \ .
\label{eqn:Generic}
\end{equation}
Since the point $w$ depends on $x,y,z$, in formulas will often write
$(xyz)$ for the point $w$ given by (\ref{eqn:Generic}) (and in drawings we will label it  just by $xyz$).
The formula is a sort of hybrid between formulas (\ref{eqn:ParUsual}) and
(\ref{eqn:ParPersp}): in these formulas,  the same line appeared twice, 
 whereas
here {\em  two different} lines appear. Thus 
the figure formed by $x,y,z,w$ is no longer a parallelogram: it is a trapezoid.
Keep the line $a$ fixed and
move the points $x$, $y$ or  $z$ around. Observe. 
Your image may look like this:\footnote{A technical hint: use the ``tool'' command of
geogebra to create a tool doing automatically this construction -- this may be useful 
later on since we are going to use it a lot of times}

\begin{center}
\newrgbcolor{dcrutc}{0.86 0.08 0.24}
\newrgbcolor{wqwqwq}{0.38 0.38 0.38}
\newrgbcolor{uququq}{0.25 0.25 0.25}
\newrgbcolor{aqaqaq}{0.63 0.63 0.63}
\newrgbcolor{qqwuqq}{0 0.39 0}
\psset{xunit=0.5cm,yunit=0.4cm,algebraic=true,dotstyle=o,dotsize=3pt 0,linewidth=0.8pt,arrowsize=3pt 2,arrowinset=0.25}
\begin{pspicture*}(-3.72,-6.97)(17.88,5.41)
\pspolygon[linecolor=aqaqaq,fillcolor=aqaqaq,fillstyle=solid,opacity=0.1](8.22,-1.46)(4.99,-4.17)(0.26,-3.43)(5.6,-1.05)
\pspolygon[linecolor=dcrutc,fillcolor=dcrutc,fillstyle=solid,opacity=0.1](5.6,-1.05)(8.22,-1.46)(4.99,-4.17)(0.26,-3.43)
\pspolygon[linecolor=qqwuqq,fillcolor=qqwuqq,fillstyle=solid,opacity=0.1](5.6,-1.05)(8.22,-1.46)(4.99,-4.17)(0.26,-3.43)
\psplot{-3.72}{17.88}{(--47.16-0.04*x)/24.4}
\psplot[linecolor=dcrutc]{-3.72}{17.88}{(-47.16--0.04*x)/-24.4}
\psplot[linecolor=dcrutc]{-3.72}{17.88}{(-47.16--0.04*x)/-24.4}
\psplot[linecolor=wqwqwq]{-3.72}{17.88}{(--18.94-2.38*x)/-5.35}
\psplot[linecolor=wqwqwq]{-3.72}{17.88}{(-16.02-0.75*x)/4.73}
\psplot[linewidth=0.4pt,linecolor=uququq]{-3.72}{17.88}{(-0.77-0.75*x)/4.73}
\psplot[linecolor=wqwqwq]{-3.72}{17.88}{(--60.66-6.09*x)/-7.26}
\psline[linecolor=aqaqaq](8.22,-1.46)(4.99,-4.17)
\psline[linecolor=aqaqaq](4.99,-4.17)(0.26,-3.43)
\psline[linecolor=aqaqaq](0.26,-3.43)(5.6,-1.05)
\psline[linecolor=aqaqaq](5.6,-1.05)(8.22,-1.46)
\psline[linecolor=dcrutc](5.6,-1.05)(8.22,-1.46)
\psline[linecolor=dcrutc](8.22,-1.46)(4.99,-4.17)
\psline[linecolor=dcrutc](4.99,-4.17)(0.26,-3.43)
\psline[linecolor=dcrutc](0.26,-3.43)(5.6,-1.05)
\psline[linecolor=qqwuqq](5.6,-1.05)(8.22,-1.46)
\psline[linecolor=qqwuqq](8.22,-1.46)(4.99,-4.17)
\psline[linecolor=qqwuqq](4.99,-4.17)(0.26,-3.43)
\psline[linecolor=qqwuqq](0.26,-3.43)(5.6,-1.05)
\begin{scriptsize}
\psdots[dotstyle=*,linecolor=blue](-4.3,1.94)
\psdots[dotstyle=*,linecolor=blue](20.1,1.9)
\psdots[dotstyle=*,linecolor=blue](0.26,-3.43)
\rput[bl](0.33,-3.32){\blue{$Y$}}
\psdots[dotstyle=*,linecolor=blue](5.6,-1.05)
\rput[bl](5.68,-0.94){\blue{$Z$}}
\psdots[dotstyle=*,linecolor=blue](4.99,-4.17)
\rput[bl](5.06,-4.07){\blue{$X$}}
\psdots[dotstyle=*,linecolor=darkgray](-34.1,1.99)
\psdots[dotstyle=*,linecolor=darkgray](12.25,1.91)
\psdots[dotstyle=*,linecolor=darkgray](8.22,-1.46)
\rput[bl](8.29,-1.36){\darkgray{$W=XYZ$}}
\end{scriptsize}
\end{pspicture*}
\end{center}

or like this (crossed trapezoid):

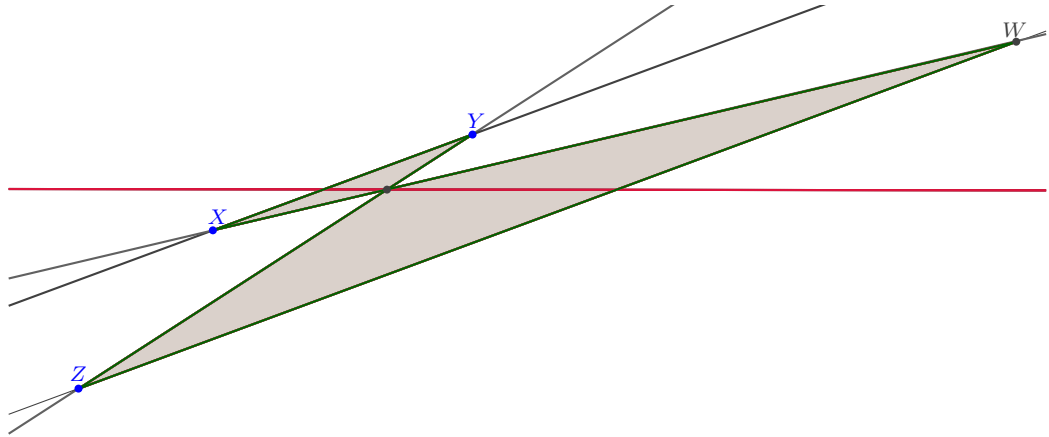
\begin{center}
\newrgbcolor{dcrutc}{0.86 0.08 0.24}
\newrgbcolor{wqwqwq}{0.38 0.38 0.38}
\newrgbcolor{uququq}{0.25 0.25 0.25}
\newrgbcolor{aqaqaq}{0.63 0.63 0.63}
\newrgbcolor{qqwuqq}{0 0.39 0}
\psset{xunit=0.5cm,yunit=0.4cm,algebraic=true,dotstyle=o,dotsize=3pt 0,linewidth=0.8pt,arrowsize=3pt 2,arrowinset=0.25}
\begin{pspicture*}(-6.51,-8.49)(20.78,7.15)
\pspolygon[linecolor=aqaqaq,fillcolor=aqaqaq,fillstyle=solid,opacity=0.1](19.99,5.95)(-1.14,-0.3)(5.69,2.88)(-4.67,-5.55)
\pspolygon[linecolor=dcrutc,fillcolor=dcrutc,fillstyle=solid,opacity=0.1](-4.67,-5.55)(19.99,5.95)(-1.14,-0.3)(5.69,2.88)
\pspolygon[linecolor=qqwuqq,fillcolor=qqwuqq,fillstyle=solid,opacity=0.1](-4.67,-5.55)(19.99,5.95)(-1.14,-0.3)(5.69,2.88)
\psplot{-6.51}{20.78}{(--25.87-0.04*x)/24.4}
\psplot[linecolor=dcrutc]{-6.51}{20.78}{(-25.87--0.04*x)/-24.4}
\psplot[linecolor=dcrutc]{-6.51}{20.78}{(-25.87--0.04*x)/-24.4}
\psplot[linecolor=wqwqwq]{-6.51}{20.78}{(-18.08--8.43*x)/10.36}
\psplot[linecolor=uququq]{-6.51}{20.78}{(-1.57-3.18*x)/-6.83}
\psplot[linewidth=0.4pt,linecolor=uququq]{-6.51}{20.78}{(--23-3.18*x)/-6.83}
\psplot[linecolor=wqwqwq]{-6.51}{20.78}{(-0.17-1.36*x)/-4.58}
\psline[linecolor=aqaqaq](19.99,5.95)(-1.14,-0.3)
\psline[linecolor=aqaqaq](-1.14,-0.3)(5.69,2.88)
\psline[linecolor=aqaqaq](5.69,2.88)(-4.67,-5.55)
\psline[linecolor=aqaqaq](-4.67,-5.55)(19.99,5.95)
\psline[linecolor=dcrutc](-4.67,-5.55)(19.99,5.95)
\psline[linecolor=dcrutc](19.99,5.95)(-1.14,-0.3)
\psline[linecolor=dcrutc](-1.14,-0.3)(5.69,2.88)
\psline[linecolor=dcrutc](5.69,2.88)(-4.67,-5.55)
\psline[linecolor=qqwuqq](-4.67,-5.55)(19.99,5.95)
\psline[linecolor=qqwuqq](19.99,5.95)(-1.14,-0.3)
\psline[linecolor=qqwuqq](-1.14,-0.3)(5.69,2.88)
\psline[linecolor=qqwuqq](5.69,2.88)(-4.67,-5.55)
\begin{scriptsize}
\psdots[dotstyle=*,linecolor=blue](5.69,2.88)
\rput[bl](5.55,3.09){\blue{$Y$}}
\psdots[dotstyle=*,linecolor=blue](-4.67,-5.55)
\rput[bl](-4.9,-5.29){\blue{$Z$}}
\psdots[dotstyle=*,linecolor=blue](-1.14,-0.3)
\rput[bl](-1.27,-0.1){\blue{$X$}}
\psdots[dotstyle=*,linecolor=darkgray](3.44,1.05)
\psdots[dotstyle=*,linecolor=darkgray](19.99,5.95)
\rput[bl](19.63,6.12){\darkgray{$W$}}
\end{scriptsize}
\end{pspicture*}
\end{center}

Note also that, if $y \lor z$ happens to be parallel to $a$, then geogebra completes the
construction in the approriate way, namely by drawing a ``true'' parallelogram.
Finally, you may also move the line $a$:
observe that, if you move the line $a$ ``far away'' from $x,y,z$, then the trapezoid 
converges to a ``true'' parallelogram; thus the parallelogram geometry from the preceding section
can be seen as a ``limit case'' where $a$ approaches the ``line at infinity'' $i$.

\subsubsection{Non-commutativity}\label{subsec:non-com}

One observation that is probably not immediately visible is the fact 
 that $(xyz)$  is {\em not} symmetric in $x$ and $z$:
construct, by using another color, also the point $(zyx)$ and observe that,
except for some special situations (which ones?), $(xyz)$ and $(zyx)$ are different:
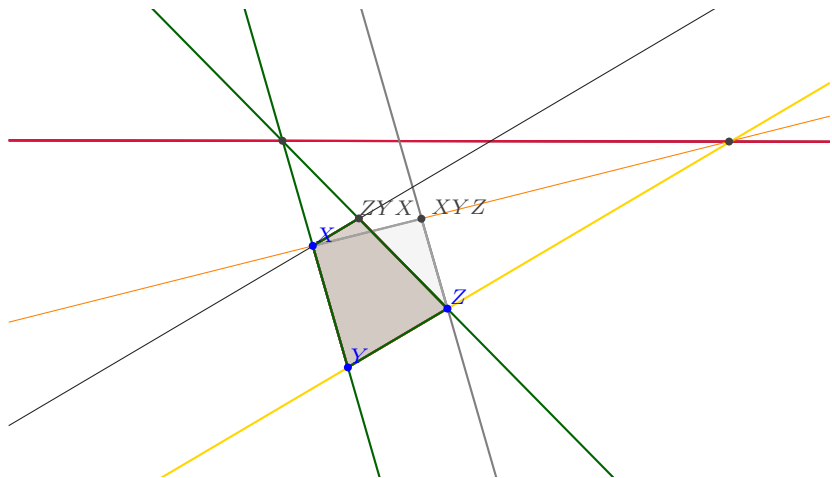
\begin{center} 
\newrgbcolor{dcrutc}{0.86 0.08 0.24}
\newrgbcolor{qqwuqq}{0 0.39 0}
\newrgbcolor{ffdxqq}{1 0.84 0}
\newrgbcolor{ffxfqq}{1 0.5 0}
\newrgbcolor{sqsqsq}{0.13 0.13 0.13}
\newrgbcolor{aqaqaq}{0.63 0.63 0.63}
\psset{xunit=0.5cm,yunit=0.5cm,algebraic=true,dotstyle=o,dotsize=3pt 0,linewidth=0.8pt,arrowsize=3pt 2,arrowinset=0.25}
\begin{pspicture*}(-3.72,-6.97)(17.88,5.41)
\pspolygon[linecolor=aqaqaq,fillcolor=aqaqaq,fillstyle=solid,opacity=0.1](7.14,-0.14)(7.82,-2.52)(5.2,-4.08)(4.28,-0.86)
\pspolygon[linecolor=aqaqaq,fillcolor=aqaqaq,fillstyle=solid,opacity=0.1](5.49,-0.13)(7.82,-2.52)(5.2,-4.08)(4.28,-0.86)
\pspolygon[linecolor=dcrutc,fillcolor=dcrutc,fillstyle=solid,opacity=0.1](4.28,-0.86)(5.49,-0.13)(7.82,-2.52)(5.2,-4.08)
\pspolygon[linecolor=qqwuqq,fillcolor=qqwuqq,fillstyle=solid,opacity=0.1](4.28,-0.86)(5.49,-0.13)(7.82,-2.52)(5.2,-4.08)
\psplot{-3.72}{17.88}{(--47.16-0.04*x)/24.4}
\psplot[linecolor=dcrutc]{-3.72}{17.88}{(-47.16--0.04*x)/-24.4}
\psplot[linecolor=dcrutc]{-3.72}{17.88}{(-47.16--0.04*x)/-24.4}
\psplot[linecolor=qqwuqq]{-3.72}{17.88}{(--13-3.22*x)/0.92}
\psplot[linecolor=gray]{-3.72}{17.88}{(--22.88-3.22*x)/0.92}
\psplot[linecolor=ffdxqq]{-3.72}{17.88}{(-18.82--1.56*x)/2.62}
\psplot[linewidth=0.4pt,linecolor=ffxfqq]{-3.72}{17.88}{(-21.26--2.77*x)/10.95}
\psplot[linewidth=0.4pt,linecolor=sqsqsq]{-3.72}{17.88}{(-8.94--1.56*x)/2.62}
\psplot[linecolor=qqwuqq]{-3.72}{17.88}{(--23.84-4.45*x)/4.34}
\psline[linecolor=aqaqaq](7.14,-0.14)(7.82,-2.52)
\psline[linecolor=aqaqaq](7.82,-2.52)(5.2,-4.08)
\psline[linecolor=aqaqaq](5.2,-4.08)(4.28,-0.86)
\psline[linecolor=aqaqaq](4.28,-0.86)(7.14,-0.14)
\psline[linecolor=aqaqaq](5.49,-0.13)(7.82,-2.52)
\psline[linecolor=aqaqaq](7.82,-2.52)(5.2,-4.08)
\psline[linecolor=aqaqaq](5.2,-4.08)(4.28,-0.86)
\psline[linecolor=aqaqaq](4.28,-0.86)(5.49,-0.13)
\psline[linecolor=dcrutc](4.28,-0.86)(5.49,-0.13)
\psline[linecolor=dcrutc](5.49,-0.13)(7.82,-2.52)
\psline[linecolor=dcrutc](7.82,-2.52)(5.2,-4.08)
\psline[linecolor=dcrutc](5.2,-4.08)(4.28,-0.86)
\psline[linecolor=qqwuqq](4.28,-0.86)(5.49,-0.13)
\psline[linecolor=qqwuqq](5.49,-0.13)(7.82,-2.52)
\psline[linecolor=qqwuqq](7.82,-2.52)(5.2,-4.08)
\psline[linecolor=qqwuqq](5.2,-4.08)(4.28,-0.86)
\begin{scriptsize}
\psdots[dotstyle=*,linecolor=blue](-4.3,1.94)
\psdots[dotstyle=*,linecolor=blue](20.1,1.9)
\psdots[dotstyle=*,linecolor=blue](5.2,-4.08)
\rput[bl](5.28,-3.97){\blue{$Y$}}
\psdots[dotstyle=*,linecolor=blue](4.28,-0.86)
\rput[bl](4.35,-0.76){\blue{$X$}}
\psdots[dotstyle=*,linecolor=blue](7.82,-2.52)
\rput[bl](7.89,-2.41){\blue{$Z$}}
\psdots[dotstyle=*,linecolor=darkgray](15.23,1.91)
\psdots[dotstyle=*,linecolor=darkgray](7.14,-0.14)
\rput[bl](7.4,-0.02){\darkgray{$XYZ$}}
\psdots[dotstyle=*,linecolor=darkgray](3.48,1.93)
\psdots[dotstyle=*,linecolor=darkgray](5.49,-0.13)
\rput[bl](5.50,-0.03){\darkgray{$ZYX$}}
\end{scriptsize}
\end{pspicture*}
\end{center}

Although the product is not symmetric in $x$ and $z$, there are, however, certain 
other ``symmetry relations'':
let $R$ be the set of quadruples $(x,y,z,w)$ of pairwise distinct points in $G$ such that 
$w = (xyz)$. Show that $(x,y,z,w) \in R$ if, and only if,
 the lines $x\lor y$ and $z \lor w$ are parallel and
$y\lor z$ and $w \lor x$ meet on $a$. 
Conclude  that the following four conditions are equivalent:
\begin{equation}\label{eqn:Perminvar}
(xyz)=w, \qquad
(zwx)=y, \qquad
(yxw)=z,\qquad
(wzy)=x  \, .
\end{equation}
In other words, $R$ is invariant under {\em double transpositions}
(e.g., exchanging simultaneaously ($x$ and $y$) and ($z$ and $w$), and so on).
We can rephrase this: {\em the relation $R$ is  invariant under the Klein four group
of double transpositions}.
In particular,  letting $w=(xyz)$ in the other three expressions from (\ref{eqn:Perminvar}),
deduce that  
\begin{equation}
\boxed{
(z(xyz)x) = y, \qquad
(yx(xyz))=z, \qquad
 ((xyz)zy)=x .}
 \label{eqn:Chasles1}
 \end{equation}

\subsubsection{Associativity}\label{subsec:Geo-Associativity}

We have seen that the law $(x,z) \mapsto w = (xyz)$ for fixed $y$ is not
commutative. It is all the more remarkable that it is {\em associative}, and even
better: the generalized associativity law (\ref{eqn:torsor1}) holds, which we write
in the new context in the form
\begin{equation}
\boxed{
(xo(upv)) = ((xou)pv)} \,  .
\label{eqn:torsor2}
\end{equation}
Exercise: construct the points $(xo(upv))$ and $((xou)pv) $, as shown in the following
illustration, and compare with the illustration in \ref{sec:ParAss2}:
the ``last vertex''  of the red trapezoid always coincides with the ``last vertex'' of
the blue one. 

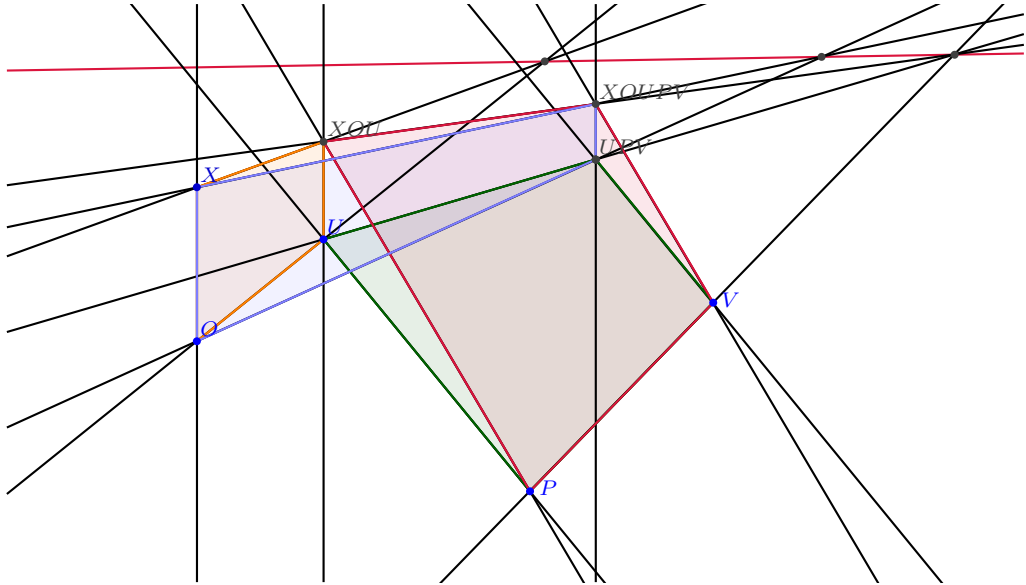
\begin{center}
\newrgbcolor{dcrutc}{0.86 0.08 0.24}
\newrgbcolor{qqqqcc}{0 0 0.8}
\newrgbcolor{ffxfqq}{1 0.5 0}
\newrgbcolor{qqwuqq}{0 0.39 0}
\newrgbcolor{xdxdff}{0.49 0.49 1}
\psset{xunit=0.5cm,yunit=0.5cm,algebraic=true,dotstyle=o,dotsize=3pt 0,linewidth=0.8pt,arrowsize=3pt 2,arrowinset=0.25}
\begin{pspicture*}(-4.7,-8.13)(22.05,7.2)
\pspolygon[linecolor=ffxfqq,fillcolor=ffxfqq,fillstyle=solid,opacity=0.1](0.3,-1.74)(0.3,2.34)(3.63,3.55)(3.63,0.96)
\pspolygon[linecolor=qqwuqq,fillcolor=qqwuqq,fillstyle=solid,opacity=0.1](9.06,-5.72)(13.88,-0.72)(10.79,3.08)(3.63,0.96)
\pspolygon[linecolor=dcrutc,fillcolor=dcrutc,fillstyle=solid,opacity=0.1](9.06,-5.72)(3.63,3.55)(10.79,4.55)(13.88,-0.72)
\pspolygon[linecolor=xdxdff,fillcolor=xdxdff,fillstyle=solid,opacity=0.1](0.3,2.34)(0.3,-1.74)(10.79,3.08)(10.79,4.55)
\psplot[linecolor=dcrutc]{-4.7}{22.05}{(--146.92--0.45*x)/26.63}
\psline(0.3,-8.13)(0.3,7.2)
\psplot{-4.7}{22.05}{(-6.6--2.7*x)/3.33}
\psplot{-4.7}{22.05}{(--20.37--3.34*x)/9.15}
\psline(3.63,-8.13)(3.63,7.2)
\psplot{-4.7}{22.05}{(-72.83--5*x)/4.82}
\psplot{-4.7}{22.05}{(--52.99-9.28*x)/5.43}
\psplot{-4.7}{22.05}{(--124.79-9.28*x)/5.43}
\psplot{-4.7}{22.05}{(--50.6--2.31*x)/16.59}
\psplot{-4.7}{22.05}{(-29.52--6.69*x)/-5.43}
\psplot{-4.7}{22.05}{(-88.84--6.69*x)/-5.43}
\psplot{-4.7}{22.05}{(-1.8--4.9*x)/16.59}
\psplot{-4.7}{22.05}{(-19.67--4.81*x)/10.49}
\psplot{-4.7}{22.05}{(--37.36--3.47*x)/16.43}
\psline(10.79,-8.13)(10.79,7.2)
\psline[linecolor=ffxfqq](0.3,-1.74)(0.3,2.34)
\psline[linecolor=ffxfqq](0.3,2.34)(3.63,3.55)
\psline[linecolor=ffxfqq](3.63,3.55)(3.63,0.96)
\psline[linecolor=ffxfqq](3.63,0.96)(0.3,-1.74)
\psline[linecolor=qqwuqq](9.06,-5.72)(13.88,-0.72)
\psline[linecolor=qqwuqq](13.88,-0.72)(10.79,3.08)
\psline[linecolor=qqwuqq](10.79,3.08)(3.63,0.96)
\psline[linecolor=qqwuqq](3.63,0.96)(9.06,-5.72)
\psline[linecolor=dcrutc](9.06,-5.72)(3.63,3.55)
\psline[linecolor=dcrutc](3.63,3.55)(10.79,4.55)
\psline[linecolor=dcrutc](10.79,4.55)(13.88,-0.72)
\psline[linecolor=dcrutc](13.88,-0.72)(9.06,-5.72)
\psline[linecolor=xdxdff](0.3,2.34)(0.3,-1.74)
\psline[linecolor=xdxdff](0.3,-1.74)(10.79,3.08)
\psline[linecolor=xdxdff](10.79,3.08)(10.79,4.55)
\psline[linecolor=xdxdff](10.79,4.55)(0.3,2.34)
\begin{scriptsize}
\psdots[dotstyle=*,linecolor=qqqqcc](0.3,2.34)
\rput[bl](0.39,2.47){\qqqqcc{$X$}}
\psdots[dotstyle=*,linecolor=blue](0.3,-1.74)
\rput[bl](0.39,-1.6){\blue{$O$}}
\psdots[dotstyle=*,linecolor=blue](3.63,0.96)
\rput[bl](3.72,1.1){\blue{$U$}}
\psdots[dotstyle=*,linecolor=darkgray](9.45,5.68)
\psdots[dotstyle=*,linecolor=darkgray](3.63,3.55)
\rput[bl](3.72,3.69){\darkgray{$XOU$}}
\psdots[dotstyle=*,linecolor=blue](9.06,-5.72)
\rput[bl](9.31,-5.81){\blue{$P$}}
\psdots[dotstyle=*,linecolor=blue](13.88,-0.72)
\rput[bl](14.1,-0.84){\blue{$V$}}
\psdots[dotstyle=*,linecolor=darkgray](20.23,5.86)
\psdots[dotstyle=*,linecolor=darkgray](10.79,4.55)
\rput[bl](10.88,4.68){\darkgray{$XOUPV$}}
\psdots[dotstyle=*,linecolor=darkgray](10.79,3.08)
\rput[bl](10.88,3.22){\darkgray{$UPV$}}
\psdots[dotstyle=*,linecolor=darkgray](10.79,3.08)
\psdots[dotstyle=*,linecolor=darkgray](10.79,3.08)
\psdots[dotstyle=*,linecolor=darkgray](16.73,5.8)
\end{scriptsize}
\end{pspicture*}
\end{center}

If you drag $o$ to $p$, you get the configuration of associativity of the
law $xz:= (xoz)$: $x(uv)=(xu)v$, that is, 
\begin{equation}
(xo(uov)) = ((xou)ov) \, .
\label{eqn:torsor3}
\end{equation}
In this case, there is a
 non-commutative analog of the ``cubic'' picture from \ref{sec:ParAss} -- we leave it to the reader
 to add the two ``diagonal trapezoids of the non-commutative cube'':
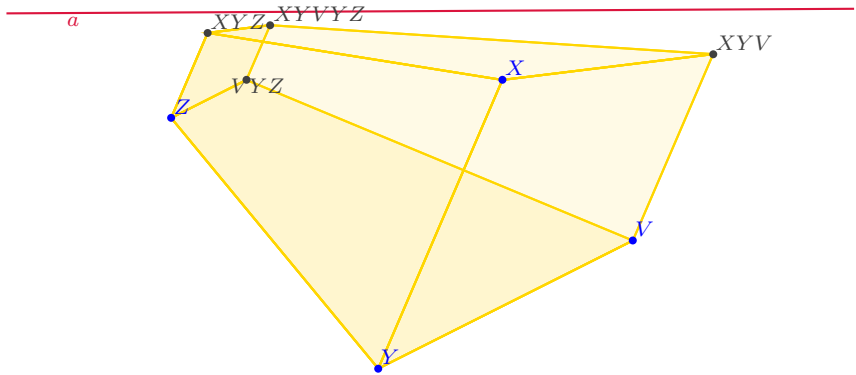
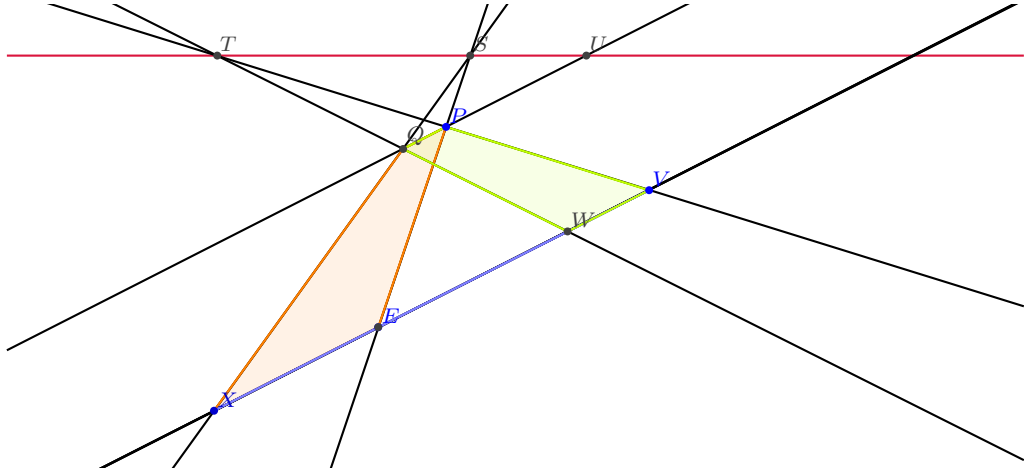
\begin{center}
\newrgbcolor{dcrutc}{0.86 0.08 0.24}
\newrgbcolor{ffdxqq}{1 0.84 0}
\psset{xunit=0.6cm,yunit=0.6cm,algebraic=true,dotstyle=o,dotsize=3pt 0,linewidth=0.8pt,arrowsize=3pt 2,arrowinset=0.25}
\begin{pspicture*}(-0.92,-5.2)(17.66,5.45)
\pspolygon[linecolor=ffdxqq,fillcolor=ffdxqq,fillstyle=solid,opacity=0.1](9.95,1.43)(7.23,-4.95)(2.69,0.59)(3.49,2.47)
\pspolygon[linecolor=ffdxqq,fillcolor=ffdxqq,fillstyle=solid,opacity=0.1](2.69,0.59)(3.49,2.47)(4.86,2.64)(4.34,1.43)
\pspolygon[linecolor=ffdxqq,fillcolor=ffdxqq,fillstyle=solid,opacity=0.1](7.23,-4.95)(12.81,-2.12)(4.34,1.43)(2.69,0.59)
\pspolygon[linecolor=ffdxqq,fillcolor=ffdxqq,fillstyle=solid,opacity=0.1](9.95,1.43)(14.57,2)(4.86,2.64)(3.49,2.47)
\pspolygon[linecolor=ffdxqq,fillcolor=ffdxqq,fillstyle=solid,opacity=0.1](9.95,1.43)(14.57,2)(12.81,-2.12)(7.23,-4.95)
\psplot[linecolor=dcrutc]{-0.92}{17.66}{(--53.18--0.1*x)/18.3}
\psline[linecolor=ffdxqq](9.95,1.43)(7.23,-4.95)
\psline[linecolor=ffdxqq](7.23,-4.95)(2.69,0.59)
\psline[linecolor=ffdxqq](2.69,0.59)(3.49,2.47)
\psline[linecolor=ffdxqq](3.49,2.47)(9.95,1.43)
\psline[linecolor=ffdxqq](2.69,0.59)(3.49,2.47)
\psline[linecolor=ffdxqq](3.49,2.47)(4.86,2.64)
\psline[linecolor=ffdxqq](4.86,2.64)(4.34,1.43)
\psline[linecolor=ffdxqq](4.34,1.43)(2.69,0.59)
\psline[linecolor=ffdxqq](7.23,-4.95)(12.81,-2.12)
\psline[linecolor=ffdxqq](12.81,-2.12)(4.34,1.43)
\psline[linecolor=ffdxqq](4.34,1.43)(2.69,0.59)
\psline[linecolor=ffdxqq](2.69,0.59)(7.23,-4.95)
\psline[linecolor=ffdxqq](9.95,1.43)(14.57,2)
\psline[linecolor=ffdxqq](14.57,2)(4.86,2.64)
\psline[linecolor=ffdxqq](4.86,2.64)(3.49,2.47)
\psline[linecolor=ffdxqq](3.49,2.47)(9.95,1.43)
\psline[linecolor=ffdxqq](9.95,1.43)(14.57,2)
\psline[linecolor=ffdxqq](14.57,2)(12.81,-2.12)
\psline[linecolor=ffdxqq](12.81,-2.12)(7.23,-4.95)
\psline[linecolor=ffdxqq](7.23,-4.95)(9.95,1.43)
\begin{scriptsize}
\rput[bl](0.41,2.62){\dcrutc{$a$}}
\psdots[dotstyle=*,linecolor=blue](9.95,1.43)
\rput[bl](10.01,1.53){\blue{$X$}}
\psdots[dotstyle=*,linecolor=blue](7.23,-4.95)
\rput[bl](7.29,-4.85){\blue{$Y$}}
\psdots[dotstyle=*,linecolor=blue](2.69,0.59)
\rput[bl](2.75,0.68){\blue{$Z$}}
\psdots[dotstyle=*,linecolor=darkgray](3.49,2.47)
\rput[bl](3.55,2.56){\darkgray{$XYZ$}}
\psdots[dotstyle=*,linecolor=blue](12.81,-2.12)
\rput[bl](12.87,-2.02){\blue{$V$}}
\psdots[dotstyle=*,linecolor=darkgray](14.57,2)
\rput[bl](14.63,2.09){\darkgray{$XYV$}}
\psdots[dotstyle=*,linecolor=darkgray](4.86,2.64)
\rput[bl](4.93,2.73){\darkgray{$XYVYZ$}}
\psdots[dotstyle=*,linecolor=darkgray](4.34,1.43)
\rput[bl](4,1.15){\darkgray{$VYZ$}}
\end{scriptsize}
\end{pspicture*}
\end{center}

\subsubsection{Unit}\label{subsec:Unit}

In your drawing of the trapezoid $(x,y,z,w)$,
drag the point $x$ to $y$ and observe that $w=(xyz)$ tends to $z$, and similarly,
if you drag $z$ to $y$, then $(xyz)$ tends to $x$.
Since our ``generic'' construction of $(xyz)$ is not defined for the cases $x=y$ or $x=z$,
we will now take as {\em definition} that
\begin{equation}
\boxed{
(xyy) := (yyx) := x } \, . \label{eqn:Idempot}
\end{equation}
This identity is called the {\em idempotent law}.
It means that, if we
fix some element $e \in G$, then  $e$ is a {\em unit element} for the binary product
$xz=(xez)$.

\subsubsection{Trapezoid geometry:  collinear case}\label{subsec:TrapezCollin}

As in the case of parallelograms, using the generalized associativity law (\ref{eqn:torsor2})
 together with
the idempotent law (\ref{eqn:Idempot}), we can
now better understand what is going on if $x,y,z$ are collinear.
Namely, assume that $x,y,z$ are on a line $\ell$ and
choose any point $p$ not on $\ell$.
Asssociativity and idempotency give us
$$
((xyp)pz) = (xy(ppz)) = (xyz),
$$
hence in the collinear case $(xyz)$ must be defined in the following way: 
{\em choose any point $p$ not on the line spanned by $x,y,z$;  then define
$(xyz):= ((xyp)pz)$} (note that the latter is defined by the generic contruction). 
Illustrate this --  
 reproduce  first  the illustration given above for the special case $u=p$ (and writing $e$
instead of $o$):

\begin{center}
\newrgbcolor{dcrutc}{0.86 0.08 0.24}
\newrgbcolor{qqqqcc}{0 0 0.8}
\newrgbcolor{ffxfqq}{1 0.5 0}
\newrgbcolor{xdxdff}{0.49 0.49 1}
\newrgbcolor{bfffqq}{0.75 1 0}
\psset{xunit=0.5cm,yunit=0.6cm,algebraic=true,dotstyle=o,dotsize=3pt 0,linewidth=0.8pt,arrowsize=3pt 2,arrowinset=0.25}
\begin{pspicture*}(-4.7,-8.13)(22.05,7.2)
\pspolygon[linecolor=ffxfqq,fillcolor=ffxfqq,fillstyle=solid,opacity=0.1](5.59,-3.34)(1.13,-5.38)(4.98,2.04)(6.4,2.7)
\pspolygon[linecolor=xdxdff,fillcolor=xdxdff,fillstyle=solid,opacity=0.1](1.13,-5.38)(5.59,-3.34)(10.5,0.42)(7.93,-0.76)
\pspolygon[linecolor=bfffqq,fillcolor=bfffqq,fillstyle=solid,opacity=0.1](6.4,2.7)(4.98,2.04)(7.93,-0.76)(10.5,0.42)
\psplot[linecolor=dcrutc]{-4.7}{22.05}{(--143.27--0.05*x)/25.96}
\psplot{-4.7}{22.05}{(-26.32--2.05*x)/4.46}
\psplot{-4.7}{22.05}{(-36.44--6.03*x)/0.81}
\psplot{-4.7}{22.05}{(-42.79--10.92*x)/5.65}
\psplot{-4.7}{22.05}{(-1.09--2.05*x)/4.46}
\psline[linecolor=ffxfqq](5.59,-3.34)(1.13,-5.38)
\psline[linecolor=ffxfqq](1.13,-5.38)(4.98,2.04)
\psline[linecolor=ffxfqq](4.98,2.04)(6.4,2.7)
\psline[linecolor=ffxfqq](6.4,2.7)(5.59,-3.34)
\psplot{-4.7}{22.05}{(-37.39--3.76*x)/4.91}
\psplot{-4.7}{22.05}{(-19.62--2.05*x)/4.46}
\psplot{-4.7}{22.05}{(--25.61-2.27*x)/4.1}
\psplot{-4.7}{22.05}{(-24.8--3.48*x)/-3.66}
\psplot{-4.7}{22.05}{(-98.85--10.93*x)/16.06}
\psline[linecolor=xdxdff](1.13,-5.38)(5.59,-3.34)
\psline[linecolor=xdxdff](5.59,-3.34)(10.5,0.42)
\psline[linecolor=xdxdff](10.5,0.42)(7.93,-0.76)
\psline[linecolor=xdxdff](7.93,-0.76)(1.13,-5.38)
\psline[linecolor=bfffqq](6.4,2.7)(4.98,2.04)
\psline[linecolor=bfffqq](4.98,2.04)(7.93,-0.76)
\psline[linecolor=bfffqq](7.93,-0.76)(10.5,0.42)
\psline[linecolor=bfffqq](10.5,0.42)(6.4,2.7)
\begin{scriptsize}
\psdots[dotstyle=*,linecolor=qqqqcc](1.13,-5.38)
\rput[bl](1.22,-5.25){\qqqqcc{$X$}}
\psdots[dotstyle=*,linecolor=blue](5.59,-3.34)
\rput[bl](5.68,-3.2){\blue{$E$}}
\psdots[dotstyle=*,linecolor=blue](6.4,2.7)
\rput[bl](6.49,2.83){\blue{$P$}}
\psdots[dotstyle=*,linecolor=darkgray](6.78,5.53)
\psdots[dotstyle=*,linecolor=darkgray](4.98,2.04)
\rput[bl](5.07,2.18){\darkgray{$XEP$}}
\psdots[dotstyle=*,linecolor=blue](10.5,0.42)
\rput[bl](10.59,0.56){\blue{$V$}}
\psdots[dotstyle=*,linecolor=darkgray](1.31,5.52)
\psdots[dotstyle=*,linecolor=darkgray](7.93,-0.76)
\rput[bl](8.02,-0.61){\darkgray{$XEPPV$}}
\psdots[dotstyle=*,linecolor=darkgray](17.19,5.55)
\psdots[dotstyle=*,linecolor=darkgray](12.59,5.54)
\end{scriptsize}
\end{pspicture*}
\end{center}

Next drag the point $e$ to the line $x \lor v$. The point $((xep)pv)$ then tends to a point
on the same line,
and observe that  this point does not depend on $p$: it remains fixed if you move $p$ around. 

\begin{center}
\newrgbcolor{dcrutc}{0.86 0.08 0.24}
\newrgbcolor{qqqqcc}{0 0 0.8}
\newrgbcolor{ffxfqq}{1 0.5 0}
\newrgbcolor{xdxdff}{0.49 0.49 1}
\newrgbcolor{bfffqq}{0.75 1 0}
\psset{xunit=0.5cm,yunit=0.4cm,algebraic=true,dotstyle=o,dotsize=3pt 0,linewidth=0.8pt,arrowsize=3pt 2,arrowinset=0.25}
\begin{pspicture*}(-4.7,-8.13)(22.05,7.2)
\pspolygon[linecolor=ffxfqq,fillcolor=ffxfqq,fillstyle=solid,opacity=0.1](5.07,-3.49)(0.75,-6.26)(5.72,2.42)(6.85,3.15)
\pspolygon[linecolor=xdxdff,fillcolor=xdxdff,fillstyle=solid,opacity=0.1](0.75,-6.26)(5.07,-3.49)(12.19,1.05)(10.05,-0.32)
\pspolygon[linecolor=bfffqq,fillcolor=bfffqq,fillstyle=solid,opacity=0.1](6.85,3.15)(5.72,2.42)(10.05,-0.32)(12.19,1.05)
\psplot[linecolor=dcrutc]{-4.7}{22.05}{(--146.41-0*x)/26.56}
\psplot{-4.7}{22.05}{(-29.15--2.77*x)/4.32}
\psplot{-4.7}{22.05}{(-39.91--6.64*x)/1.78}
\psplot{-4.7}{22.05}{(-51.01--11.77*x)/6.73}
\psplot{-4.7}{22.05}{(-5.37--2.77*x)/4.32}
\psline[linecolor=ffxfqq](5.07,-3.49)(0.75,-6.26)
\psline[linecolor=ffxfqq](0.75,-6.26)(5.72,2.42)
\psline[linecolor=ffxfqq](5.72,2.42)(6.85,3.15)
\psline[linecolor=ffxfqq](6.85,3.15)(5.07,-3.49)
\psplot{-4.7}{22.05}{(-47.92--4.55*x)/7.11}
\psplot{-4.7}{22.05}{(-29.19--2.77*x)/4.32}
\psplot{-4.7}{22.05}{(--31.14-2.09*x)/5.34}
\psplot{-4.7}{22.05}{(-29.51--3.09*x)/-4.89}
\psplot{-4.7}{22.05}{(-124.12--11.77*x)/18.41}
\psline[linecolor=xdxdff](0.75,-6.26)(5.07,-3.49)
\psline[linecolor=xdxdff](5.07,-3.49)(12.19,1.05)
\psline[linecolor=xdxdff](12.19,1.05)(10.05,-0.32)
\psline[linecolor=xdxdff](10.05,-0.32)(0.75,-6.26)
\psline[linecolor=bfffqq](6.85,3.15)(5.72,2.42)
\psline[linecolor=bfffqq](5.72,2.42)(10.05,-0.32)
\psline[linecolor=bfffqq](10.05,-0.32)(12.19,1.05)
\psline[linecolor=bfffqq](12.19,1.05)(6.85,3.15)
\begin{scriptsize}
\psdots[dotstyle=*,linecolor=qqqqcc](0.75,-6.26)
\rput[bl](0.84,-6.13){\qqqqcc{$X$}}
\psdots[dotstyle=*,linecolor=blue](5.07,-3.49)
\rput[bl](5.16,-3.36){\blue{$E$}}
\psdots[dotstyle=*,linecolor=blue](6.85,3.15)
\rput[bl](6.94,3.28){\blue{$P$}}
\psdots[dotstyle=*,linecolor=darkgray](7.49,5.51)
\rput[bl](7.57,5.65){\darkgray{$S$}}
\psdots[dotstyle=*,linecolor=darkgray](5.72,2.42)
\rput[bl](5.82,2.56){\darkgray{$Q$}}
\psdots[dotstyle=*,linecolor=blue](12.19,1.05)
\rput[bl](12.28,1.19){\blue{$V$}}
\psdots[dotstyle=*,linecolor=darkgray](0.83,5.51)
\rput[bl](0.91,5.65){\darkgray{$T$}}
\psdots[dotstyle=*,linecolor=darkgray](10.05,-0.32)
\rput[bl](10.14,-0.18){\darkgray{$W$}}
\psdots[dotstyle=*,linecolor=darkgray](10.54,5.51)
\rput[bl](10.63,5.65){\darkgray{$U$}}
\psdots[dotstyle=*,linecolor=darkgray](5.07,-3.49)
\end{scriptsize}
\end{pspicture*}
\end{center}

Summing up, for three points $x,e,v$ on a line, we have defined a fourth point $w=(xev)$
on the same line.
Such constructions are very important in geometry, since they are related to the
usual {\em product} and {\em sum} on the real line $\R$ (or over any other field). 
We will give later a computational proof (see \ref{subsec:ComputationTrapezoid}), 
but it is probably quite
insctructive to indicate here a geometric argument explaining this, as follows. 

\ssk
Using the preceding drawing, show:
if we identify the line $\ell = x \lor v$ with $\R$ by choosing the intersection point
of $\ell$ with the horizontal line as origin $0 \in \R$ (so $x,e,v$ correspond to certain
non-zero real numbers), then the point $w = (xev)$ is identified with
the element $x e\inv v$ of $\R$.
(Hints: let $q:=(xep)$;
 use the intercept (or Thales) theorem twice:
using $s$ as center, show that $x/e = \vert q - u \vert / \vert p - u \vert$, and
using $t$ as center, that $w/v =  \vert q - u \vert / \vert p - u \vert$,
whence $w/v = x/e$.)
In particular, if we identity $e$ with the element $1 \in \R$, this gives us
$w = xv$. 

\smallskip
There is, however, one special position of a different kind, namely
 when $x,e,v$ all lie on a line that is parallel to $a$:
 
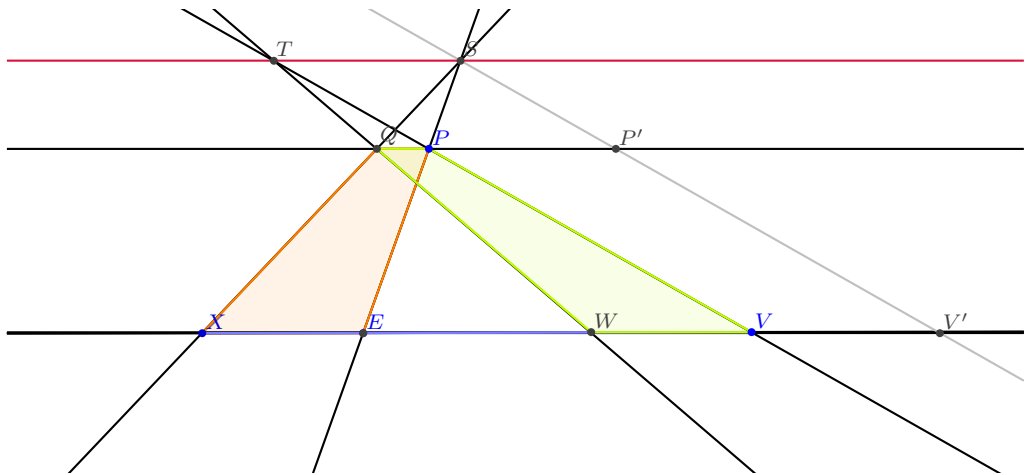
\begin{center}
\newrgbcolor{dcrutc}{0.86 0.08 0.24}
\newrgbcolor{qqqqcc}{0 0 0.8}
\newrgbcolor{ffxfqq}{1 0.5 0}
\newrgbcolor{xdxdff}{0.49 0.49 1}
\newrgbcolor{bfffqq}{0.75 1 0}
\psset{xunit=0.5cm,yunit=0.4cm,algebraic=true,dotstyle=o,dotsize=3pt 0,linewidth=0.8pt,arrowsize=3pt 2,arrowinset=0.25}
\begin{pspicture*}(-4.7,-8.13)(22.05,7.2)
\pspolygon[linecolor=ffxfqq,fillcolor=ffxfqq,fillstyle=solid,opacity=0.1](4.67,-3.52)(0.44,-3.52)(5.03,2.59)(6.4,2.59)
\pspolygon[linecolor=xdxdff,fillcolor=xdxdff,fillstyle=solid,opacity=0.1](0.44,-3.52)(4.67,-3.52)(14.89,-3.49)(10.67,-3.49)
\pspolygon[linecolor=bfffqq,fillcolor=bfffqq,fillstyle=solid,opacity=0.1](6.4,2.59)(5.03,2.59)(10.67,-3.49)(14.89,-3.49)
\psplot[linecolor=dcrutc]{-4.7}{22.05}{(--146.41-0*x)/26.56}
\psplot{-4.7}{22.05}{(-14.88-0*x)/4.23}
\psplot{-4.7}{22.05}{(-34.57--6.1*x)/1.73}
\psplot{-4.7}{22.05}{(-27.83--9.03*x)/6.8}
\psplot{-4.7}{22.05}{(--10.94-0*x)/4.23}
\psline[linecolor=ffxfqq](4.67,-3.52)(0.44,-3.52)
\psline[linecolor=ffxfqq](0.44,-3.52)(5.03,2.59)
\psline[linecolor=ffxfqq](5.03,2.59)(6.4,2.59)
\psline[linecolor=ffxfqq](6.4,2.59)(4.67,-3.52)
\psplot{-4.7}{22.05}{(-36.04--0.02*x)/10.22}
\psplot{-4.7}{22.05}{(-14.78-0*x)/4.23}
\psplot{-4.7}{22.05}{(--60.85-6.08*x)/8.49}
\psplot{-4.7}{22.05}{(-21.74--2.93*x)/-2.71}
\psplot{-4.7}{22.05}{(-14426.84--9.03*x)/4102.56}
\psline[linecolor=xdxdff](0.44,-3.52)(4.67,-3.52)
\psline[linecolor=xdxdff](4.67,-3.52)(14.89,-3.49)
\psline[linecolor=xdxdff](14.89,-3.49)(10.67,-3.49)
\psline[linecolor=xdxdff](10.67,-3.49)(0.44,-3.52)
\psline[linecolor=bfffqq](6.4,2.59)(5.03,2.59)
\psline[linecolor=bfffqq](5.03,2.59)(10.67,-3.49)
\psline[linecolor=bfffqq](10.67,-3.49)(14.89,-3.49)
\psline[linecolor=bfffqq](14.89,-3.49)(6.4,2.59)
\psplot[linecolor=lightgray]{-4.7}{22.05}{(--14705.23-985*x)/1375.35}
\begin{scriptsize}
\psdots[dotstyle=*,linecolor=qqqqcc](0.44,-3.52)
\rput[bl](0.53,-3.38){\qqqqcc{$X$}}
\psdots[dotstyle=*,linecolor=blue](4.67,-3.52)
\rput[bl](4.76,-3.38){\blue{$E$}}
\psdots[dotstyle=*,linecolor=blue](6.4,2.59)
\rput[bl](6.49,2.72){\blue{$P$}}
\psdots[dotstyle=*,linecolor=darkgray](7.23,5.51)
\rput[bl](7.32,5.65){\darkgray{$S$}}
\psdots[dotstyle=*,linecolor=darkgray](5.03,2.59)
\rput[bl](5.12,2.72){\darkgray{$Q$}}
\psdots[dotstyle=*,linecolor=blue](14.89,-3.49)
\rput[bl](14.98,-3.36){\blue{$V$}}
\psdots[dotstyle=*,linecolor=darkgray](2.32,5.51)
\rput[bl](2.39,5.65){\darkgray{$T$}}
\psdots[dotstyle=*,linecolor=darkgray](10.67,-3.49)
\rput[bl](10.75,-3.36){\darkgray{$W$}}
\psdots[dotstyle=*,linecolor=darkgray](4.67,-3.52)
\psdots[dotstyle=*,linecolor=darkgray](11.32,2.59)
\rput[bl](11.4,2.72){\darkgray{$P'$}}
\psdots[dotstyle=*,linecolor=darkgray](19.84,-3.52)
\rput[bl](19.93,-3.38){\darkgray{$V'$}}
\end{scriptsize}
\end{pspicture*}
\end{center}
Show: if, in this situation,
 we identify the line $\ell$ with $\R$  by choosing arbitrarily two points on $\ell$
 and identifying them with $0$, resp.\ $1$, then 
the point $w=(xev)$ corresponds to the sum $x-e+v$ in $\R$. 
(Hints: 
for the purpose of the proof, draw also the line parallel to $p \lor v$ through $s$
and note that $\vert t p \vert = \vert sp' \vert$ and
$\vert t v \vert = \vert s v' \vert$, then
 use again twice the intercept theorem, with center $T$, resp.\ $S$ and conclude
 that $\vert x e \vert = \vert w v \vert$, whence 
 $x-e = w-v$ in $\R$.)
In particular, if we choose $e$ as origin $0$ in $\R$, we will have
$w = x + v$.

\subsubsection{Group law}\label{subsec:Geo-Grouplaw}

Let us put  things together:  convince yourself that, for any choice of origin $e \in G$,
the set $G = P \setminus a$, is a group with product $xz:=(xez)$.
Indeed,
associativity and property of the unit element have already been seen;
and finally, every element $x$ has an inverse element: 
check this first by direct inspection of your drawing of the trapezoid $(x,e,z,w)$:
that is, check by trying that for each $x$ there is some point $z$ such
that $(xez)=e=(zex)$. 

\ssk
By using the generalized associativity law, show that $u= (exe)$ is the inverse of $x$.
(Just check that $(xeu)=e=(uex)$.) Below we will inspect the
 geometric construction of the inverse.

\ssk
Using associativity and idempotency, compute also $((xe(eye))ez)$ and conclude that
\begin{equation}
\boxed{
(xyz) = x y\inv z }
\end{equation}
where the right hand side is computed in the group $(G,e)$.
In other words, the ternary law $(xyz)$ can be recovered from the binary product
in the group with neutral element $e$.

\subsubsection{One-parameter subgroups}\label{subsec:One-parameter}

We wish to study the structure of the group $(G,e)$ in more detail.
First of all, note that $G$ has {\em two connected components}: the ``upper half plane''
and the ``lower half plane''; show that the half plane containing $e$ (in our images: the lower one)
is a subgroup, and the other is not (how are points from this other half plane multiplied?). 

\ssk
Show that
every line $\ell$ through $e$ is a subgroup of $G$ (it suffices to observe that, if
$x,e,z$ are on such a line, then $(xez)$ is on the same line, and so are their inverses).
Let us call a subgroup of this type a {\em one-parameter subgroup}.
To be more precise: if $\ell$ is not parallel to $a$, then the point 
$\ell \land a$ has to be excluded, so the subgroup is really $H:=\ell \setminus a$
(as noted above, as a group, this is the multiplicative group $\R^\times$).
If $\ell$ is parallel to $a$, then $H:=\ell$ is a subgroup (isomorphic to $(\R,+)$).

\ssk
Fix an origin $o$ and the line $a$, and let $x,o,h$ not be collinear points.
We have seen that the line $H$ spanned by $o$ and $h$ is a subgroup of $G$.
Draw the {\em left} and {\em right cosets} of $H$ through $x$,
$$
xH := \{ xz \mid z \in H \}, \qquad
Hx:= \{ zx \mid z \in H \}.
$$
First of all, convince yourself that $xH$ and $Hx$ are lines passing through $x$.
In your drawing, they may look as follows (the thick grey and blue lines, where
$z \in \ell$ is a typical element of the group $H$):

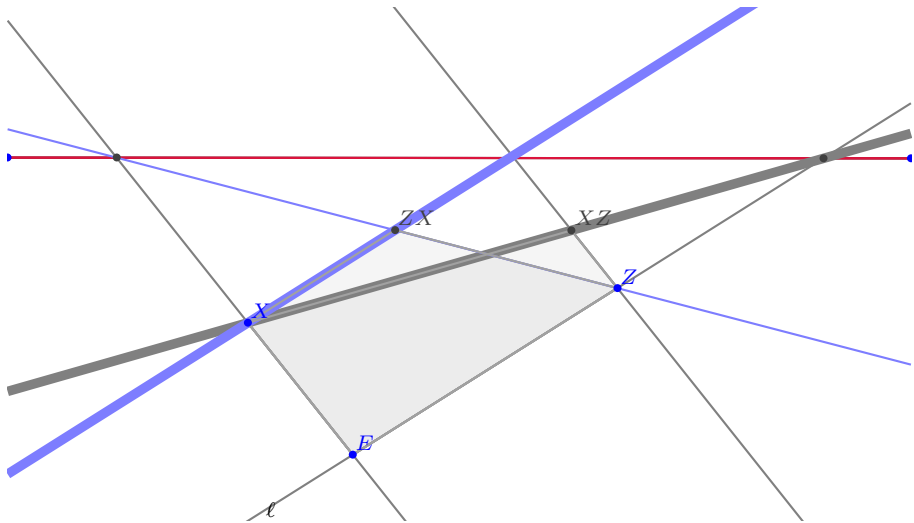
\begin{center} 
\newrgbcolor{dcrutc}{0.86 0.08 0.24}
\newrgbcolor{xdxdff}{0.49 0.49 1}
\newrgbcolor{aqaqaq}{0.63 0.63 0.63}
\psset{xunit=0.5cm,yunit=0.5cm,algebraic=true,dotstyle=o,dotsize=3pt 0,linewidth=0.8pt,arrowsize=3pt 2,arrowinset=0.25}
\begin{pspicture*}(-4.3,-7.32)(19.46,6.3)
\pspolygon[linecolor=aqaqaq,fillcolor=aqaqaq,fillstyle=solid,opacity=0.1](10.53,0.39)(11.74,-1.14)(4.78,-5.56)(2.02,-2.06)
\pspolygon[linecolor=aqaqaq,fillcolor=aqaqaq,fillstyle=solid,opacity=0.1](5.89,0.4)(11.74,-1.14)(4.78,-5.56)(2.02,-2.06)
\psplot{-4.3}{19.46}{(--55.04-0.02*x)/23.76}
\psplot[linecolor=dcrutc]{-4.3}{19.46}{(-55.04--0.02*x)/-23.76}
\psplot[linecolor=dcrutc]{-4.3}{19.46}{(-55.04--0.02*x)/-23.76}
\psplot[linecolor=gray]{-4.3}{19.46}{(--1.38-3.5*x)/2.76}
\psplot[linecolor=gray]{-4.3}{19.46}{(--37.94-3.5*x)/2.76}
\psplot[linecolor=gray]{-4.3}{19.46}{(-59.83--4.42*x)/6.96}
\psplot[linewidth=3.6pt,linecolor=gray]{-4.3}{19.46}{(-40--4.36*x)/15.14}
\psplot[linewidth=3.6pt,linecolor=xdxdff]{-4.3}{19.46}{(-23.27--4.42*x)/6.96}
\psplot[linecolor=xdxdff]{-4.3}{19.46}{(--25.58-3.46*x)/13.17}
\psline[linecolor=aqaqaq](10.53,0.39)(11.74,-1.14)
\psline[linecolor=aqaqaq](11.74,-1.14)(4.78,-5.56)
\psline[linecolor=aqaqaq](4.78,-5.56)(2.02,-2.06)
\psline[linecolor=aqaqaq](2.02,-2.06)(10.53,0.39)
\psline[linecolor=aqaqaq](5.89,0.4)(11.74,-1.14)
\psline[linecolor=aqaqaq](11.74,-1.14)(4.78,-5.56)
\psline[linecolor=aqaqaq](4.78,-5.56)(2.02,-2.06)
\psline[linecolor=aqaqaq](2.02,-2.06)(5.89,0.4)
\begin{scriptsize}
\psdots[dotstyle=*,linecolor=blue](-4.3,2.32)
\psdots[dotstyle=*,linecolor=blue](19.46,2.3)
\psdots[dotstyle=*,linecolor=blue](4.78,-5.56)
\rput[bl](4.86,-5.44){\blue{$E$}}
\psdots[dotstyle=*,linecolor=blue](2.02,-2.06)
\rput[bl](2.1,-1.94){\blue{$X$}}
\psdots[dotstyle=*,linecolor=blue](11.74,-1.14)
\rput[bl](11.82,-1.02){\blue{$Z$}}
\rput[bl](2.5,-7.22){\black{$\ell$}}
\psdots[dotstyle=*,linecolor=darkgray](17.16,2.3)
\psdots[dotstyle=*,linecolor=darkgray](10.53,0.39)
\rput[bl](10.62,0.52){\darkgray{$XZ$}}
\psdots[dotstyle=*,linecolor=darkgray](-1.43,2.32)
\psdots[dotstyle=*,linecolor=darkgray](5.89,0.4)
\rput[bl](5.96,0.52){\darkgray{$ZX$}}
\end{scriptsize}
\end{pspicture*}
\end{center}

Note that, in the drawing, $xH$ and $Hx$ are different lines.
However, if you move the point $z$ around (while $e$ remains fixed), then
there is one remarkable position for which we do have $xH = Hx$:
show that this  happens when the line $\ell$ is parallel to $a$
(this does not mean that the points $xz$ and $zx$ coincide, as you see here):

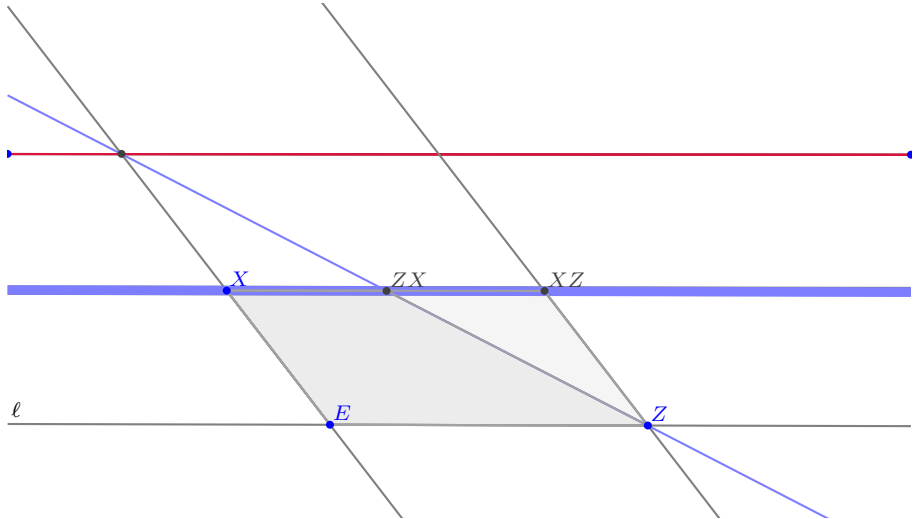
\begin{center} 
\newrgbcolor{dcrutc}{0.86 0.08 0.24}
\newrgbcolor{xdxdff}{0.49 0.49 1}
\newrgbcolor{aqaqaq}{0.63 0.63 0.63}
\psset{xunit=0.5cm,yunit=0.5cm,algebraic=true,dotstyle=o,dotsize=3pt 0,linewidth=0.8pt,arrowsize=3pt 2,arrowinset=0.25}
\begin{pspicture*}(-4.3,-7.32)(19.46,6.3)
\pspolygon[linecolor=aqaqaq,fillcolor=aqaqaq,fillstyle=solid,opacity=0.1](9.82,-1.31)(12.54,-4.88)(4.18,-4.86)(1.46,-1.3)
\pspolygon[linecolor=aqaqaq,fillcolor=aqaqaq,fillstyle=solid,opacity=0.1](5.67,-1.31)(12.54,-4.88)(4.18,-4.86)(1.46,-1.3)
\psplot{-4.3}{19.46}{(--55.04-0.02*x)/23.76}
\psplot[linecolor=dcrutc]{-4.3}{19.46}{(-55.04--0.02*x)/-23.76}
\psplot[linecolor=dcrutc]{-4.3}{19.46}{(-55.04--0.02*x)/-23.76}
\psplot[linecolor=gray]{-4.3}{19.46}{(--1.66-3.56*x)/2.72}
\psplot[linecolor=gray]{-4.3}{19.46}{(--31.37-3.56*x)/2.72}
\psplot[linecolor=gray]{-4.3}{19.46}{(-40.55-0.02*x)/8.36}
\psplot[linewidth=3.6pt,linecolor=gray]{-4.3}{19.46}{(--5999.15--7.51*x)/-4623.16}
\psplot[linewidth=3.6pt,linecolor=xdxdff]{-4.3}{19.46}{(-10.84-0.02*x)/8.36}
\psplot[linecolor=xdxdff]{-4.3}{19.46}{(--22.7-7.2*x)/13.84}
\psline[linecolor=aqaqaq](9.82,-1.31)(12.54,-4.88)
\psline[linecolor=aqaqaq](12.54,-4.88)(4.18,-4.86)
\psline[linecolor=aqaqaq](4.18,-4.86)(1.46,-1.3)
\psline[linecolor=aqaqaq](1.46,-1.3)(9.82,-1.31)
\psline[linecolor=aqaqaq](5.67,-1.31)(12.54,-4.88)
\psline[linecolor=aqaqaq](12.54,-4.88)(4.18,-4.86)
\psline[linecolor=aqaqaq](4.18,-4.86)(1.46,-1.3)
\psline[linecolor=aqaqaq](1.46,-1.3)(5.67,-1.31)
\begin{scriptsize}
\psdots[dotstyle=*,linecolor=blue](-4.3,2.32)
\psdots[dotstyle=*,linecolor=blue](19.46,2.3)
\psdots[dotstyle=*,linecolor=blue](4.18,-4.86)
\rput[bl](4.26,-4.74){\blue{$E$}}
\psdots[dotstyle=*,linecolor=blue](1.46,-1.3)
\rput[bl](1.54,-1.18){\blue{$X$}}
\psdots[dotstyle=*,linecolor=blue](12.54,-4.88)
\rput[bl](12.62,-4.76){\blue{$Z$}}
\rput[bl](-4.2,-4.66){\black{$\ell$}}
\psdots[dotstyle=*,linecolor=darkgray](9.82,-1.31)
\rput[bl](9.9,-1.2){\darkgray{$XZ$}}
\psdots[dotstyle=*,linecolor=darkgray](-1.3,2.32)
\psdots[dotstyle=*,linecolor=darkgray](5.67,-1.31)
\rput[bl](5.76,-1.2){\darkgray{$ZX$}}
\end{scriptsize}
\end{pspicture*}
\end{center}

Recall from group theory: the subgroup $H$ is called a {\em normal subgroup} if
$xH=Hx$ for all $x \in G$. We have seen: $H$ is normal if, and only if, $H$ is parallel
to $a$. 
Thus there is precisely one one-parameter subgroup which is a normal subgroup; all others are not.
Call this normal subgroup $N$,  and fix some other one-parameter subgroup, say $K$.
From group theory it is known that the quotient $G/N$ is a group, and in fact this group
is isomorphic to the complementary subgroup $K$. Illustrate this as follows:
choose two lines parallel to $N$ (and different from $a$), and confine the variables $x$ and
$z$ to these subgroups (that is, draw first the lines and then use the command ``point on object''). 
Observe that the products $x z$ then always lie on a third parallel
to $N$, and hence we can ``multiply'' these parallels in a well-defined way. Since each
such parallel intersects $K$ in a unique point, the group $K$ is indeed a group of 
representatives of these parallels, and the representatives $x',z'$ are multiplied exactly
as the parallels through them:

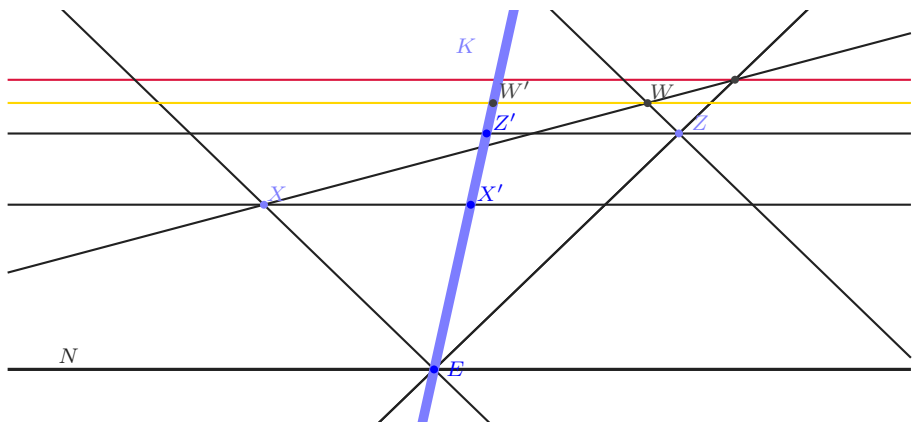
\begin{center}
\newrgbcolor{dcrutc}{0.86 0.08 0.24}
\newrgbcolor{sqsqsq}{0.13 0.13 0.13}
\newrgbcolor{xdxdff}{0.49 0.49 1}
\newrgbcolor{ffdxqq}{1 0.84 0}
\psset{xunit=0.5cm,yunit=0.4cm,algebraic=true,dotstyle=o,dotsize=3pt 0,linewidth=0.8pt,arrowsize=3pt 2,arrowinset=0.25}
\begin{pspicture*}(-4.3,-7.32)(19.46,6.3)
\psplot[linecolor=dcrutc]{-4.3}{19.46}{(--72-0*x)/18}
\psplot[linecolor=sqsqsq]{-4.3}{19.46}{(-2.52-0*x)/18}
\psplot[linewidth=1.2pt,linecolor=sqsqsq]{-4.3}{19.46}{(-100.8-0*x)/18}
\psplot[linecolor=sqsqsq]{-4.3}{19.46}{(--39.96-0*x)/18}
\psplot[linecolor=sqsqsq]{-4.3}{19.46}{(--12.7-5.46*x)/4.48}
\psplot[linecolor=sqsqsq]{-4.3}{19.46}{(-90.18--7.82*x)/6.44}
\psplot[linecolor=sqsqsq]{-4.3}{19.46}{(-11.84--4.14*x)/12.39}
\psplot[linecolor=sqsqsq]{-4.3}{19.46}{(--82.89-5.46*x)/4.48}
\psplot[linecolor=sqsqsq]{-4.3}{19.46}{(-90.18--7.82*x)/6.44}
\psplot[linecolor=ffdxqq]{-4.3}{19.46}{(--5013.85-0*x)/1551.13}
\psplot[linewidth=3.6pt,linecolor=xdxdff]{-4.3}{19.46}{(-85.38--10.8*x)/1.9}
\begin{scriptsize}
\psdots[dotstyle=*,linecolor=blue](6.92,-5.6)
\rput[bl](7.24,-5.82){\blue{$E$}}
\psdots[dotstyle=*,linecolor=blue](7.88,-0.14)
\rput[bl](8.04,-0.02){\blue{$X'$}}
\rput[bl](-2.96,-5.4){\sqsqsq{$N$}}
\psdots[dotstyle=*,linecolor=blue](8.3,2.22)
\rput[bl](8.46,2.34){\blue{$Z'$}}
\psdots[dotstyle=*,linecolor=xdxdff](2.44,-0.14)
\rput[bl](2.52,-0.02){\xdxdff{$X$}}
\psdots[dotstyle=*,linecolor=xdxdff](13.36,2.22)
\rput[bl](13.7,2.32){\xdxdff{$Z$}}
\psdots[dotstyle=*,linecolor=darkgray](14.83,4)
\psdots[dotstyle=*,linecolor=darkgray](12.53,3.23)
\rput[bl](12.6,3.36){\darkgray{$W$}}
\rput[bl](7.48,4.88){\xdxdff{$K$}}
\psdots[dotstyle=*,linecolor=darkgray](8.47,3.23)
\rput[bl](8.62,3.4){\darkgray{$W'$}}
\end{scriptsize}
\end{pspicture*}
\end{center}

Finally,
show that, for every element $g \in G$, there is exactly one $z \in K$ and exactly one $x \in N$
such that $g= xz$ (we write then $G=KN$),
respectively such that $g=z' x'$ (so $G=NK$; show that then $z=z'$ but in general
$x \not= x'$). Hint: given $g=w$, do the following construction:

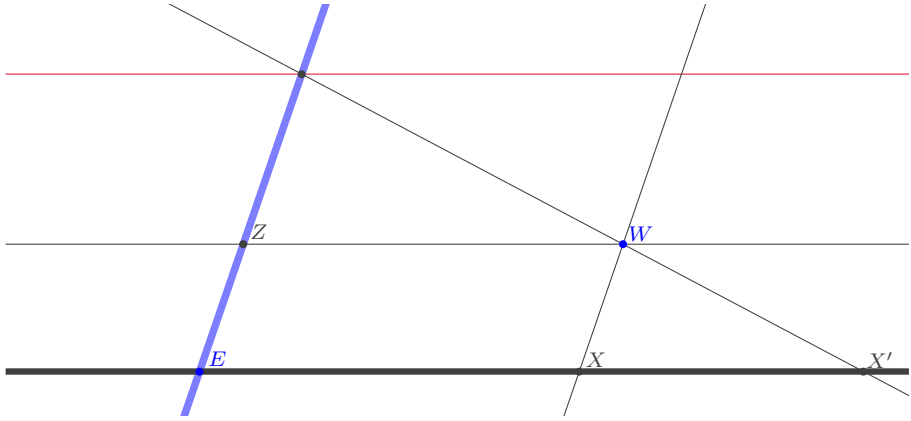
\begin{center}
\newrgbcolor{dcrutc}{0.86 0.08 0.24}
\newrgbcolor{xdxdff}{0.49 0.49 1}
\newrgbcolor{uququq}{0.25 0.25 0.25}
\psset{xunit=0.5cm,yunit=0.4cm,algebraic=true,dotstyle=o,dotsize=3pt 0,linewidth=0.8pt,arrowsize=3pt 2,arrowinset=0.25}
\begin{pspicture*}(-4.3,-7.32)(19.46,6.3)
\psplot[linewidth=0.4pt,linecolor=dcrutc]{-4.3}{19.46}{(--84-0*x)/21}
\psplot[linewidth=2.8pt,linecolor=xdxdff]{-4.3}{19.46}{(-26.38--11*x)/3}
\psplot[linewidth=2.4pt,linecolor=uququq]{-4.3}{19.46}{(-123.06-0*x)/21}
\psplot[linewidth=0.4pt,linecolor=uququq]{-4.3}{19.46}{(-136.26--11*x)/3}
\psplot[linewidth=0.4pt,linecolor=uququq]{-4.3}{19.46}{(-34.44-0*x)/21}
\psplot[linewidth=0.4pt,linecolor=uququq]{-4.3}{19.46}{(--53.48-5.64*x)/8.45}
\begin{scriptsize}
\psdots[dotstyle=*,linecolor=blue](0.8,-5.86)
\rput[bl](1.04,-5.68){\blue{$E$}}
\psdots[dotstyle=*,linecolor=blue](11.94,-1.64)
\rput[bl](12.08,-1.52){\blue{$W$}}
\psdots[dotstyle=*,linecolor=darkgray](1.95,-1.64)
\rput[bl](2.14,-1.48){\darkgray{$Z$}}
\psdots[dotstyle=*,linecolor=darkgray](10.79,-5.86)
\rput[bl](10.96,-5.7){\darkgray{$X$}}
\psdots[dotstyle=*,linecolor=darkgray](3.49,4)
\psdots[dotstyle=*,linecolor=darkgray](18.26,-5.86)
\rput[bl](18.34,-5.74){\darkgray{$X'$}}
\end{scriptsize}
\end{pspicture*}
\end{center}

Remark. 
In algebraic language, $G$ is a {\em semidirect product} of $N$ and $K$, see chapter \ref{sec:Compute}.
Thus $G$ is isomorphic to a semidirect product of $\R^\times$ with $\R$.


\subsubsection{Boundary values and affine group of the line}\label{subsec:Geo-Boundary}

So far we have excluded the line $a$ from our space: we have assumed that $x,y,z$ are
not on $a$. However, the construction continues to make sense if one of the points
$x$ or $z$ (but not $y$) belongs to $a$:
fix a point $e \in G$ as origin and consider the product $(xez)$ defined by formula
(\ref{eqn:Generic}).
Using your drawing from section \ref{subsec:GenericConstruction} , observe what is going to happen:

\ssk
-- if $z \in a$, show that always $(xez)=z$ (explain why!),

-- if $x \in a$, show that $(xez)$ always  lies on $a$, and, more precisely,  show that we get
\begin{equation}
(xez) = (((x \lor e) \land i) \lor z) \land a \, ,
\end{equation}
and this is the intersection of $a$ with the parallel of $x \lor e$ through $z$:


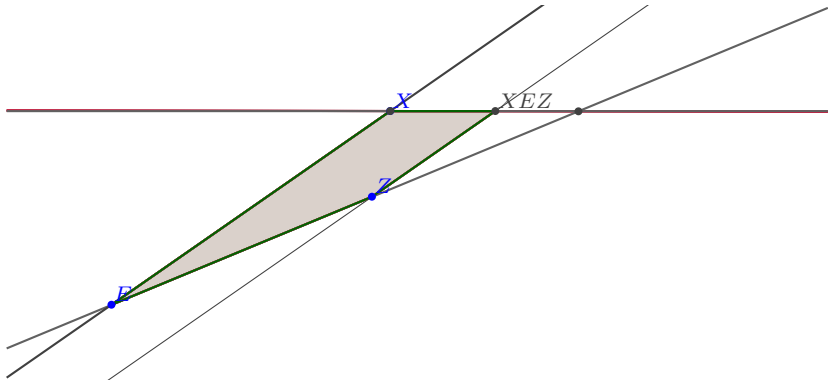
\begin{center}
\newrgbcolor{dcrutc}{0.86 0.08 0.24}
\newrgbcolor{wqwqwq}{0.38 0.38 0.38}
\newrgbcolor{uququq}{0.25 0.25 0.25}
\newrgbcolor{aqaqaq}{0.63 0.63 0.63}
\newrgbcolor{qqwuqq}{0 0.39 0}
\psset{xunit=0.5cm,yunit=0.4cm,algebraic=true,dotstyle=o,dotsize=3pt 0,linewidth=0.8pt,arrowsize=3pt 2,arrowinset=0.25}
\begin{pspicture*}(-3.72,-6.97)(17.88,5.41)
\pspolygon[linecolor=aqaqaq,fillcolor=aqaqaq,fillstyle=solid,opacity=0.1](9.13,1.92)(6.37,1.92)(-0.96,-4.5)(5.89,-0.92)
\pspolygon[linecolor=dcrutc,fillcolor=dcrutc,fillstyle=solid,opacity=0.1](5.89,-0.92)(9.13,1.92)(6.37,1.92)(-0.96,-4.5)
\pspolygon[linecolor=qqwuqq,fillcolor=qqwuqq,fillstyle=solid,opacity=0.1](5.89,-0.92)(9.13,1.92)(6.37,1.92)(-0.96,-4.5)
\psplot{-3.72}{17.88}{(--47.16-0.04*x)/24.4}
\psplot[linecolor=dcrutc]{-3.72}{17.88}{(-47.16--0.04*x)/-24.4}
\psplot[linecolor=dcrutc]{-3.72}{17.88}{(-47.16--0.04*x)/-24.4}
\psplot[linecolor=wqwqwq]{-3.72}{17.88}{(--27.42-3.58*x)/-6.85}
\psplot[linecolor=uququq]{-3.72}{17.88}{(-26.82--6.42*x)/7.33}
\psplot[linewidth=0.4pt,linecolor=uququq]{-3.72}{17.88}{(-44.57--6.42*x)/7.33}
\psplot[linecolor=wqwqwq]{-3.72}{17.88}{(-9.5-0*x)/-4.95}
\psline[linecolor=aqaqaq](9.13,1.92)(6.37,1.92)
\psline[linecolor=aqaqaq](6.37,1.92)(-0.96,-4.5)
\psline[linecolor=aqaqaq](-0.96,-4.5)(5.89,-0.92)
\psline[linecolor=aqaqaq](5.89,-0.92)(9.13,1.92)
\psline[linecolor=dcrutc](5.89,-0.92)(9.13,1.92)
\psline[linecolor=dcrutc](9.13,1.92)(6.37,1.92)
\psline[linecolor=dcrutc](6.37,1.92)(-0.96,-4.5)
\psline[linecolor=dcrutc](-0.96,-4.5)(5.89,-0.92)
\psline[linecolor=qqwuqq](5.89,-0.92)(9.13,1.92)
\psline[linecolor=qqwuqq](9.13,1.92)(6.37,1.92)
\psline[linecolor=qqwuqq](6.37,1.92)(-0.96,-4.5)
\psline[linecolor=qqwuqq](-0.96,-4.5)(5.89,-0.92)
\begin{scriptsize}
\psdots[dotstyle=*,linecolor=blue](-4.3,1.94)
\psdots[dotstyle=*,linecolor=blue](20.1,1.9)
\psdots[dotstyle=*,linecolor=blue](-0.96,-4.5)
\rput[bl](-0.89,-4.39){\blue{$E$}}
\psdots[dotstyle=*,linecolor=blue](5.89,-0.92)
\rput[bl](5.97,-0.81){\blue{$Z$}}
\psdots[dotstyle=*,linecolor=blue](6.37,1.92)
\rput[bl](6.44,2.03){\blue{$X$}}
\psdots[dotstyle=*,linecolor=darkgray](6.37,1.92)
\psdots[dotstyle=*,linecolor=darkgray](11.32,1.91)
\psdots[dotstyle=*,linecolor=darkgray](9.13,1.92)
\rput[bl](9.2,2.03){\darkgray{$XEZ$}}
\end{scriptsize}
\end{pspicture*}
\end{center}

We have thus defined a map $\rho(z):a \to a$, $x \mapsto (xez)$.
Oberve that the map $\rho(z)$ described in the drawing has a {\em fixed point}:
which one? 
Using the intercept theorem with respect to that fixed point, 
show that $\rho(z)$ is a {\em dilation of the line $a$}.
For which position of $z$ does $\rho(z)$ not admit any fixed points?
Show that in this case $\rho(z)$ is a {\em translation of $a$}.
Summing up, $\rho(z)$ is an {\em affine self-map of $a$}.
Using the generalized associativity law, show that
$$
\rho( z z') = \rho(z') \circ \rho(z)
$$
(thus
{\em the group $(G,e)$ acts, via $\rho$, from the right on the line $a$}).
Altogether, conclude that this establishes an isomorphism
between the group $(G,e)$ and the affine group of the line. 
Since the affine group of a line is  a semidirect product of 
the additive group of a field with its multiplicative group, this is in keeping with what 
we observed above.

\subsubsection{Transformations}

In your drawing of the trapezoid configuration $x,y,z,w$,  make all lines and
auxiliary points invisible (in geogebra, click  ``object properties'', then 
turn off ``show object''), so that only the four points remain
(we won't reproduce that image here since its ``static'' version is not interesting at all). 
Then move $x$ and study the behaviour of $w$,
and similarly for $z$ and $y$: what can you say about the map $x \mapsto w$,
resp.\ $y \mapsto w$, resp.\ $z \mapsto w$?
For instance, are there fixed points of these maps in $G$, or in $P$?
Are these maps bijections, and can they be of finite order?
Can you reconstruct the line $a$ from the behavior of these maps?
    
\subsubsection{Powers and lattices}\label{subsec:Geo-Powers}

Let us come back to our group $(G,e)$.
In a group with neutral element $e$, the {\em powers} of an element $x$ are defined by
$x^{n+1}=x x^n$ and $x^{-n}=(x\inv)^n$.

\ssk
Exercise: fix an origin $e$ and construct  $x^2 = (xex)$.
Note that this  is a special case of the collinear
construction from subsection \ref{subsec:TrapezCollin}: you have to choose an auxiliary element $z$.
Hence  $x^2 = (xex)=(xe(zzx))= ((xez)zx)$ is given by
taking the product of $(xez)$ and $x$ with respect to $z$:

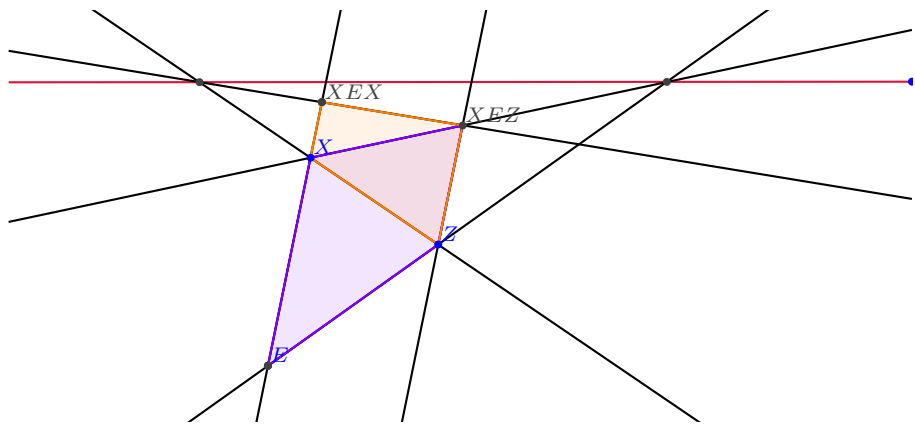
\begin{center}
\newrgbcolor{dcrutc}{0.86 0.08 0.24}
\newrgbcolor{xfqqff}{0.5 0 1}
\newrgbcolor{ffxfqq}{1 0.5 0}
\psset{xunit=0.5cm,yunit=0.4cm,algebraic=true,dotstyle=o,dotsize=3pt 0,linewidth=0.8pt,arrowsize=3pt 2,arrowinset=0.25}
\begin{pspicture*}(-4.3,-7.32)(19.46,6.3)
\pspolygon[linecolor=xfqqff,fillcolor=xfqqff,fillstyle=solid,opacity=0.1](3.64,1.42)(2.52,-5.48)(7,-1.46)(7.64,2.49)
\pspolygon[linecolor=ffxfqq,fillcolor=ffxfqq,fillstyle=solid,opacity=0.1](3.94,3.26)(3.64,1.42)(7,-1.46)(7.64,2.49)
\psplot[linecolor=dcrutc]{-4.3}{19.46}{(--94.17--0.02*x)/24}
\psplot{-4.3}{19.46}{(-23.53--6.9*x)/1.12}
\psplot{-4.3}{19.46}{(-34.68--4.02*x)/4.48}
\psplot{-4.3}{19.46}{(-49.94--6.9*x)/1.12}
\psplot{-4.3}{19.46}{(--4.15--2.51*x)/9.37}
\psplot{-4.3}{19.46}{(--15.25-2.88*x)/3.36}
\psplot{-4.3}{19.46}{(--28.2-1.43*x)/6.92}
\psline[linecolor=xfqqff](3.64,1.42)(2.52,-5.48)
\psline[linecolor=xfqqff](2.52,-5.48)(7,-1.46)
\psline[linecolor=xfqqff](7,-1.46)(7.64,2.49)
\psline[linecolor=xfqqff](7.64,2.49)(3.64,1.42)
\psline[linecolor=ffxfqq](3.94,3.26)(3.64,1.42)
\psline[linecolor=ffxfqq](3.64,1.42)(7,-1.46)
\psline[linecolor=ffxfqq](7,-1.46)(7.64,2.49)
\psline[linecolor=ffxfqq](7.64,2.49)(3.94,3.26)
\begin{scriptsize}
\psdots[dotstyle=*,linecolor=blue](-4.54,3.92)
\psdots[dotstyle=*,linecolor=blue](19.46,3.94)
\psdots[dotstyle=*,linecolor=blue](2.52,-5.48)
\rput[bl](2.6,-5.36){\blue{$E$}}
\psdots[dotstyle=*,linecolor=blue](3.64,1.42)
\rput[bl](3.72,1.54){\blue{$X$}}
\psdots[dotstyle=*,linecolor=blue](7,-1.46)
\rput[bl](7.08,-1.34){\blue{$Z$}}
\psdots[dotstyle=*,linecolor=darkgray](13.01,3.93)
\psdots[dotstyle=*,linecolor=darkgray](7.64,2.49)
\rput[bl](7.72,2.62){\darkgray{$XEZ$}}
\psdots[dotstyle=*,linecolor=darkgray](2.52,-5.48)
\psdots[dotstyle=*,linecolor=darkgray](2.52,-5.48)
\psdots[dotstyle=*,linecolor=darkgray](0.72,3.92)
\psdots[dotstyle=*,linecolor=darkgray](3.94,3.26)
\rput[bl](4.02,3.38){\darkgray{$XEX$}}
\psdots[dotstyle=*,linecolor=darkgray](3.94,3.26)
\psdots[dotstyle=*,linecolor=darkgray](3.94,3.26)
\end{scriptsize}
\end{pspicture*}
\end{center}

Move the point $z$: you will notice that the point $x^2$ does not depend on the 
choixe of $z$. 
Next construct $x\inv$. Since  $x\inv = (exe)$ is nothing but the square of $e$ with
respect to the origin $x$, we can use the same kind of construction as above:

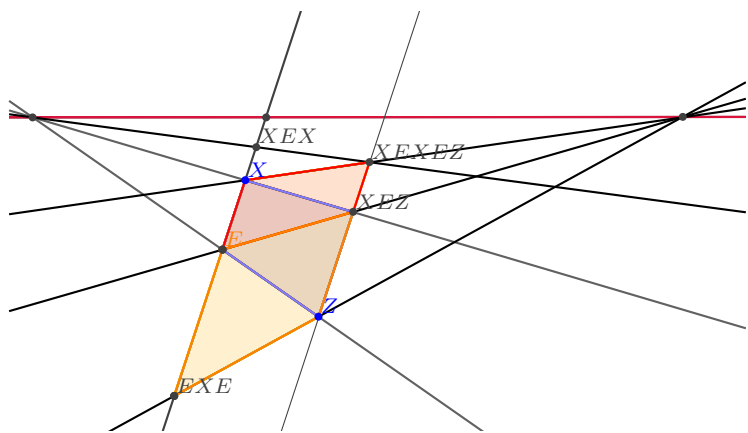
\begin{center}
\newrgbcolor{dcrutc}{0.86 0.08 0.24}
\newrgbcolor{ffxfqq}{1 0.5 0}
\newrgbcolor{wqwqwq}{0.38 0.38 0.38}
\newrgbcolor{uququq}{0.25 0.25 0.25}
\newrgbcolor{zzttqq}{0.6 0.2 0}
\newrgbcolor{ffdxqq}{1 0.84 0}
\newrgbcolor{xdxdff}{0.49 0.49 1}
\psset{xunit=0.5cm,yunit=0.5cm,algebraic=true,dotstyle=o,dotsize=3pt 0,linewidth=0.8pt,arrowsize=3pt 2,arrowinset=0.25}
\begin{pspicture*}(-3.18,-2.41)(16.2,8.71)
\pspolygon[linecolor=zzttqq,fillcolor=zzttqq,fillstyle=solid,opacity=0.1](3.03,4.23)(2.43,2.39)(4.96,0.61)(5.87,3.39)
\pspolygon[linecolor=ffdxqq,fillcolor=ffdxqq,fillstyle=solid,opacity=0.1](2.43,2.39)(5.87,3.39)(6.3,4.71)(3.03,4.23)
\pspolygon[linecolor=yellow,fillcolor=yellow,fillstyle=solid,opacity=0.1](1.16,-1.49)(4.96,0.61)(5.87,3.39)(2.43,2.39)
\pspolygon[linecolor=yellow,fillcolor=yellow,fillstyle=solid,opacity=0.1](2.43,2.39)(2.43,2.39)(2.43,2.39)
\pspolygon[linecolor=xdxdff,fillcolor=xdxdff,fillstyle=solid,opacity=0.1](3.03,4.23)(2.43,2.39)(4.96,0.61)(5.87,3.39)
\pspolygon[linecolor=red,fillcolor=red,fillstyle=solid,opacity=0.1](3.03,4.23)(6.3,4.71)(5.87,3.39)(2.43,2.39)
\pspolygon[linecolor=ffxfqq,fillcolor=ffxfqq,fillstyle=solid,opacity=0.1](2.43,2.39)(1.16,-1.49)(4.96,0.61)(5.87,3.39)
\psplot{-3.18}{16.2}{(--145.92--0.01*x)/24.72}
\psplot[linecolor=dcrutc]{-3.18}{16.2}{(-145.92-0.01*x)/-24.72}
\psplot[linecolor=dcrutc]{-3.18}{16.2}{(-145.92-0.01*x)/-24.72}
\psplot[linecolor=wqwqwq]{-3.18}{16.2}{(-10.37--1.78*x)/-2.53}
\psplot[linecolor=uququq]{-3.18}{16.2}{(-3.04--1.84*x)/0.6}
\psplot[linewidth=0.4pt,linecolor=uququq]{-3.18}{16.2}{(-8.78--1.84*x)/0.6}
\psplot[linecolor=wqwqwq]{-3.18}{16.2}{(--28.76-1.67*x)/5.6}
\psplot{-3.18}{16.2}{(--5.79--1*x)/3.44}
\psplot{-3.18}{16.2}{(--43.63--1.67*x)/11.51}
\psplot{-3.18}{16.2}{(-49.27--1.19*x)/-8.87}
\psplot{-3.18}{16.2}{(-20.43--5.3*x)/9.58}
\psline[linecolor=zzttqq](3.03,4.23)(2.43,2.39)
\psline[linecolor=zzttqq](2.43,2.39)(4.96,0.61)
\psline[linecolor=zzttqq](4.96,0.61)(5.87,3.39)
\psline[linecolor=zzttqq](5.87,3.39)(3.03,4.23)
\psline[linecolor=ffdxqq](2.43,2.39)(5.87,3.39)
\psline[linecolor=ffdxqq](5.87,3.39)(6.3,4.71)
\psline[linecolor=ffdxqq](6.3,4.71)(3.03,4.23)
\psline[linecolor=ffdxqq](3.03,4.23)(2.43,2.39)
\psline[linecolor=yellow](1.16,-1.49)(4.96,0.61)
\psline[linecolor=yellow](4.96,0.61)(5.87,3.39)
\psline[linecolor=yellow](5.87,3.39)(2.43,2.39)
\psline[linecolor=yellow](2.43,2.39)(1.16,-1.49)
\psline[linecolor=yellow](2.43,2.39)(2.43,2.39)
\psline[linecolor=yellow](2.43,2.39)(2.43,2.39)
\psline[linecolor=yellow](2.43,2.39)(2.43,2.39)
\psline[linecolor=xdxdff](3.03,4.23)(2.43,2.39)
\psline[linecolor=xdxdff](2.43,2.39)(4.96,0.61)
\psline[linecolor=xdxdff](4.96,0.61)(5.87,3.39)
\psline[linecolor=xdxdff](5.87,3.39)(3.03,4.23)
\psline[linecolor=red](3.03,4.23)(6.3,4.71)
\psline[linecolor=red](6.3,4.71)(5.87,3.39)
\psline[linecolor=red](5.87,3.39)(2.43,2.39)
\psline[linecolor=red](2.43,2.39)(3.03,4.23)
\psline[linecolor=ffxfqq](2.43,2.39)(1.16,-1.49)
\psline[linecolor=ffxfqq](1.16,-1.49)(4.96,0.61)
\psline[linecolor=ffxfqq](4.96,0.61)(5.87,3.39)
\psline[linecolor=ffxfqq](5.87,3.39)(2.43,2.39)
\begin{scriptsize}
\psdots[dotstyle=*,linecolor=blue](-5.33,5.9)
\psdots[dotstyle=*,linecolor=blue](19.39,5.91)
\psdots[dotstyle=*,linecolor=ffxfqq](2.43,2.39)
\rput[bl](2.5,2.49){\ffxfqq{$E$}}
\psdots[dotstyle=*,linecolor=blue](4.96,0.61)
\rput[bl](5.03,0.71){\blue{$Z$}}
\psdots[dotstyle=*,linecolor=blue](3.03,4.23)
\rput[bl](3.1,4.33){\blue{$X$}}
\psdots[dotstyle=*,linecolor=darkgray](3.58,5.9)
\psdots[dotstyle=*,linecolor=darkgray](-2.56,5.9)
\psdots[dotstyle=*,linecolor=darkgray](5.87,3.39)
\rput[bl](5.94,3.48){\darkgray{$XEZ$}}
\psdots[dotstyle=*,linecolor=darkgray](14.54,5.91)
\psdots[dotstyle=*,linecolor=darkgray](6.3,4.71)
\rput[bl](6.36,4.8){\darkgray{$XEXEZ$}}
\psdots[dotstyle=*,linecolor=darkgray](3.32,5.11)
\rput[bl](3.39,5.21){\darkgray{$XEX$}}
\psdots[dotstyle=*,linecolor=darkgray](2.43,2.39)
\psdots[dotstyle=*,linecolor=darkgray](2.43,2.39)
\psdots[dotstyle=*,linecolor=darkgray](2.43,2.39)
\psdots[dotstyle=*,linecolor=darkgray](1.16,-1.49)
\rput[bl](1.22,-1.4){\darkgray{$EXE$}}
\psdots[dotstyle=*,linecolor=darkgray](1.16,-1.49)
\psdots[dotstyle=*,linecolor=darkgray](1.16,-1.49)
\end{scriptsize}
\end{pspicture*}
\end{center}

Now you have seen the principle, construct at least ten powers 
of $x$, and then move again $x,e,z$:

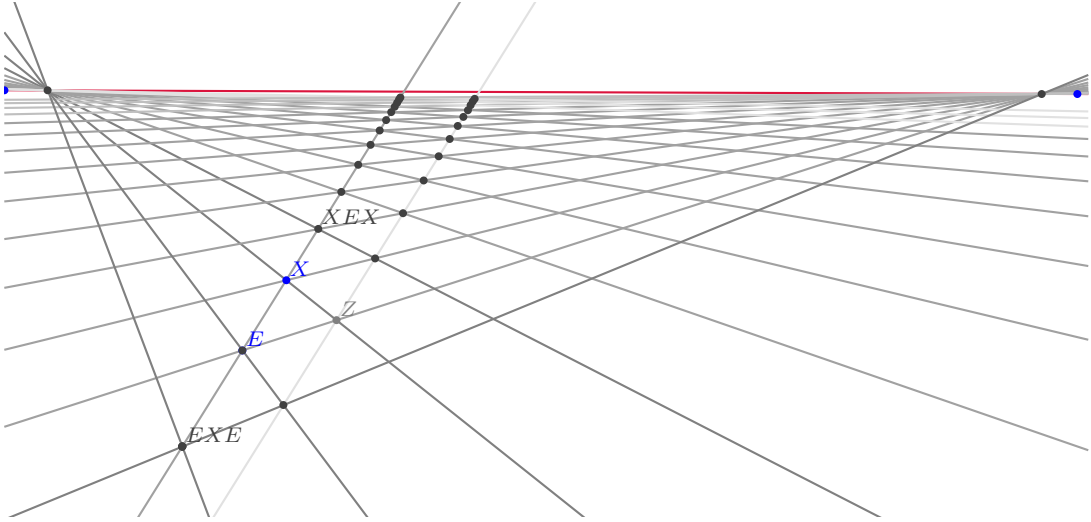
\begin{center}
\newrgbcolor{dcrutc}{0.86 0.08 0.24}
\newrgbcolor{aqaqaq}{0.63 0.63 0.63}
\newrgbcolor{eqeqeq}{0.88 0.88 0.88}
\psset{xunit=0.6cm,yunit=0.5cm,algebraic=true,dotstyle=o,dotsize=3pt 0,linewidth=0.8pt,arrowsize=3pt 2,arrowinset=0.25}
\begin{pspicture*}(-4.3,-7.32)(19.46,6.3)
\psplot[linecolor=dcrutc]{-4.3}{19.46}{(--93.18-0.1*x)/23.52}
\psplot[linecolor=aqaqaq]{-4.3}{19.46}{(--4.51-1.86*x)/-0.96}
\psplot[linecolor=aqaqaq]{-4.3}{19.46}{(-6.75--0.8*x)/2.06}
\psplot[linecolor=eqeqeq]{-4.3}{19.46}{(--7.58-1.86*x)/-0.96}
\psplot[linecolor=aqaqaq]{-4.3}{19.46}{(-26.85--4.94*x)/16.56}
\psplot[linecolor=gray]{-4.3}{19.46}{(--0.83-1.06*x)/1.1}
\psplot[linecolor=gray]{-4.3}{19.46}{(--13.61-4.45*x)/7.17}
\psplot[linecolor=aqaqaq]{-4.3}{19.46}{(-4.56--3.59*x)/15.86}
\psplot[linecolor=aqaqaq]{-4.3}{19.46}{(-20.07--3.26*x)/-7.79}
\psplot[linecolor=aqaqaq]{-4.3}{19.46}{(--11.68--2.6*x)/15.35}
\psplot[linecolor=aqaqaq]{-4.3}{19.46}{(--24.78-2.39*x)/8.24}
\psplot[linecolor=aqaqaq]{-4.3}{19.46}{(--23.52--1.88*x)/14.98}
\psplot[linecolor=aqaqaq]{-4.3}{19.46}{(-28.21--1.75*x)/-8.57}
\psplot[linecolor=aqaqaq]{-4.3}{19.46}{(--32.15--1.35*x)/14.71}
\psplot[linecolor=aqaqaq]{-4.3}{19.46}{(-30.72--1.29*x)/-8.81}
\psplot[linecolor=aqaqaq]{-4.3}{19.46}{(--38.45--0.97*x)/14.51}
\psplot[linecolor=aqaqaq]{-4.3}{19.46}{(-32.54--0.95*x)/-8.98}
\psplot[linecolor=aqaqaq]{-4.3}{19.46}{(--43.03--0.69*x)/14.36}
\psplot[linecolor=aqaqaq]{-4.3}{19.46}{(-33.87--0.7*x)/-9.11}
\psplot[linecolor=aqaqaq]{-4.3}{19.46}{(--46.38--0.49*x)/14.26}
\psplot[linecolor=aqaqaq]{-4.3}{19.46}{(-34.84--0.52*x)/-9.2}
\psplot[linecolor=lightgray]{-4.3}{19.46}{(--48.81--0.34*x)/14.18}
\psplot[linecolor=eqeqeq]{-4.3}{19.46}{(-35.54--0.39*x)/-9.27}
\psplot[linecolor=lightgray]{-4.3}{19.46}{(--50.59--0.23*x)/14.13}
\psplot[linecolor=eqeqeq]{-4.3}{19.46}{(-36.06--0.3*x)/-9.32}
\psplot[linecolor=lightgray]{-4.3}{19.46}{(--51.89--0.15*x)/14.09}
\psplot[linecolor=eqeqeq]{-4.3}{19.46}{(-36.43--0.23*x)/-9.35}
\psplot[linecolor=lightgray]{-4.3}{19.46}{(--52.83--0.1*x)/14.06}
\psplot[linecolor=gray]{-4.3}{19.46}{(-6.11-6.9*x)/4.27}
\psplot[linecolor=gray]{-4.3}{19.46}{(-87.66--8.25*x)/16.62}
\psplot[linecolor=gray]{-4.3}{19.46}{(-19.89-9.45*x)/2.95}
\begin{scriptsize}
\psdots[dotstyle=*,linecolor=blue](-4.3,3.98)
\psdots[dotstyle=*,linecolor=blue](19.22,3.88)
\psdots[dotstyle=*,linecolor=blue](0.92,-2.92)
\rput[bl](1,-2.8){\blue{$E$}}
\psdots[dotstyle=*,linecolor=blue](1.88,-1.06)
\rput[bl](1.96,-0.94){\blue{$X$}}
\psdots[dotstyle=*,linecolor=gray](2.98,-2.12)
\rput[bl](3.06,-2){\gray{$Z$}}
\psdots[dotstyle=*,linecolor=darkgray](18.44,3.88)
\psdots[dotstyle=*,linecolor=darkgray](3.83,-0.48)
\psdots[dotstyle=*,linecolor=darkgray](-3.35,3.98)
\psdots[dotstyle=*,linecolor=darkgray](2.58,0.3)
\rput[bl](2.66,0.42){\darkgray{$XEX$}}
\psdots[dotstyle=*,linecolor=darkgray](0.92,-2.92)
\psdots[dotstyle=*,linecolor=darkgray](4.44,0.72)
\psdots[dotstyle=*,linecolor=darkgray](3.09,1.28)
\psdots[dotstyle=*,linecolor=darkgray](4.89,1.59)
\psdots[dotstyle=*,linecolor=darkgray](3.46,2)
\psdots[dotstyle=*,linecolor=darkgray](5.22,2.23)
\psdots[dotstyle=*,linecolor=darkgray](3.73,2.53)
\psdots[dotstyle=*,linecolor=darkgray](5.46,2.69)
\psdots[dotstyle=*,linecolor=darkgray](3.93,2.91)
\psdots[dotstyle=*,linecolor=darkgray](5.64,3.03)
\psdots[dotstyle=*,linecolor=darkgray](4.07,3.19)
\psdots[dotstyle=*,linecolor=darkgray](5.76,3.27)
\psdots[dotstyle=*,linecolor=darkgray](4.18,3.4)
\psdots[dotstyle=*,linecolor=darkgray](5.86,3.45)
\psdots[dotstyle=*,linecolor=darkgray](4.26,3.54)
\psdots[dotstyle=*,linecolor=darkgray](5.92,3.58)
\psdots[dotstyle=*,linecolor=darkgray](4.31,3.65)
\psdots[dotstyle=*,linecolor=darkgray](5.97,3.68)
\psdots[dotstyle=*,linecolor=darkgray](4.35,3.73)
\psdots[dotstyle=*,linecolor=darkgray](6.01,3.75)
\psdots[dotstyle=*,linecolor=darkgray](4.38,3.79)
\psdots[dotstyle=*,linecolor=darkgray](1.82,-4.37)
\psdots[dotstyle=*,linecolor=darkgray](-0.4,-5.47)
\rput[bl](-0.32,-5.36){\darkgray{$EXE$}}
\psdots[dotstyle=*,linecolor=darkgray](-0.4,-5.47)
\psdots[dotstyle=*,linecolor=darkgray](-0.4,-5.47)
\psdots[dotstyle=*,linecolor=darkgray](-0.4,-5.47)
\psdots[dotstyle=*,linecolor=darkgray](-0.4,-5.47)
\end{scriptsize}
\end{pspicture*}
\end{center}
    
The powers $\{ x^n \mid n\in \Z \}$ of an element $x \in G$ form a
{\em discrete subgroup of $G$}.

\section{Computational approach: linear algebra}\label{sec:Compute}

Geometry can be translated into ``analytic formulas'' by identifying the plane $E$ with
the real vector space $\R^2$. 
Vectors are written $x=(x_1,x_2)$, $y=(y_1,y_2)$, and  their sum is defined by
$x+y=(x_1+y_1,x_2+y_2)$.
In the following, the reader should avoid as much as possible to use the components
(that is, avoid to use the standard basis of $\R^2$), but just the vector space properties of  $E$.

\subsection{Some reminders}\label{subsec:ComputeReminders}

For two distinct points $x,y \in E$, recall that the {\em affine line spanned by $x$ and $y$} is
\begin{equation}
x \lor y = \bigl\{ t x + (1-t) y \mid \, t \in \R \bigr\} \, .
\end{equation}
Let $z \lor w$ be another such line. Prove that the following are equivalent:

\begin{description}
\item{(1)}
the vector $z-w$ is a multiple of the vector $x-y$,
\item{(2)}
either $x \lor y = z \lor w$, or $(x\lor y) \cap (z\lor w) = \emptyset$.
\end{description}

We then say that the lines $x \lor y$ and $z\lor w$ are {\em parallel}, and we write
$(x\lor y) \parallel (z \lor w)$.

\subsection{Parallelograms and vector addition}\label{subsec:ComputePar}

\subsubsection{Usual parallelograms}

Let $(x,y,z,w)$ a quadruple of points in the plane $E$, no two of them equal. 
Show that the following are equivalent:

\begin{description}
\item{(1)} $(x\lor y) \parallel (z \lor w)$ and $(x\lor w) \parallel (y\lor z)$,
\item{(2)} $w = x - y + z$,
\item{(3)} the quadrangle with vertices $xyzw$ (in this order) is a parallelogram.
\end{description}

Show that, for a fixed element $y \in E$, the law $x+_y z:= x-y+z$ defines 
a commutative group law on $E$ with neutral element $y$.
For $y=0$, we get back  usual vector addition.

\msk
{\sl Some analysis.}
Recall the notion of {\em open} subsets of $\R^n$. We may identify $E^3$ with $\R^6$. 
Show   that the set $U:=\{ (x,y,z) \in E^3 \vert \, x \not= y \not= z, x \not= z \}$ is {\em open and dense}
in $E^3$ and that the map
$$
U  \to E, \quad (x,y,z) \mapsto x-y+z
$$
is {\em continuous}. Show that this map admits a unique continuous extension to a map
$E^3 \to E$, given by the same formula.

\subsubsection{Perspective view}\label{subsec:ComputationPerspective}

This exercise is optional (it is not needed in the sequel).
Let $\alpha : E \to \R$ to be a non-zero  linear form and
$a:= \ker(\alpha)$. Recall that $\dim a = \dim E - 1 =1$, so $a$ is a line.
(You may take $\alpha(x) = x_2$ to get a horizontal line.)
Let  $x,y,z$ in the plane such that the lines
$x \lor y$ and $z \lor y$ are not parallel to $a$. 
Show that the point $w$ defined by Equation (\ref{eqn:ParPersp}) is given by the formula
\begin{equation}
w = \frac{\alpha(y) \alpha(z)}{d} x  -
\frac{\alpha(x)\alpha(z)}{d} y + \frac{\alpha(y)\alpha(x)}{d} z \, ,
\end{equation}
with denominator $d = \alpha(y)\alpha(x) -\alpha(x)\alpha(z) +\alpha(y)\alpha(z)$.
(Hint: computations may be simplified by observing that the given expression of $w$
is a certain {\em barycenter} of the three points $x,y,z$.)  
Note that this formula makes sense also if $x,y,z$ are collinear, and also
if one  of these points belongs  to the line $a$
(what do you get if $x \in a$?  if $z \in a$? if $y \in a$?). 
  
  \ssk
Describe the set of all $(x,y,z) \in E^3 = \R^6$ such that $d\not=0$ and
show that it is open and dense in $E^3$, and that the map
$(x,y,z) \mapsto w$ is continuous on this set.
Compare with your image from section \ref{sec:ParPersp}: what happens
when $d=0$?

\subsection{The linear algebra of trapezoid geometry}\label{subsec:ComputationTrapezoid}

In this section, we fix a line $a$ given by $a=\{ x \in \R^2 \mid \, \alpha(x)=0 \}$ for some non-zero
linear form $\alpha:\R^2 \to \R$; 
to keep close to the drawings, you may choose $\alpha(x)=-x_2$ (but try to
avoid using  this explicit formula).
The set $G$ is defined by
\begin{equation}
G := \{ x \in E \mid \alpha(x) \not= 0 \} .
\label{eqn:G-def}
\end{equation}

\subsubsection{Main formula}

Show that the point $w$ defined  by Equation (\ref{eqn:Generic}) is given by the formula
\begin{equation}
\boxed{
(xyz) = w = \alpha(z) \alpha(y)\inv (x-y) + z } \  .\label{eqn:Main}
\end{equation}
Hints:
show first that $w$ is of the form $z+t(x-y)$ with some $t\in \R$; then determine $t$ by
using the intercept theorem with center $u$ being the intersection $a \land (y \lor z)$  -- 
show that the ratio $\vert w-z \vert / \vert x-y \vert$ is equal to the ratio
$\vert u - z \vert / \vert u - y\vert$, which in turn is equal to $\alpha(z)/\alpha(y)$,
and conclude that $t = \alpha(z)\alpha(y)\inv$.

\subsubsection{Group structure}

Fix a point $e \in G$ such that $\alpha(e)=1$. Show that $G$ is a group with
neutral element $e$ and product
\begin{equation}
x \cdot z = (xez)= \alpha(z) (x-e) + z . \label{eqn:GroupProduct}
\end{equation}
Hint: start by proving the very useful formula
\begin{equation}
\boxed{
\alpha \bigl( (xyz) \bigr) = \alpha(x) \alpha(y)\inv  \alpha(z) } 
\label{eqn:HomMult}
\end{equation}
which gives you
\begin{equation}
\alpha(x \cdot z) = \alpha(x) \alpha(z),
\label{eqn:Grouphom}
\end{equation}
and deduce from this that $G$ is indeed stable under the product and next that
the product is associative.
Then show that $e \cdot x = x =x\cdot e$, and that the inverse of $x$ is given by
$$
x\inv = \frac{1}{\alpha(x)} (e-x) + x .
$$

\msk
{\sl Some analysis.}
Show that $G$ is open dense in $E$ and that the group law
$G \times G \to G$ and the inversion map $G\to G$, $x \mapsto x\inv$ are
continuous and even smooth of class $C^\infty$. 

\msk
{\sl Group structure: semidirect product.}
Show that $\alpha\vert_G :G \to \R^\times$ is a group homomorphism
and that its kernel $N=\{ x \in G \mid \alpha(x)=1 \}$ is a subgroup isomorphic to $(\R,+)$.
Show that the set
$$
H:= \R^\times e = \{ r e \mid \, r \in \R^\times \}
$$
is a subgroup of $G$, and that $\phi: \R^\times \to H$, $r \mapsto r e$ is a group
isomorphism such that $\alpha (\phi (t)) = t$. 
This means that $G$ is a semidirect product of $N \cong (\R,+)$ with $H \cong (\R^\times,\cdot)$.
For a complete description of $G$ it remains to determine the action of $H$ on the normal subgroup
$N$: show that this action is the ``standard action up to an inversion'', that is,
$$
\R^\times \times \R \to \R, \qquad (r,t) \mapsto r\inv t .
$$
Hint: either use a direct computation to show that 
$$
re \cdot x \cdot (re)\inv = e + r\inv (x-e) ,
$$
or change the origin as explained below (which is better adapted to the additive structure
of the  normal subgroup). In any case, you may now
 conclude that $G$ is isomorphic to a semidirect product $\R \rtimes \R^\times$,
 and hence $G$ is isomorphic to the 
``$ax+b$-group'', i.e., to the affine group of the real line, in matrix form:
\begin{equation}\label{eqn:Ga1}
\Ga(1,\R) := 
\Big\{
\begin{pmatrix} a & b \cr 0 & 1 \end{pmatrix} \mid \, 
a \in \R^\times, b \in \R \Big\} .
\end{equation}

\ssk
{\sl Change of origin.}
Sometimes it may be useful to present the group $G$ in a different way: we translate
the neutral element $e$ to the origin $0$ of the vector space $E$, that is, we
translate $G$ to the set $G':= \{ x \in E \mid \alpha(x) \not= 1 \}$, with the new
product
$(x,y) \mapsto (x+e)\cdot (y+e) - e$. Show that this new product is given by 
\begin{equation}
G' \times G' \to G' , \qquad (x,z) \mapsto  x + \alpha(x)z + z
\end{equation}
and the new inverse by
\begin{equation}
x\inv = - \frac{1}{\alpha(x)+1} x .
\end{equation}
Compare with formula (\ref{eqn:GroupProduct}):
in (\ref{eqn:GroupProduct}) it is not possible to take easily a 
``limit $\alpha \to 0$'', whereas here this is possible: for $\alpha = 0$ one gets
vector addition as limit case!

\subsubsection{Generalized associativity: torsors}\label{subsec:ComputationAssoc}

Using (\ref{eqn:HomMult}), show that $G$ is stable under the ternary
product $(x,y,z) \mapsto (xyz)$, and  that the ternary product $(xyz)$ satisfies the generalized associativity law
(\ref{eqn:torsor2}) and the idempotency law (\ref{eqn:Idempot}).
In other words, show  that $G$ with product $(xyz)$ given by formula
(\ref{eqn:Main}) is a {\em torsor}: 

\msk
{\bf Definition.}
{\em A {\em torsor}\footnote{A remark on terminology is in order here: 
torsors are for groups what affine spaces are for vector spaces, and
the concept named here
``torsor'' is at least 70 years old, see \cite{Cer43}; however, it has not really become part
of the mathematical mainstream, possibly because there is no universally accepted terminology --
other terms such as {\em groud, heap, flock, principal homogeneous space} are also used,
see \cite{BeKi10} for some more comments on the terminology we use.}
 is a set $G$ together with a ternary ``product map''
$G^3 \to G$, $(x,y,z) \mapsto (xyz)$
satisfying the idempotency law (\ref{eqn:Idempot})
and the generalized associative law (\ref{eqn:torsor2}). 
A {\em homomorphism of torsors} is a map between torsors respecting the ternary
products.
A torsor is called {\em commutative}  if it satisfies the identity  $(xyz)=(zyx)$.
}

\msk
In the following exercises let us collect some properties that hold for {\em all}
torsors: 

\msk
{\sl
From torsors to groups.} Let $G$ be a torsor, and an arbitrary point $e \in G$.
Show that $(G,e)$ is a group with product $x \cdot z:= (xez)$.
(Hint: show that the inverse element of $x$ is $x\inv = (exe)$.)

\msk
{\sl From groups to torsors.}
Let $(G,e)$ be a group. Show that $G$ with ternary law 
$(xyz)=xy\inv z$ is a torsor. 
(Note: when $G$ is commutative, the group law is often written additively, and
then torsor law becomes $x-y+z$. Example: the torsors $\R^n$. On the other hand,
for the commutative torsor $\R^\times$ one uses of course a multiplicative notation.)

\msk
{\sl Back and forth.} Show that the preceding two constructions are inverse to each other.
Conclude that torsors are ``the same as groups where one has forgotten the origin''.
(If you prefer more formal statements:  there is an equivalence of categories
between the categories of groups and of torsors with a distinguished base point.) 

\msk
{\sl Para-associativity.}
Show that in any torsor the {\em para-associative law}
\begin{equation}
(xw(vuz)) = (x(uvw)z) = ((xwv)uz) \label{eqn:ParaAss} 
\end{equation}
holds. (Remark: the analog of a semigroup is a {\em semitorsor} -- one just demands
(\ref{eqn:ParaAss}) to hold, but not necessarily the idempotency law.) 

\msk
{\sl Subtorsors.}
A {\em subtorsor} is a subset $H$ of a torsor $G$ which is stable under the ternary law.
Fix a point $e \in G$. Show that
 $H \subset G$ is a subgroup iff $H$ is a subtorsor containing $e$.
 (Remark: this characterization is simpler than the one for groups, since we do not
have to investigate separately that inverses again belong to $H$.

\msk
{\sl Symmetry relations.} 
Show that for a quadruple of elements of a torsor  the following four conditions are equivalent:
\begin{equation}
(xyz)=w, \qquad
(zwx)=y, \qquad
(yxw)=x,\qquad
(wzy)=x ¥, .
\end{equation}
In other words, the graph $R$ of the ternary product map $G^3 \to G$,
$$
R = \{ (x,y,z,w)\in G^4 \mid \, w = (xyz) \},
$$
is invariant under double
permutations, i.e., under the {\em Klein four group}. 
Deduce that in any torsor the identities
(\ref{eqn:Perminvar}) are satisfied.
Assume that the torsor is commutative:
which subgroup of permutations of the four letters $(x,y,z,w)$ leaves $R$ invariant then? 
Compare with the proof given in section \ref{subsec:non-com}.

\msk
{\sl Left- and Right translations.}
In any torsor, define {\em left and right translation operators} by
\begin{equation}
L_{xy}:G \to G, \quad z \mapsto (xyz), \qquad \quad
R_{zy}:G \to G, \quad x \mapsto (xyz) \, .
\end{equation}
Show that the left translations and the right translations form groups of bijections of $G$,
both isomorphic to $(G,e)$, for any choice of origin $e \in G$.

\subsubsection{Collinear case: one-parameter subgroups and powers}\label{subsec:powers}

Now let again  $G$  be defined by (\ref{eqn:G-def}) with ternary torsor law
given by (\ref{eqn:Main}).
Let $e \in G$ be an arbitrary point.
Show that $G$ with product $(x,z) \mapsto (xez)$ is a group, and that the inverse in this group is
\begin{equation}
x\inv = (exe) =  \alpha(e) \alpha(x)\inv (e-x) + e .
\end{equation}
Give a computational proof of the affirmations from subsection \ref{subsec:TrapezCollin}:
let $\ell$ be an affine line in the plane $E$ and $e \in \ell$.
Show that $\ell$ is parallel to $a$ iff $\alpha(x)=\alpha(y)$ for all $x,y \in \ell$, 
and $\ell$ is not parallel to $a$ iff the restriction 
$\phi: \ell \to \R$, $x \mapsto \alpha(x)$ is a bijection.
Using this, show:

\begin{description}
\item{(1)}
if $\ell$ is parallel to $a$, then $\ell$ is a subgroup of $G$, isomorphic to $(\R,+)$
\item{(2)}
if $\ell$ is not parallel to $a$, then $\ell \cap G$ is a subgroup of $G$, isomorphic to $(\R^\times,
\cdot)$.
(Hint: use equation (\ref{eqn:HomMult}) in order to show that 
$\phi$, restricted to $\ell \cap G$, defines  an isomorphism of groups.)
\end{description}

Recall that the powers in the group $(G,e)$ are defined by
$x^2 =  (xex)$,
$x^3 = (xe(xex))$,
$x^{n+1} = (xe x^n)$. 
Prove the following explicit formula for the powers:
\begin{equation}\label{eqn:Power}
x^n =
\frac{\alpha(x)^n \alpha(e)^{-n+1} - \alpha(e)}{\alpha(x-e)} x +
\frac{\alpha(x) - \alpha(x)^n \alpha(e)^{-n+1}}{\alpha(x-e)} e 
\end{equation}
Hints:
we have seen that $\alpha:G \to \R^\times$ is a torsor homomorphism, and it
is an isomorphism when restricted to the line $(x \lor e) \cap G$ (if it is not parallel to $a$).
Show that the map
\begin{equation} \label{eqn:Section}
\phi : \R^\times \to G, \quad t \mapsto
\frac{t - \alpha(e)}{\alpha(x-e)} x + \frac{\alpha(x) - t}{\alpha(x-e)} e
\end{equation}
is a section of $\alpha$, i.e., it is a torsor homomorphism such that $\alpha(\phi(t))=t$.
It maps $\alpha(x)$ to $x$ and $\alpha(e)$ to $e$, and hence defines an isomorphism
$\R^\times \to (x \lor e)\cap G$. 
Deduce that $x^n = \phi ((\alpha(x)^n)$.

\msk
It is possible to plot the powers of $x$ by programming (\ref{eqn:Power}) in geogebra: 
here is the geogebra command for drawing the powers  $X^{-50},\ldots,X^{50}$ 
with respect to the origin $Y$, and using
the linear form $\alpha(x)=-x_2$:

\ssk
\begin{center}
Sequence[DynamicCoordinates$[X, (x(X) y(Y) - y(X) x(Y) + y(X)^{(n + 1)} y(Y)^{(-n)}
 (x(Y) - x(X))) (y(Y) - y(X))^{(-1)}, y(X)^{(n + 1)} y(Y)^{(-n)}], n, -50, 50]$
\end{center}

The illustration shows the powers $X^n$, $Z^n$, $U^n$ with respect to a common
origin $Y$.
The line $X \lor Y$ is close to being parallel to $a$, and hence the powers look
like an arithmetic sequence; the powers of the two other elements clearly form
a geometric sequence converging to a point of $a$ as $n$ tends to $+\infty$:

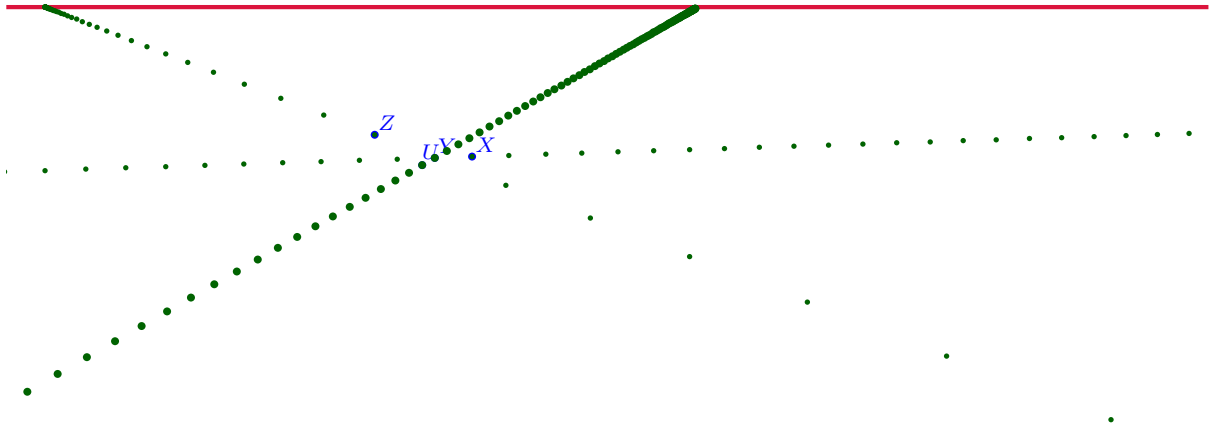
\begin{center}
\newrgbcolor{dcrutc}{0.86 0.08 0.24}
\newrgbcolor{qqwuqq}{0 0.39 0}
\psset{xunit=0.5cm,yunit=0.4cm,algebraic=true,dotstyle=o,dotsize=3pt 0,linewidth=0.8pt,arrowsize=3pt 2,arrowinset=0.25}
\begin{pspicture*}(-15.61,-14.21)(16.01,3.92)
\psplot[linewidth=1.6pt,linecolor=dcrutc]{-15.61}{16.01}{(-0-0*x)/20.02}
\begin{scriptsize}
\psdots[dotstyle=*,linecolor=blue](-4.34,-5)
\rput[bl](-4.25,-4.84){\blue{$Y$}}
\psdots[dotstyle=*,linecolor=blue](-3.36,-4.96)
\rput[bl](-3.26,-4.79){\blue{$X$}}
\psdots[dotsize=2pt 0,dotstyle=*,linecolor=qqwuqq](-15.66,-5.46)
\psdots[dotsize=2pt 0,dotstyle=*,linecolor=qqwuqq](-14.59,-5.42)
\psdots[dotsize=2pt 0,dotstyle=*,linecolor=qqwuqq](-13.52,-5.37)
\psdots[dotsize=2pt 0,dotstyle=*,linecolor=qqwuqq](-12.47,-5.33)
\psdots[dotsize=2pt 0,dotstyle=*,linecolor=qqwuqq](-11.42,-5.29)
\psdots[dotsize=2pt 0,dotstyle=*,linecolor=qqwuqq](-10.39,-5.25)
\psdots[dotsize=2pt 0,dotstyle=*,linecolor=qqwuqq](-9.36,-5.2)
\psdots[dotsize=2pt 0,dotstyle=*,linecolor=qqwuqq](-8.34,-5.16)
\psdots[dotsize=2pt 0,dotstyle=*,linecolor=qqwuqq](-7.33,-5.12)
\psdots[dotsize=2pt 0,dotstyle=*,linecolor=qqwuqq](-6.32,-5.08)
\psdots[dotsize=2pt 0,dotstyle=*,linecolor=qqwuqq](-5.33,-5.04)
\psdots[dotsize=2pt 0,dotstyle=*,linecolor=qqwuqq](-4.34,-5)
\psdots[dotsize=2pt 0,dotstyle=*,linecolor=qqwuqq](-3.36,-4.96)
\psdots[dotsize=2pt 0,dotstyle=*,linecolor=qqwuqq](-2.39,-4.92)
\psdots[dotsize=2pt 0,dotstyle=*,linecolor=qqwuqq](-1.42,-4.88)
\psdots[dotsize=2pt 0,dotstyle=*,linecolor=qqwuqq](-0.47,-4.84)
\psdots[dotsize=2pt 0,dotstyle=*,linecolor=qqwuqq](0.48,-4.8)
\psdots[dotsize=2pt 0,dotstyle=*,linecolor=qqwuqq](1.42,-4.76)
\psdots[dotsize=2pt 0,dotstyle=*,linecolor=qqwuqq](2.36,-4.73)
\psdots[dotsize=2pt 0,dotstyle=*,linecolor=qqwuqq](3.28,-4.69)
\psdots[dotsize=2pt 0,dotstyle=*,linecolor=qqwuqq](4.2,-4.65)
\psdots[dotsize=2pt 0,dotstyle=*,linecolor=qqwuqq](5.11,-4.61)
\psdots[dotsize=2pt 0,dotstyle=*,linecolor=qqwuqq](6.02,-4.58)
\psdots[dotsize=2pt 0,dotstyle=*,linecolor=qqwuqq](6.92,-4.54)
\psdots[dotsize=2pt 0,dotstyle=*,linecolor=qqwuqq](7.81,-4.5)
\psdots[dotsize=2pt 0,dotstyle=*,linecolor=qqwuqq](8.69,-4.47)
\psdots[dotsize=2pt 0,dotstyle=*,linecolor=qqwuqq](9.56,-4.43)
\psdots[dotsize=2pt 0,dotstyle=*,linecolor=qqwuqq](10.43,-4.4)
\psdots[dotsize=2pt 0,dotstyle=*,linecolor=qqwuqq](11.3,-4.36)
\psdots[dotsize=2pt 0,dotstyle=*,linecolor=qqwuqq](12.15,-4.33)
\psdots[dotsize=2pt 0,dotstyle=*,linecolor=qqwuqq](13,-4.29)
\psdots[dotsize=2pt 0,dotstyle=*,linecolor=qqwuqq](13.84,-4.26)
\psdots[dotsize=2pt 0,dotstyle=*,linecolor=qqwuqq](14.67,-4.22)
\psdots[dotsize=2pt 0,dotstyle=*,linecolor=qqwuqq](15.5,-4.19)
\psdots[dotsize=2pt 0,dotstyle=*,linecolor=qqwuqq](16.32,-4.16)
\psdots[dotstyle=*,linecolor=blue](-5.92,-4.23)
\rput[bl](-5.82,-4.07){\blue{$Z$}}
\psdots[dotsize=2pt 0,dotstyle=*,linecolor=qqwuqq](13.45,-13.68)
\psdots[dotsize=2pt 0,dotstyle=*,linecolor=qqwuqq](9.12,-11.57)
\psdots[dotsize=2pt 0,dotstyle=*,linecolor=qqwuqq](5.46,-9.78)
\psdots[dotsize=2pt 0,dotstyle=*,linecolor=qqwuqq](2.36,-8.27)
\psdots[dotsize=2pt 0,dotstyle=*,linecolor=qqwuqq](-0.25,-6.99)
\psdots[dotsize=2pt 0,dotstyle=*,linecolor=qqwuqq](-2.47,-5.91)
\psdots[dotsize=2pt 0,dotstyle=*,linecolor=qqwuqq](-4.34,-5)
\psdots[dotsize=2pt 0,dotstyle=*,linecolor=qqwuqq](-5.92,-4.23)
\psdots[dotsize=2pt 0,dotstyle=*,linecolor=qqwuqq](-7.26,-3.58)
\psdots[dotsize=2pt 0,dotstyle=*,linecolor=qqwuqq](-8.39,-3.02)
\psdots[dotsize=2pt 0,dotstyle=*,linecolor=qqwuqq](-9.35,-2.56)
\psdots[dotsize=2pt 0,dotstyle=*,linecolor=qqwuqq](-10.16,-2.16)
\psdots[dotsize=2pt 0,dotstyle=*,linecolor=qqwuqq](-10.84,-1.83)
\psdots[dotsize=2pt 0,dotstyle=*,linecolor=qqwuqq](-11.42,-1.55)
\psdots[dotsize=2pt 0,dotstyle=*,linecolor=qqwuqq](-11.91,-1.31)
\psdots[dotsize=2pt 0,dotstyle=*,linecolor=qqwuqq](-12.32,-1.11)
\psdots[dotsize=2pt 0,dotstyle=*,linecolor=qqwuqq](-12.67,-0.93)
\psdots[dotsize=2pt 0,dotstyle=*,linecolor=qqwuqq](-12.97,-0.79)
\psdots[dotsize=2pt 0,dotstyle=*,linecolor=qqwuqq](-13.22,-0.67)
\psdots[dotsize=2pt 0,dotstyle=*,linecolor=qqwuqq](-13.43,-0.57)
\psdots[dotsize=2pt 0,dotstyle=*,linecolor=qqwuqq](-13.61,-0.48)
\psdots[dotsize=2pt 0,dotstyle=*,linecolor=qqwuqq](-13.76,-0.4)
\psdots[dotsize=2pt 0,dotstyle=*,linecolor=qqwuqq](-13.89,-0.34)
\psdots[dotsize=2pt 0,dotstyle=*,linecolor=qqwuqq](-14,-0.29)
\psdots[dotsize=2pt 0,dotstyle=*,linecolor=qqwuqq](-14.09,-0.24)
\psdots[dotsize=2pt 0,dotstyle=*,linecolor=qqwuqq](-14.17,-0.21)
\psdots[dotsize=2pt 0,dotstyle=*,linecolor=qqwuqq](-14.23,-0.17)
\psdots[dotsize=2pt 0,dotstyle=*,linecolor=qqwuqq](-14.29,-0.15)
\psdots[dotsize=2pt 0,dotstyle=*,linecolor=qqwuqq](-14.33,-0.12)
\psdots[dotsize=2pt 0,dotstyle=*,linecolor=qqwuqq](-14.37,-0.11)
\psdots[dotsize=2pt 0,dotstyle=*,linecolor=qqwuqq](-14.41,-0.09)
\psdots[dotsize=2pt 0,dotstyle=*,linecolor=qqwuqq](-14.44,-0.08)
\psdots[dotsize=2pt 0,dotstyle=*,linecolor=qqwuqq](-14.46,-0.06)
\psdots[dotsize=2pt 0,dotstyle=*,linecolor=qqwuqq](-14.48,-0.05)
\psdots[dotsize=2pt 0,dotstyle=*,linecolor=qqwuqq](-14.5,-0.05)
\psdots[dotsize=2pt 0,dotstyle=*,linecolor=qqwuqq](-14.51,-0.04)
\psdots[dotsize=2pt 0,dotstyle=*,linecolor=qqwuqq](-14.52,-0.03)
\psdots[dotsize=2pt 0,dotstyle=*,linecolor=qqwuqq](-14.53,-0.03)
\psdots[dotsize=2pt 0,dotstyle=*,linecolor=qqwuqq](-14.54,-0.02)
\psdots[dotsize=2pt 0,dotstyle=*,linecolor=qqwuqq](-14.55,-0.02)
\psdots[dotsize=2pt 0,dotstyle=*,linecolor=qqwuqq](-14.56,-0.02)
\psdots[dotsize=2pt 0,dotstyle=*,linecolor=qqwuqq](-14.56,-0.01)
\psdots[dotsize=2pt 0,dotstyle=*,linecolor=qqwuqq](-14.57,-0.01)
\psdots[dotsize=2pt 0,dotstyle=*,linecolor=qqwuqq](-14.57,-0.01)
\psdots[dotsize=2pt 0,dotstyle=*,linecolor=qqwuqq](-14.57,-0.01)
\psdots[dotsize=2pt 0,dotstyle=*,linecolor=qqwuqq](-14.58,-0.01)
\psdots[dotsize=2pt 0,dotstyle=*,linecolor=qqwuqq](-14.58,-0.01)
\psdots[dotsize=2pt 0,dotstyle=*,linecolor=qqwuqq](-14.58,-0.01)
\psdots[dotsize=2pt 0,dotstyle=*,linecolor=qqwuqq](-14.58,0)
\psdots[dotsize=2pt 0,dotstyle=*,linecolor=qqwuqq](-14.58,0)
\psdots[dotsize=2pt 0,dotstyle=*,linecolor=qqwuqq](-14.58,0)
\psdots[dotsize=2pt 0,dotstyle=*,linecolor=qqwuqq](-14.58,0)
\psdots[dotsize=2pt 0,dotstyle=*,linecolor=qqwuqq](-14.59,0)
\psdots[dotsize=2pt 0,dotstyle=*,linecolor=qqwuqq](-14.59,0)
\psdots[dotsize=2pt 0,dotstyle=*,linecolor=qqwuqq](-14.59,0)
\psdots[dotsize=2pt 0,dotstyle=*,linecolor=qqwuqq](-14.59,0)
\psdots[dotsize=2pt 0,dotstyle=*,linecolor=qqwuqq](-14.59,0)
\psdots[dotsize=2pt 0,dotstyle=*,linecolor=qqwuqq](-14.59,0)
\psdots[dotstyle=*,linecolor=blue](-4.67,-5.24)
\rput[bl](-4.67,-5.03){\blue{$U$}}
\psdots[dotstyle=*,linecolor=qqwuqq](2.51,-0.05)
\psdots[dotstyle=*,linecolor=qqwuqq](2.51,-0.05)
\psdots[dotstyle=*,linecolor=qqwuqq](2.51,-0.05)
\psdots[dotstyle=*,linecolor=qqwuqq](2.5,-0.06)
\psdots[dotstyle=*,linecolor=qqwuqq](2.5,-0.06)
\psdots[dotstyle=*,linecolor=qqwuqq](2.5,-0.06)
\psdots[dotstyle=*,linecolor=qqwuqq](2.49,-0.06)
\psdots[dotstyle=*,linecolor=qqwuqq](2.49,-0.07)
\psdots[dotstyle=*,linecolor=qqwuqq](2.48,-0.07)
\psdots[dotstyle=*,linecolor=qqwuqq](2.48,-0.07)
\psdots[dotstyle=*,linecolor=qqwuqq](2.47,-0.08)
\psdots[dotstyle=*,linecolor=qqwuqq](2.47,-0.08)
\psdots[dotstyle=*,linecolor=qqwuqq](2.46,-0.09)
\psdots[dotstyle=*,linecolor=qqwuqq](2.46,-0.09)
\psdots[dotstyle=*,linecolor=qqwuqq](2.45,-0.09)
\psdots[dotstyle=*,linecolor=qqwuqq](2.45,-0.1)
\psdots[dotstyle=*,linecolor=qqwuqq](2.44,-0.1)
\psdots[dotstyle=*,linecolor=qqwuqq](2.43,-0.11)
\psdots[dotstyle=*,linecolor=qqwuqq](2.43,-0.11)
\psdots[dotstyle=*,linecolor=qqwuqq](2.42,-0.12)
\psdots[dotstyle=*,linecolor=qqwuqq](2.41,-0.12)
\psdots[dotstyle=*,linecolor=qqwuqq](2.4,-0.13)
\psdots[dotstyle=*,linecolor=qqwuqq](2.39,-0.14)
\psdots[dotstyle=*,linecolor=qqwuqq](2.38,-0.14)
\psdots[dotstyle=*,linecolor=qqwuqq](2.38,-0.15)
\psdots[dotstyle=*,linecolor=qqwuqq](2.37,-0.16)
\psdots[dotstyle=*,linecolor=qqwuqq](2.35,-0.16)
\psdots[dotstyle=*,linecolor=qqwuqq](2.34,-0.17)
\psdots[dotstyle=*,linecolor=qqwuqq](2.33,-0.18)
\psdots[dotstyle=*,linecolor=qqwuqq](2.32,-0.19)
\psdots[dotstyle=*,linecolor=qqwuqq](2.31,-0.2)
\psdots[dotstyle=*,linecolor=qqwuqq](2.29,-0.21)
\psdots[dotstyle=*,linecolor=qqwuqq](2.28,-0.22)
\psdots[dotstyle=*,linecolor=qqwuqq](2.27,-0.23)
\psdots[dotstyle=*,linecolor=qqwuqq](2.25,-0.24)
\psdots[dotstyle=*,linecolor=qqwuqq](2.24,-0.25)
\psdots[dotstyle=*,linecolor=qqwuqq](2.22,-0.26)
\psdots[dotstyle=*,linecolor=qqwuqq](2.2,-0.27)
\psdots[dotstyle=*,linecolor=qqwuqq](2.18,-0.29)
\psdots[dotstyle=*,linecolor=qqwuqq](2.16,-0.3)
\psdots[dotstyle=*,linecolor=qqwuqq](2.14,-0.32)
\psdots[dotstyle=*,linecolor=qqwuqq](2.12,-0.33)
\psdots[dotstyle=*,linecolor=qqwuqq](2.1,-0.35)
\psdots[dotstyle=*,linecolor=qqwuqq](2.08,-0.36)
\psdots[dotstyle=*,linecolor=qqwuqq](2.05,-0.38)
\psdots[dotstyle=*,linecolor=qqwuqq](2.03,-0.4)
\psdots[dotstyle=*,linecolor=qqwuqq](2,-0.42)
\psdots[dotstyle=*,linecolor=qqwuqq](1.97,-0.44)
\psdots[dotstyle=*,linecolor=qqwuqq](1.95,-0.46)
\psdots[dotstyle=*,linecolor=qqwuqq](1.92,-0.48)
\psdots[dotstyle=*,linecolor=qqwuqq](1.88,-0.5)
\psdots[dotstyle=*,linecolor=qqwuqq](1.85,-0.53)
\psdots[dotstyle=*,linecolor=qqwuqq](1.81,-0.55)
\psdots[dotstyle=*,linecolor=qqwuqq](1.78,-0.58)
\psdots[dotstyle=*,linecolor=qqwuqq](1.74,-0.61)
\psdots[dotstyle=*,linecolor=qqwuqq](1.7,-0.64)
\psdots[dotstyle=*,linecolor=qqwuqq](1.66,-0.67)
\psdots[dotstyle=*,linecolor=qqwuqq](1.61,-0.7)
\psdots[dotstyle=*,linecolor=qqwuqq](1.57,-0.73)
\psdots[dotstyle=*,linecolor=qqwuqq](1.52,-0.77)
\psdots[dotstyle=*,linecolor=qqwuqq](1.47,-0.81)
\psdots[dotstyle=*,linecolor=qqwuqq](1.41,-0.84)
\psdots[dotstyle=*,linecolor=qqwuqq](1.36,-0.89)
\psdots[dotstyle=*,linecolor=qqwuqq](1.3,-0.93)
\psdots[dotstyle=*,linecolor=qqwuqq](1.24,-0.97)
\psdots[dotstyle=*,linecolor=qqwuqq](1.17,-1.02)
\psdots[dotstyle=*,linecolor=qqwuqq](1.1,-1.07)
\psdots[dotstyle=*,linecolor=qqwuqq](1.03,-1.12)
\psdots[dotstyle=*,linecolor=qqwuqq](0.96,-1.17)
\psdots[dotstyle=*,linecolor=qqwuqq](0.88,-1.23)
\psdots[dotstyle=*,linecolor=qqwuqq](0.8,-1.29)
\psdots[dotstyle=*,linecolor=qqwuqq](0.71,-1.35)
\psdots[dotstyle=*,linecolor=qqwuqq](0.63,-1.41)
\psdots[dotstyle=*,linecolor=qqwuqq](0.53,-1.48)
\psdots[dotstyle=*,linecolor=qqwuqq](0.43,-1.55)
\psdots[dotstyle=*,linecolor=qqwuqq](0.33,-1.63)
\psdots[dotstyle=*,linecolor=qqwuqq](0.22,-1.7)
\psdots[dotstyle=*,linecolor=qqwuqq](0.11,-1.79)
\psdots[dotstyle=*,linecolor=qqwuqq](-0.01,-1.87)
\psdots[dotstyle=*,linecolor=qqwuqq](-0.13,-1.96)
\psdots[dotstyle=*,linecolor=qqwuqq](-0.26,-2.06)
\psdots[dotstyle=*,linecolor=qqwuqq](-0.4,-2.15)
\psdots[dotstyle=*,linecolor=qqwuqq](-0.54,-2.26)
\psdots[dotstyle=*,linecolor=qqwuqq](-0.69,-2.36)
\psdots[dotstyle=*,linecolor=qqwuqq](-0.85,-2.48)
\psdots[dotstyle=*,linecolor=qqwuqq](-1.01,-2.6)
\psdots[dotstyle=*,linecolor=qqwuqq](-1.19,-2.72)
\psdots[dotstyle=*,linecolor=qqwuqq](-1.37,-2.85)
\psdots[dotstyle=*,linecolor=qqwuqq](-1.56,-2.99)
\psdots[dotstyle=*,linecolor=qqwuqq](-1.75,-3.13)
\psdots[dotstyle=*,linecolor=qqwuqq](-1.96,-3.28)
\psdots[dotstyle=*,linecolor=qqwuqq](-2.18,-3.44)
\psdots[dotstyle=*,linecolor=qqwuqq](-2.41,-3.6)
\psdots[dotstyle=*,linecolor=qqwuqq](-2.65,-3.78)
\psdots[dotstyle=*,linecolor=qqwuqq](-2.9,-3.96)
\psdots[dotstyle=*,linecolor=qqwuqq](-3.16,-4.15)
\psdots[dotstyle=*,linecolor=qqwuqq](-3.43,-4.35)
\psdots[dotstyle=*,linecolor=qqwuqq](-3.72,-4.55)
\psdots[dotstyle=*,linecolor=qqwuqq](-4.02,-4.77)
\psdots[dotstyle=*,linecolor=qqwuqq](-4.34,-5)
\psdots[dotstyle=*,linecolor=qqwuqq](-4.67,-5.24)
\psdots[dotstyle=*,linecolor=qqwuqq](-5.02,-5.49)
\psdots[dotstyle=*,linecolor=qqwuqq](-5.38,-5.75)
\psdots[dotstyle=*,linecolor=qqwuqq](-5.76,-6.03)
\psdots[dotstyle=*,linecolor=qqwuqq](-6.16,-6.32)
\psdots[dotstyle=*,linecolor=qqwuqq](-6.58,-6.62)
\psdots[dotstyle=*,linecolor=qqwuqq](-7.02,-6.94)
\psdots[dotstyle=*,linecolor=qqwuqq](-7.48,-7.27)
\psdots[dotstyle=*,linecolor=qqwuqq](-7.96,-7.62)
\psdots[dotstyle=*,linecolor=qqwuqq](-8.47,-7.98)
\psdots[dotstyle=*,linecolor=qqwuqq](-9,-8.37)
\psdots[dotstyle=*,linecolor=qqwuqq](-9.55,-8.77)
\psdots[dotstyle=*,linecolor=qqwuqq](-10.14,-9.19)
\psdots[dotstyle=*,linecolor=qqwuqq](-10.75,-9.63)
\psdots[dotstyle=*,linecolor=qqwuqq](-11.38,-10.09)
\psdots[dotstyle=*,linecolor=qqwuqq](-12.05,-10.57)
\psdots[dotstyle=*,linecolor=qqwuqq](-12.75,-11.08)
\psdots[dotstyle=*,linecolor=qqwuqq](-13.49,-11.61)
\psdots[dotstyle=*,linecolor=qqwuqq](-14.26,-12.16)
\psdots[dotstyle=*,linecolor=qqwuqq](-15.06,-12.75)
\psdots[dotstyle=*,linecolor=qqwuqq](-15.91,-13.36)
\psdots[dotstyle=*,linecolor=qqwuqq](-16.8,-14)
\end{scriptsize}
\end{pspicture*}
\end{center}

Using the construction of the powers $X^n$ and $Z^m$ as above, we may
complete the image and  construct their products $X^n Z^m$ (with respect to some
origin $Y$): 

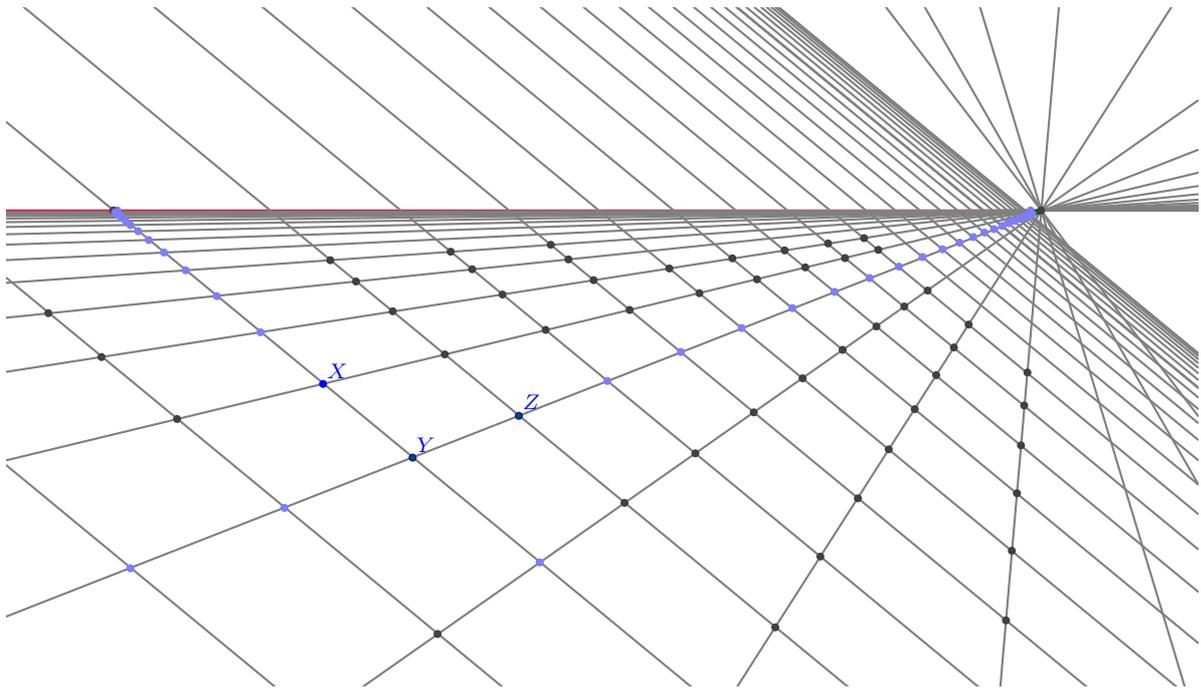
\begin{center}
\newrgbcolor{xdxdff}{0.49 0.49 1}
\newrgbcolor{dcrutc}{0.86 0.08 0.24}
\newrgbcolor{qqwuqq}{0 0.39 0}
\psset{xunit=0.6cm,yunit=0.6cm,algebraic=true,dotstyle=o,dotsize=3pt 0,linewidth=0.8pt,arrowsize=3pt 2,arrowinset=0.25}
\begin{pspicture*}(-12.44,-10.5)(13.7,4.48)
\psplot[linecolor=dcrutc]{-12.44}{13.7}{(-0-0*x)/20.02}
\psplot[linecolor=gray]{-12.44}{13.7}{(-16.43-1.63*x)/1.96}
\psplot[linecolor=gray]{-12.44}{13.7}{(-10.83-1.63*x)/1.96}
\psplot[linecolor=gray]{-12.44}{13.7}{(-6.17-1.63*x)/1.96}
\psplot[linecolor=gray]{-12.44}{13.7}{(-2.3-1.63*x)/1.96}
\psplot[linecolor=gray]{-12.44}{13.7}{(--0.91-1.63*x)/1.96}
\psplot[linecolor=gray]{-12.44}{13.7}{(--3.58-1.63*x)/1.96}
\psplot[linecolor=gray]{-12.44}{13.7}{(--5.8-1.63*x)/1.96}
\psplot[linecolor=gray]{-12.44}{13.7}{(--7.64-1.63*x)/1.96}
\psplot[linecolor=gray]{-12.44}{13.7}{(--9.17-1.63*x)/1.96}
\psplot[linecolor=gray]{-12.44}{13.7}{(--10.45-1.63*x)/1.96}
\psplot[linecolor=gray]{-12.44}{13.7}{(--11.5-1.63*x)/1.96}
\psplot[linecolor=gray]{-12.44}{13.7}{(--12.38-1.63*x)/1.96}
\psplot[linecolor=gray]{-12.44}{13.7}{(--13.11-1.63*x)/1.96}
\psplot[linecolor=gray]{-12.44}{13.7}{(--13.72-1.63*x)/1.96}
\psplot[linecolor=gray]{-12.44}{13.7}{(--14.22-1.63*x)/1.96}
\psplot[linecolor=gray]{-12.44}{13.7}{(--14.64-1.63*x)/1.96}
\psplot[linecolor=gray]{-12.44}{13.7}{(--14.64-1.63*x)/1.96}
\psplot[linecolor=gray]{-12.44}{13.7}{(-9.47--0.92*x)/2.33}
\psplot[linecolor=gray]{-12.44}{13.7}{(-27.57--2.69*x)/17.11}
\psplot[linecolor=gray]{-12.44}{13.7}{(-19.35--1.89*x)/18.08}
\psplot[linecolor=gray]{-12.44}{13.7}{(-13.58--1.32*x)/18.75}
\psplot[linecolor=gray]{-12.44}{13.7}{(-9.53--0.93*x)/19.23}
\psplot[linecolor=gray]{-12.44}{13.7}{(-6.69--0.65*x)/19.56}
\psplot[linecolor=gray]{-12.44}{13.7}{(-4.69--0.46*x)/19.8}
\psplot[linecolor=gray]{-12.44}{13.7}{(-3.3--0.32*x)/19.96}
\psplot[linecolor=gray]{-12.44}{13.7}{(-2.31--0.23*x)/20.08}
\psplot[linecolor=gray]{-12.44}{13.7}{(-1.62--0.16*x)/20.16}
\psplot[linecolor=gray]{-12.44}{13.7}{(--1.14-0.11*x)/-20.21}
\psplot[linecolor=gray]{-12.44}{13.7}{(-39.29--3.83*x)/15.74}
\psplot[linecolor=gray]{-12.44}{13.7}{(--0.8-0.08*x)/-20.25}
\psplot[linecolor=gray]{-12.44}{13.7}{(--0.56-0.05*x)/-20.28}
\psplot[linecolor=gray]{-12.44}{13.7}{(--0.39-0.04*x)/-20.3}
\psplot[linecolor=gray]{-12.44}{13.7}{(-23.18-1.63*x)/1.96}
\psplot[linecolor=gray]{-12.44}{13.7}{(-31.3-1.63*x)/1.96}
\psplot[linecolor=gray]{-12.44}{13.7}{(-41.08-1.63*x)/1.96}
\psplot[linecolor=gray]{-12.44}{13.7}{(-52.85-1.63*x)/1.96}
\psplot[linecolor=gray]{-12.44}{13.7}{(-67.01-1.63*x)/1.96}
\psplot[linecolor=gray]{-12.44}{13.7}{(-84.06-1.63*x)/1.96}
\psplot[linecolor=gray]{-12.44}{13.7}{(-104.59-1.63*x)/1.96}
\psplot[linecolor=gray]{-12.44}{13.7}{(--79.76-7.78*x)/-10.99}
\psplot[linecolor=gray]{-12.44}{13.7}{(--113.64-11.08*x)/-7.02}
\psplot[linecolor=gray]{-12.44}{13.7}{(--161.92-15.79*x)/-1.35}
\psplot[linecolor=gray]{-12.44}{13.7}{(--230.71-22.5*x)/6.72}
\psplot[linecolor=gray]{-12.44}{13.7}{(--328.73-32.06*x)/18.22}
\psplot[linecolor=gray]{-12.44}{13.7}{(--468.39-45.69*x)/34.6}
\psplot[linecolor=gray]{-12.44}{13.7}{(--14.98-1.63*x)/1.96}
\psplot[linecolor=gray]{-12.44}{13.7}{(--15.27-1.63*x)/1.96}
\psplot[linecolor=gray]{-12.44}{13.7}{(--15.51-1.63*x)/1.96}
\psplot[linecolor=gray]{-12.44}{13.7}{(--15.71-1.63*x)/1.96}
\psplot[linecolor=gray]{-12.44}{13.7}{(--14.64-1.63*x)/1.96}
\psplot[linecolor=gray]{-12.44}{13.7}{(--15.88-1.63*x)/1.96}
\psplot[linecolor=gray]{-12.44}{13.7}{(--16.02-1.63*x)/1.96}
\psplot[linecolor=gray]{-12.44}{13.7}{(--16.13-1.63*x)/1.96}
\psplot[linecolor=gray]{-12.44}{13.7}{(--16.13-1.63*x)/1.96}
\begin{scriptsize}
\psdots[dotstyle=*,linecolor=blue](-3.53,-5.46)
\rput[bl](-3.44,-5.33){\blue{$Y$}}
\psdots[dotstyle=*,linecolor=blue](-5.49,-3.83)
\rput[bl](-5.4,-3.7){\blue{$X$}}
\psdots[dotstyle=*,linecolor=xdxdff](-10,0)
\psdots[dotstyle=*,linecolor=xdxdff](10.02,0)
\psdots[dotsize=2pt 0,dotstyle=*,linecolor=qqwuqq](-10.09,0)
\psdots[dotsize=2pt 0,dotstyle=*,linecolor=qqwuqq](-10.09,0)
\psdots[dotsize=2pt 0,dotstyle=*,linecolor=qqwuqq](-10.09,0)
\psdots[dotsize=2pt 0,dotstyle=*,linecolor=qqwuqq](-10.09,0)
\psdots[dotsize=2pt 0,dotstyle=*,linecolor=qqwuqq](-10.09,0)
\psdots[dotsize=2pt 0,dotstyle=*,linecolor=qqwuqq](-10.09,0)
\psdots[dotsize=2pt 0,dotstyle=*,linecolor=qqwuqq](-10.09,0)
\psdots[dotsize=2pt 0,dotstyle=*,linecolor=qqwuqq](-10.09,0)
\psdots[dotsize=2pt 0,dotstyle=*,linecolor=qqwuqq](-10.09,0)
\psdots[dotsize=2pt 0,dotstyle=*,linecolor=qqwuqq](-10.09,0)
\psdots[dotsize=2pt 0,dotstyle=*,linecolor=qqwuqq](-10.09,0)
\psdots[dotsize=2pt 0,dotstyle=*,linecolor=qqwuqq](-10.09,0)
\psdots[dotsize=2pt 0,dotstyle=*,linecolor=qqwuqq](-10.09,0)
\psdots[dotsize=2pt 0,dotstyle=*,linecolor=qqwuqq](-10.09,0)
\psdots[dotsize=2pt 0,dotstyle=*,linecolor=qqwuqq](-10.09,0)
\psdots[dotsize=2pt 0,dotstyle=*,linecolor=qqwuqq](-10.09,0)
\psdots[dotsize=2pt 0,dotstyle=*,linecolor=qqwuqq](-10.09,0)
\psdots[dotsize=2pt 0,dotstyle=*,linecolor=qqwuqq](-10.09,0)
\psdots[dotsize=2pt 0,dotstyle=*,linecolor=qqwuqq](-10.09,0)
\psdots[dotsize=2pt 0,dotstyle=*,linecolor=qqwuqq](-10.09,0)
\psdots[dotsize=2pt 0,dotstyle=*,linecolor=qqwuqq](-10.09,0)
\psdots[dotsize=2pt 0,dotstyle=*,linecolor=qqwuqq](-10.09,0)
\psdots[dotsize=2pt 0,dotstyle=*,linecolor=qqwuqq](-10.09,0)
\psdots[dotsize=2pt 0,dotstyle=*,linecolor=qqwuqq](-10.09,0)
\psdots[dotsize=2pt 0,dotstyle=*,linecolor=qqwuqq](-10.09,0)
\psdots[dotsize=2pt 0,dotstyle=*,linecolor=qqwuqq](-10.09,0)
\psdots[dotsize=2pt 0,dotstyle=*,linecolor=qqwuqq](-10.09,0)
\psdots[dotsize=2pt 0,dotstyle=*,linecolor=qqwuqq](-10.09,0)
\psdots[dotsize=2pt 0,dotstyle=*,linecolor=qqwuqq](-10.09,0)
\psdots[dotsize=2pt 0,dotstyle=*,linecolor=qqwuqq](-10.09,0)
\psdots[dotsize=2pt 0,dotstyle=*,linecolor=qqwuqq](-10.09,0)
\psdots[dotsize=2pt 0,dotstyle=*,linecolor=qqwuqq](-10.09,0)
\psdots[dotstyle=*,linecolor=blue](-1.2,-4.54)
\rput[bl](-1.11,-4.4){\blue{$Z$}}
\psdots[dotsize=2pt 0,dotstyle=*,linecolor=qqwuqq](-18.68,-11.46)
\psdots[dotsize=2pt 0,dotstyle=*,linecolor=qqwuqq](-13.78,-9.52)
\psdots[dotsize=2pt 0,dotstyle=*,linecolor=qqwuqq](-9.71,-7.91)
\psdots[dotsize=2pt 0,dotstyle=*,linecolor=qqwuqq](-6.34,-6.57)
\psdots[dotsize=2pt 0,dotstyle=*,linecolor=qqwuqq](-3.53,-5.46)
\psdots[dotsize=2pt 0,dotstyle=*,linecolor=qqwuqq](-1.2,-4.54)
\psdots[dotsize=2pt 0,dotstyle=*,linecolor=qqwuqq](0.74,-3.77)
\psdots[dotsize=2pt 0,dotstyle=*,linecolor=qqwuqq](2.35,-3.13)
\psdots[dotsize=2pt 0,dotstyle=*,linecolor=qqwuqq](3.69,-2.6)
\psdots[dotsize=2pt 0,dotstyle=*,linecolor=qqwuqq](4.8,-2.16)
\psdots[dotsize=2pt 0,dotstyle=*,linecolor=qqwuqq](5.72,-1.8)
\psdots[dotsize=2pt 0,dotstyle=*,linecolor=qqwuqq](6.49,-1.49)
\psdots[dotsize=2pt 0,dotstyle=*,linecolor=qqwuqq](7.13,-1.24)
\psdots[dotsize=2pt 0,dotstyle=*,linecolor=qqwuqq](7.65,-1.03)
\psdots[dotsize=2pt 0,dotstyle=*,linecolor=qqwuqq](8.09,-0.86)
\psdots[dotsize=2pt 0,dotstyle=*,linecolor=qqwuqq](8.46,-0.71)
\psdots[dotsize=2pt 0,dotstyle=*,linecolor=qqwuqq](8.76,-0.59)
\psdots[dotsize=2pt 0,dotstyle=*,linecolor=qqwuqq](9.01,-0.49)
\psdots[dotsize=2pt 0,dotstyle=*,linecolor=qqwuqq](9.22,-0.41)
\psdots[dotsize=2pt 0,dotstyle=*,linecolor=qqwuqq](9.4,-0.34)
\psdots[dotsize=2pt 0,dotstyle=*,linecolor=qqwuqq](9.54,-0.28)
\psdots[dotsize=2pt 0,dotstyle=*,linecolor=qqwuqq](9.66,-0.23)
\psdots[dotsize=2pt 0,dotstyle=*,linecolor=qqwuqq](9.76,-0.19)
\psdots[dotsize=2pt 0,dotstyle=*,linecolor=qqwuqq](9.85,-0.16)
\psdots[dotsize=2pt 0,dotstyle=*,linecolor=qqwuqq](9.91,-0.13)
\psdots[dotsize=2pt 0,dotstyle=*,linecolor=qqwuqq](9.97,-0.11)
\psdots[dotsize=2pt 0,dotstyle=*,linecolor=qqwuqq](10.02,-0.09)
\psdots[dotsize=2pt 0,dotstyle=*,linecolor=qqwuqq](10.06,-0.08)
\psdots[dotsize=2pt 0,dotstyle=*,linecolor=qqwuqq](10.09,-0.06)
\psdots[dotsize=2pt 0,dotstyle=*,linecolor=qqwuqq](10.12,-0.05)
\psdots[dotsize=2pt 0,dotstyle=*,linecolor=qqwuqq](10.14,-0.04)
\psdots[dotsize=2pt 0,dotstyle=*,linecolor=qqwuqq](10.16,-0.04)
\psdots[dotsize=2pt 0,dotstyle=*,linecolor=qqwuqq](10.18,-0.03)
\psdots[dotsize=2pt 0,dotstyle=*,linecolor=qqwuqq](10.19,-0.03)
\psdots[dotsize=2pt 0,dotstyle=*,linecolor=qqwuqq](10.2,-0.02)
\psdots[dotsize=2pt 0,dotstyle=*,linecolor=qqwuqq](10.21,-0.02)
\psdots[dotsize=2pt 0,dotstyle=*,linecolor=qqwuqq](10.22,-0.01)
\psdots[dotsize=2pt 0,dotstyle=*,linecolor=qqwuqq](10.22,-0.01)
\psdots[dotsize=2pt 0,dotstyle=*,linecolor=qqwuqq](10.23,-0.01)
\psdots[dotsize=2pt 0,dotstyle=*,linecolor=qqwuqq](10.23,-0.01)
\psdots[dotsize=2pt 0,dotstyle=*,linecolor=qqwuqq](10.23,-0.01)
\psdots[dotsize=2pt 0,dotstyle=*,linecolor=qqwuqq](10.24,-0.01)
\psdots[dotsize=2pt 0,dotstyle=*,linecolor=qqwuqq](10.24,0)
\psdots[dotsize=2pt 0,dotstyle=*,linecolor=qqwuqq](10.24,0)
\psdots[dotsize=2pt 0,dotstyle=*,linecolor=qqwuqq](10.24,0)
\psdots[dotsize=2pt 0,dotstyle=*,linecolor=qqwuqq](10.25,0)
\psdots[dotsize=2pt 0,dotstyle=*,linecolor=qqwuqq](10.25,0)
\psdots[dotsize=2pt 0,dotstyle=*,linecolor=qqwuqq](10.25,0)
\psdots[dotsize=2pt 0,dotstyle=*,linecolor=qqwuqq](10.25,0)
\psdots[dotsize=2pt 0,dotstyle=*,linecolor=qqwuqq](10.25,0)
\psdots[dotsize=2pt 0,dotstyle=*,linecolor=qqwuqq](10.25,0)
\psdots[dotsize=2pt 0,dotstyle=*,linecolor=qqwuqq](10.25,0)
\psdots[dotsize=2pt 0,dotstyle=*,linecolor=qqwuqq](10.25,0)
\psdots[dotsize=2pt 0,dotstyle=*,linecolor=qqwuqq](10.25,0)
\psdots[dotsize=2pt 0,dotstyle=*,linecolor=qqwuqq](10.25,0)
\psdots[dotsize=2pt 0,dotstyle=*,linecolor=qqwuqq](10.25,0)
\psdots[dotstyle=*,linecolor=darkgray](-10.09,0)
\psdots[dotstyle=*,linecolor=xdxdff](0.74,-3.77)
\psdots[dotstyle=*,linecolor=xdxdff](2.35,-3.13)
\psdots[dotstyle=*,linecolor=xdxdff](3.69,-2.6)
\psdots[dotstyle=*,linecolor=xdxdff](4.8,-2.16)
\psdots[dotstyle=*,linecolor=xdxdff](5.72,-1.8)
\psdots[dotstyle=*,linecolor=xdxdff](6.49,-1.49)
\psdots[dotstyle=*,linecolor=xdxdff](7.13,-1.24)
\psdots[dotstyle=*,linecolor=xdxdff](7.65,-1.03)
\psdots[dotstyle=*,linecolor=xdxdff](8.09,-0.86)
\psdots[dotstyle=*,linecolor=xdxdff](8.46,-0.71)
\psdots[dotstyle=*,linecolor=xdxdff](8.76,-0.59)
\psdots[dotstyle=*,linecolor=xdxdff](9.01,-0.49)
\psdots[dotstyle=*,linecolor=xdxdff](9.22,-0.41)
\psdots[dotstyle=*,linecolor=xdxdff](9.4,-0.34)
\psdots[dotstyle=*,linecolor=darkgray](10.25,0)
\psdots[dotstyle=*,linecolor=xdxdff](-6.86,-2.69)
\psdots[dotstyle=*,linecolor=xdxdff](-7.82,-1.89)
\psdots[dotstyle=*,linecolor=xdxdff](-8.5,-1.32)
\psdots[dotstyle=*,linecolor=xdxdff](-8.98,-0.93)
\psdots[dotstyle=*,linecolor=xdxdff](-9.31,-0.65)
\psdots[dotstyle=*,linecolor=xdxdff](-9.54,-0.46)
\psdots[dotstyle=*,linecolor=xdxdff](-9.71,-0.32)
\psdots[dotstyle=*,linecolor=xdxdff](-9.82,-0.23)
\psdots[dotstyle=*,linecolor=xdxdff](-9.9,-0.16)
\psdots[dotstyle=*,linecolor=xdxdff](-9.96,-0.11)
\psdots[dotstyle=*,linecolor=xdxdff](-10,-0.08)
\psdots[dotstyle=*,linecolor=xdxdff](-10.03,-0.05)
\psdots[dotstyle=*,linecolor=xdxdff](-10.05,-0.04)
\psdots[dotstyle=*,linecolor=xdxdff](-6.34,-6.57)
\psdots[dotstyle=*,linecolor=xdxdff](-9.71,-7.91)
\psdots[dotstyle=*,linecolor=xdxdff](-13.78,-9.52)
\psdots[dotstyle=*,linecolor=xdxdff](-18.68,-11.46)
\psdots[dotstyle=*,linecolor=xdxdff](-24.57,-13.8)
\psdots[dotstyle=*,linecolor=xdxdff](-31.66,-16.61)
\psdots[dotstyle=*,linecolor=xdxdff](-40.2,-19.99)
\psdots[dotstyle=*,linecolor=xdxdff](-0.74,-7.78)
\psdots[dotstyle=*,linecolor=xdxdff](3.24,-11.08)
\psdots[dotstyle=*,linecolor=xdxdff](8.9,-15.79)
\psdots[dotstyle=*,linecolor=xdxdff](16.97,-22.5)
\psdots[dotstyle=*,linecolor=xdxdff](28.47,-32.06)
\psdots[dotstyle=*,linecolor=xdxdff](44.85,-45.69)
\psdots[dotstyle=*,linecolor=xdxdff](9.54,-0.28)
\psdots[dotstyle=*,linecolor=xdxdff](9.66,-0.23)
\psdots[dotstyle=*,linecolor=xdxdff](9.76,-0.19)
\psdots[dotstyle=*,linecolor=xdxdff](9.85,-0.16)
\psdots[dotstyle=*,linecolor=xdxdff](9.91,-0.13)
\psdots[dotstyle=*,linecolor=xdxdff](9.97,-0.11)
\psdots[dotstyle=*,linecolor=xdxdff](10.02,-0.09)
\psdots[dotstyle=*,linecolor=xdxdff](10.06,-0.08)
\psdots[dotstyle=*,linecolor=darkgray](-2.82,-3.18)
\psdots[dotstyle=*,linecolor=darkgray](-0.61,-2.64)
\psdots[dotstyle=*,linecolor=darkgray](1.23,-2.2)
\psdots[dotstyle=*,linecolor=darkgray](2.76,-1.83)
\psdots[dotstyle=*,linecolor=darkgray](4.02,-1.52)
\psdots[dotstyle=*,linecolor=darkgray](5.08,-1.26)
\psdots[dotstyle=*,linecolor=darkgray](5.95,-1.05)
\psdots[dotstyle=*,linecolor=darkgray](6.68,-0.87)
\psdots[dotstyle=*,linecolor=darkgray](-3.96,-2.23)
\psdots[dotstyle=*,linecolor=darkgray](-1.56,-1.86)
\psdots[dotstyle=*,linecolor=darkgray](0.44,-1.54)
\psdots[dotstyle=*,linecolor=darkgray](2.1,-1.28)
\psdots[dotstyle=*,linecolor=darkgray](3.48,-1.06)
\psdots[dotstyle=*,linecolor=darkgray](4.63,-0.88)
\psdots[dotstyle=*,linecolor=darkgray](5.58,-0.73)
\psdots[dotstyle=*,linecolor=darkgray](6.37,-0.61)
\psdots[dotstyle=*,linecolor=darkgray](-4.77,-1.57)
\psdots[dotstyle=*,linecolor=darkgray](-2.22,-1.3)
\psdots[dotstyle=*,linecolor=darkgray](-0.11,-1.08)
\psdots[dotstyle=*,linecolor=darkgray](1.12,-6.46)
\psdots[dotstyle=*,linecolor=darkgray](2.67,-5.37)
\psdots[dotstyle=*,linecolor=darkgray](3.95,-4.46)
\psdots[dotstyle=*,linecolor=darkgray](5.9,-3.08)
\psdots[dotstyle=*,linecolor=darkgray](5.02,-3.71)
\psdots[dotstyle=*,linecolor=darkgray](6.64,-2.56)
\psdots[dotstyle=*,linecolor=darkgray](7.25,-2.12)
\psdots[dotstyle=*,linecolor=darkgray](7.76,-1.77)
\psdots[dotstyle=*,linecolor=darkgray](8.66,-2.52)
\psdots[dotstyle=*,linecolor=darkgray](8.34,-3.03)
\psdots[dotstyle=*,linecolor=darkgray](7.95,-3.64)
\psdots[dotstyle=*,linecolor=darkgray](7.48,-4.39)
\psdots[dotstyle=*,linecolor=darkgray](6.91,-5.28)
\psdots[dotstyle=*,linecolor=darkgray](6.23,-6.36)
\psdots[dotstyle=*,linecolor=darkgray](5.41,-7.65)
\psdots[dotstyle=*,linecolor=darkgray](4.42,-9.21)
\psdots[dotstyle=*,linecolor=darkgray](-2.98,-9.36)
\psdots[dotstyle=*,linecolor=darkgray](-8.69,-4.61)
\psdots[dotstyle=*,linecolor=darkgray](-10.35,-3.24)
\psdots[dotstyle=*,linecolor=darkgray](-11.51,-2.27)
\psdots[dotstyle=*,linecolor=darkgray](-14.54,-3.9)
\psdots[dotstyle=*,linecolor=darkgray](-12.55,-5.55)
\psdots[dotstyle=*,linecolor=darkgray](-5.67,-11.27)
\psdots[dotstyle=*,linecolor=darkgray](1.81,-13.34)
\psdots[dotstyle=*,linecolor=darkgray](9.13,-13.12)
\psdots[dotstyle=*,linecolor=darkgray](9.32,-10.9)
\psdots[dotstyle=*,linecolor=darkgray](9.48,-9.06)
\psdots[dotstyle=*,linecolor=darkgray](9.61,-7.52)
\psdots[dotstyle=*,linecolor=darkgray](9.72,-6.25)
\psdots[dotstyle=*,linecolor=darkgray](9.81,-5.19)
\psdots[dotstyle=*,linecolor=darkgray](9.88,-4.31)
\psdots[dotstyle=*,linecolor=darkgray](9.95,-3.58)
\psdots[dotstyle=*,linecolor=darkgray](-5.33,-1.1)
\psdots[dotstyle=*,linecolor=darkgray](-2.69,-0.91)
\psdots[dotstyle=*,linecolor=darkgray](-0.5,-0.76)
\end{scriptsize}
\end{pspicture*}
\end{center}


\subsubsection{Boundary values}\label{subsec:ComputeBoundary}

The line $a$ does not belong to $G$: we may consider it as a sort of {\em boundary} of $G$.
It is clear that the group $G$ acts on itself from the left and from the right.
Show that $G$ also acts from the left and from the right on the boundary $a$, just by
extending the domain of definition of our main formula (\ref{eqn:Main}):

\begin{description}
\item{(1)} show that, if $z \in a$ and $x,y \in G$, we have $(xyz)=z$, and conclude that
$G$ acts trivially on $a$ from the left,
\item{(2)} show that, if $x \in a$ and $y,z \in G$, we have $\alpha((xyz))=0$, so 
 $w$ lies again on $a$.
\end{description}

Show that (2) defines an action of $G$ on $a$ from the right: 
fix an origin $e \in G$ with $\alpha(e)=1$, and show that  
the group $(G,e)$ acts on the line $a$ from the right by the formula
$$
a \times G  \to a, \quad (x,z) \mapsto (xez) .
$$
Moreover, show
that the map $\rho(z):a \to a$, $x \mapsto (xez)$ is affine, and that every element
of the group $\Ga(1,\R)$ of affine bijectijons of $\R$ is obtained in this way.
Show that the map
$$
G \to \Ga(1,\R), \quad z \mapsto \rho(z)
$$
is an isomorphism of groups.

\ssk
From a slightly different viewpoint, show that $(G,e)$ acts on $E$ from the left and from
the right, that both of these actions commute, and determine the orbits of both actions.

\subsubsection{Connected components}

The line $a$ divides the set $G$ into two half-planes which we call the two
{\em connected components of $G$}.
Let $e \in G$; show that the connected component containing $e$ is a subgroup
of $(G,e)$, and that the other connected component is not a subgroup, but both components are
subtorsors of $G$.

\subsection{Generalization I: affine spaces, a hyperplane taken out}
\label{subsec:GeneralizationI}

The constructions described above admit far-reaching generalizations.
Unfortunately, it is no longer possible to illustrate them in a such a direct way since
the most interesting new features need, to emerge,  at least a four dimensional surrounding (the
space of $2 \times 2$-matrices). In this sense experimental geometry has reached its limits --
however, this does not mean that geometry dissappears from the picture, see the last chapter.

\ssk
The following is an immediate generalization of the preceding exercises. 
Let $E = \R^n$ and $a$ a hyperplane of $E$, which we write in the form
$a=\ker(\alpha)$, for some non-zero linear form $\alpha:E \to \R$.
Show that formula (\ref{eqn:Main}) defines a torsor structure on the set
$G = E \setminus a$, and that $\alpha:G \to \R^\times $ is a torsor homomorphism.
Show that, after choice of an origin $e \in G$, the group $(G,e)$ is isomorphic
to a semidirect product $\R^{n-1} \rtimes \R^\times$.
This group has the structure of the {\em group of dilations} of $\R^{n+1}$ (sugroup
of the affine group of $\R^{n-1}$ generated by all dilations and translations).
A geometric way to identify $G$ with this group is by letting act $G$
on the ``boundary'' $a$ of $G$, as above. 
(The only difference with the case $n=2$ is that, for $n>2$, the dilation group
is strictly smaller than the full affine group, and that it is more difficult to make drawings.)

\msk
Show that formula (\ref{eqn:Main}) is related to incidence geometry in the same way
as for $n=2$: just replace the term ``line'' by ``hyperplane''.  
As a consquence, show that we can project the $n+1$-dimensional space homomorphically
onto the $n$-dimensional one:
let $E_{n+1}=\R^{n+1}$ with the linear form $\alpha(x)=x_1$; show that then the map
$E_{n+1}\to E_n$, $(x_1,\ldots,x_{n+1}) \mapsto (x_1,\ldots,x_n)$ is a torsor-homomorphism.
Its kernel is a one-dimensional subgroup. 
For instance, if $n=2$, we get a homomorphic projection $E_3 \to E_2$, that is,
a two-dimensional homomorphic image of a three-dimensional situation.
This sheds some light on the impression that many of our drawings (and even more
their dynamic versions) seem to have a ``spacial interpretation''.

\msk
Finally, show that the field $\R$ may everywhere be replaced by an
arbitrary field $\K$.

\subsection{Generalization II: spaces of rectangular matrices}
\label{subsec:GeneralizationII}

The second generalization is most important:
take as space $E$ the space $M(p,q;\R)$ of $q \times p$-matrices, where $p,q \in \N$.
Fix some  $p\times q$-matrix $A$.
For  three $q \times p$-matrices $X,Y,Z \in E$, 
such that the $p \times p$-matrix $AY$ is invertible,
define another $q \times p$-matrix
\begin{equation}\label{eqn:Main!}
\boxed{
W :=  (XYZ)_A:=(X-Y) (AY)\inv AZ  + Z = X(AY)\inv AZ - Y (AY)\inv AZ  + Z } \,  . 
\end{equation}
(Convince yourself that this is well-defined, i.e., all matrix products are indeed defined!)
Show that the ternary law $(XYZ)_A$ defines a torsor structure on the set
\begin{equation}
\boxed{
G := \Gl(p,q;A,\R) := \big\{ X \in M(p,q;\R) \mid \, \det (AY) \not= 0 \big\} } \, .
\end{equation}
Hint: proceed as above (which is indeed the special case $p=1,q=n$):
start by showing that $A (XYZ)_A = (AX)(AY)\inv(AX)$, i.e., the map
\begin{equation}
\alpha: \Gl(p,q;A,\R) \to \Gl(p,\R), \quad X \mapsto AX
\end{equation}
is a homomorphism of ternary products. 

\msk
A detailed study of the possible cases would lead us too far away from the topics
we started with. Let us just mention that the most interesting and richest cases
arise for $p=q >1$: 
for instance, if $p=q=2$, we may choose for $A$ a matrix of rank 1 or 2.
If the rank of $A$ is 1, then the resulting geometry looks very much like the geometry of
$\R^4$ we mentioned in the preceding subsection; however, if the rank of $A$ is $2$, then
it is quite different: 
in this case the homomorphism $\alpha$ becomes an isomorphism, and hence our
torsor $\Gl(2,2;A,\R)$ is isomorphic to the general linear group $\Gl(2,\R)$.

\msk
We add some short comments for the more expert reader:
in the language of Lie group theory, we see that in the last mentioned case
 a {\em reductive} Lie group shows up,
whereas in the preceding subsection only solvable Lie groups appeared.
All other classical Lie groups can be realized in a similar way (but the constructions are
more complicated, using the Jordan algebras of
Hermitian or skew-Hermitian matrices),  see \cite{BeKi10, BeKi10b, BeKi12, Be12}.

\section{Geometry revisited: axiomatic approach}\label{sec:Geometry}

In this final section, we say some words on the relation between {\em geometry} and {\em algebra}.
The power of algebra consists, for instance, in the possibility of the vast generalizations 
mentioned in the preceding subsections.
The power of geometry manifests itself in another kind of generalization: the various
kinds of {\em non-Euclidean geometries}. 
In this section, which is far from being exhaustive, we would like to explain some 
of the challenges  arising in this context.


\subsection{Two-parallel planes}\label{subsec:Two-par}

The most primitive structure appearing in an axiomatic approach to plane geometry is
{\em incidence}: we have a set $\cP$ of {\em points} and a set $\cL$ of subsets of $\cP$
whose elements are called {\em lines}; a point $p$ and a line $\ell$ are {\em incident}
if $p \in \ell$. Two lines $\ell$ and $m$ are called {\em parallel} (notation: $m\parallel \ell$)
if they are equal or disjoint.
These data are called an {\em affine plane} if they satisfy the following axioms:

\begin{description}
\item{(A1)} Any two distinct points lie on a unique line.
\item{(A2)} For a point $p$ and  a line $m$, there is a unique line $\ell$ through $p$
and parallel to $m$.
\item{(A3)} There are three points which are not collinear.
\end{description}

Show that parallism is an {\em equivalence relation} on the set of lines, and show
 that, for any field $\K$, the plane $\K^2$ with its usual incidence structure is an
affine plane in the sense of the definition.
There are interesting affine planes that are {\em not} of this form; we come to this item
below. 
Non-Euclidean geometry arose from modifying some of Euclide's axioms; in particular,
it was questioned whether the uniqueness of parallels were a consequence of the other
axioms or not -- indeed, it turned out to be independent. In this subsection, we propose
some exercises about a simple model of a geometry in which there are always
{\em one or two parallels to a given line} through each point. 
Here are the axioms:

\msk
{\bf Definition.} {\em 
A \emph{two-parallel plane} is a set $U$ of ``points'' together with a set $\cL$ of subsets of $U$,
 called ``lines''; we assume that there are two kinds of lines, that is, we have
a partition $\cL = \cL_1 \dot\cup \cL_2$ into a set $\cL_1$ of {\em lines of the first kind} and
a set $\cL_2$ of {\em lines of the second kind},
such that:
\begin{description}
\item{(TP1)}
Any two distinct points lie on a unique line.
\item{(TP2)}
Given a line $\ell \in \cL_i$ ($i=1,2$) and a point $p$ not on $\ell$,
 there are exactly $i$ different lines through $p$ and  parallel to $\ell$.
\item{(TP3)}
(Consistency axiom.)
The set $P := \{ (\ell,\ell') \in \cL_2 \times \cL_2 \mid \ell \parallel \ell', \ \ell \not= \ell' \}$ admits a partition
$P = P_a  \dot\cup P_b$ (and call a pair of lines $(\ell,\ell')$ {\em $a$-parallel}, resp.\
{\em $b$-parallel},  if $(\ell,\ell') \in P_a$, resp.\ $(\ell,m) \in P_b$)
 such that:  for a point $p$ and a line $\ell \in \cL_2$,
there is exactly one $a$-parallel through $p$ and exactly one $b$-parallel through $p$.
\item{(TP4)}
There exist three non-collinear points (points not on a single line).
\end{description}
}

\msk
Show:
if $\cP$ is an affine plane and $a$ some line of this plane, then the set
$U:= \cP \setminus a$ is a two-parallel plane, where a {\em line in $U$} is the trace
(=intersection with $U$) of a line in $\cP$.
Hints:
a line $\ell$ is of the first kind if it is  (in $\cP$) parallel to $a$, and 
of the second kind else; $P_b$ is the set of pairs of lines that are in $\cP$ parallel, and
$P_a$ the set of pairs of lines that (in $\cP$) intersect on $a$ -- so in $U$ they are parallel, too!

\msk 
(Optional.) What about the converse: is every two-parallel plane obtained by the preceding
construction? Hints: if you know how to construct a projective plane out of an affine plane,
you may apply the same construction twice -- so the answer will be positive. 

\msk
In a two-parallel plane, we may mimick our fundamental definition 
(\ref{eqn:Generic}): for a generic triple of points $(x,y,z) \in U^3$, define a fourth point:

\begin{center}
$w:=(xyz):=$ intersection of [the $a$-parallel of $y \lor x$ through $z$] 

${ }\qquad$ with [the $b$-parallel of $y \lor z$ through $x$].
\end{center}

What can you say about the ternay law $(x,y,z) \mapsto w$:
does it satisfy idempotency and general associativity, or other algebraic laws?
And can we extend this definition in some ``nice way'' also to the case of collinear
triples $(x,y,z)$? 
These are difficult questions -- in fact, they are essentially equivalent to similar questions
from the theory of {\em projective planes} -- but the viewpoint presented here seems to be new. 
In the following, let us explain that the answer to the first question is related to 
{\em Desargues' theorem}.

\subsection{Associativity and Desargues theorem}

Preliminaries: inform yourself about {\em Desargues theorem},
either in its affine or projective version,
saying that
{\em two triangles are in perspective axially if and only if they are in perspective centrally},
and using geogebra, 
convince yourself experimentally of the validity of this theorem:


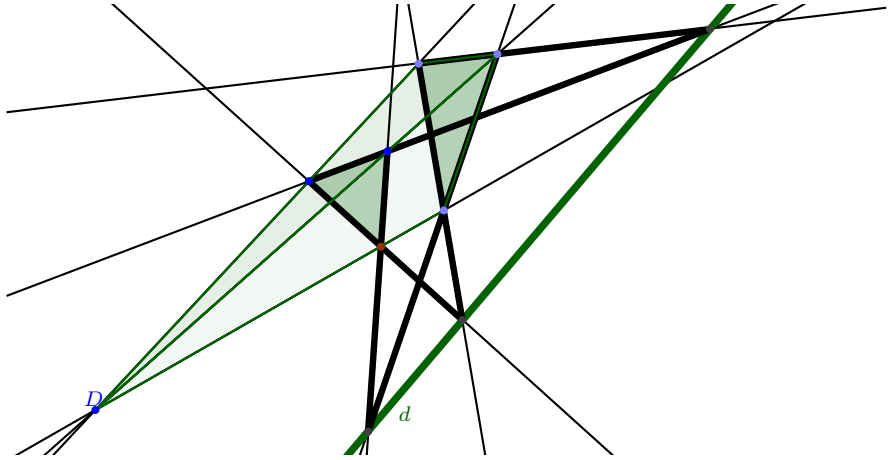
\begin{center}
\newrgbcolor{zzttqq}{0.6 0.2 0}
\newrgbcolor{xdxdff}{0.49 0.49 1}
\newrgbcolor{qqwuqq}{0 0.39 0}
\newrgbcolor{wqwqwq}{0.38 0.38 0.38}
\psset{xunit=0.5cm,yunit=0.4cm,algebraic=true,dotstyle=o,dotsize=3pt 0,linewidth=0.8pt,arrowsize=3pt 2,arrowinset=0.25}
\begin{pspicture*}(-7.11,-10.21)(16.03,4.73)
\pspolygon[linecolor=qqwuqq,fillcolor=qqwuqq,fillstyle=solid,opacity=0.25](0.84,-1.14)(2.91,-0.15)(2.74,-3.31)
\pspolygon[linecolor=qqwuqq,fillcolor=qqwuqq,fillstyle=solid,opacity=0.25](3.73,2.76)(5.8,3.07)(4.39,-2.11)
\pspolygon[linecolor=qqwuqq,fillcolor=qqwuqq,fillstyle=solid,opacity=0.1](3.73,2.76)(5.8,3.07)(-4.77,-8.73)
\pspolygon[linecolor=qqwuqq,fillcolor=qqwuqq,fillstyle=solid,opacity=0.05](5.8,3.07)(4.39,-2.11)(-4.77,-8.73)
\psplot{-7.11}{16.03}{(-39.66--5.42*x)/7.51}
\psplot{-7.11}{16.03}{(-26.13--8.57*x)/7.68}
\psplot{-7.11}{16.03}{(-12.77--7.58*x)/5.61}
\psplot{-7.11}{16.03}{(-3.2--0.99*x)/2.07}
\psplot{-7.11}{16.03}{(--9.21-3.16*x)/-0.17}
\psplot{-7.11}{16.03}{(-0.35-2.16*x)/1.9}
\psplot{-7.11}{16.03}{(--4.57--0.31*x)/2.07}
\psplot{-7.11}{16.03}{(--20-4.88*x)/0.66}
\psplot{-7.11}{16.03}{(--25.73-5.18*x)/-1.41}
\psline[linecolor=qqwuqq](0.84,-1.14)(2.91,-0.15)
\psline[linecolor=qqwuqq](2.91,-0.15)(2.74,-3.31)
\psline[linecolor=qqwuqq](2.74,-3.31)(0.84,-1.14)
\psline[linecolor=qqwuqq](3.73,2.76)(5.8,3.07)
\psline[linecolor=qqwuqq](5.8,3.07)(4.39,-2.11)
\psline[linecolor=qqwuqq](4.39,-2.11)(3.73,2.76)
\psplot[linewidth=2.8pt,linecolor=qqwuqq]{-7.11}{16.03}{(-84.56--9.65*x)/6.51}
\psline[linewidth=1.2pt,linecolor=wqwqwq](4.39,-2.11)(2.4,-9.44)
\psline[linewidth=2pt,linecolor=wqwqwq](2.4,-9.44)(2.74,-3.31)
\psline[linewidth=2pt,linecolor=wqwqwq](4.39,-2.11)(2.4,-9.44)
\psline[linewidth=2pt,linecolor=wqwqwq](4.88,-5.75)(4.39,-2.11)
\psline[linewidth=2pt,linecolor=wqwqwq](4.88,-5.75)(2.74,-3.31)
\psline[linewidth=2pt,linecolor=wqwqwq](3.73,2.76)(11.39,3.9)
\psline[linewidth=2pt,linecolor=wqwqwq](0.84,-1.14)(11.39,3.9)
\psline[linewidth=2.4pt](3.73,2.76)(11.39,3.9)
\psline[linewidth=2.4pt](0.84,-1.14)(11.39,3.9)
\psline[linewidth=2.4pt](3.73,2.76)(4.88,-5.75)
\psline[linewidth=2.4pt](0.84,-1.14)(4.88,-5.75)
\psline[linewidth=2.4pt](2.91,-0.15)(2.4,-9.44)
\psline[linewidth=2.4pt](5.8,3.07)(2.4,-9.44)
\psline[linecolor=qqwuqq](3.73,2.76)(5.8,3.07)
\psline[linecolor=qqwuqq](5.8,3.07)(-4.77,-8.73)
\psline[linecolor=qqwuqq](-4.77,-8.73)(3.73,2.76)
\psline[linecolor=qqwuqq](5.8,3.07)(4.39,-2.11)
\psline[linecolor=qqwuqq](4.39,-2.11)(-4.77,-8.73)
\psline[linecolor=qqwuqq](-4.77,-8.73)(5.8,3.07)
\begin{scriptsize}
\psdots[dotstyle=*,linecolor=zzttqq](2.74,-3.31)
\psdots[dotstyle=*,linecolor=blue](-4.77,-8.73)
\rput[bl](-5.07,-8.61){\blue{$D$}}
\psdots[dotstyle=*,linecolor=blue](2.91,-0.15)
\psdots[dotstyle=*,linecolor=blue](0.84,-1.14)
\psdots[dotstyle=*,linecolor=xdxdff](3.73,2.76)
\psdots[dotstyle=*,linecolor=xdxdff](5.8,3.07)
\psdots[dotstyle=*,linecolor=xdxdff](4.39,-2.11)
\psdots[dotstyle=*,linecolor=darkgray](2.4,-9.44)
\psdots[dotstyle=*,linecolor=darkgray](4.88,-5.75)
\psdots[dotstyle=*,linecolor=darkgray](11.39,3.9)
\rput[bl](3.21,-9.1){\qqwuqq{$d$}}
\end{scriptsize}
\end{pspicture*}
\end{center}

Now
assume $\cP$ is an affine plane, $a$ a line of $\cP$ and $U:=\cP \setminus a$.
For a generic triple $(x,y,z) \in U^3$, we define $w:=(xyz)$ as explained above.
Under which conditions does the identity $((xou)pv)=(xo(upv))$ hold?
To get an idea, reproduce, by using goegebra, the drawing illustrating this identity from
subsection \ref{subsec:Geo-Associativity}, by highlighting four triangles appearing in this
construction, as follows:


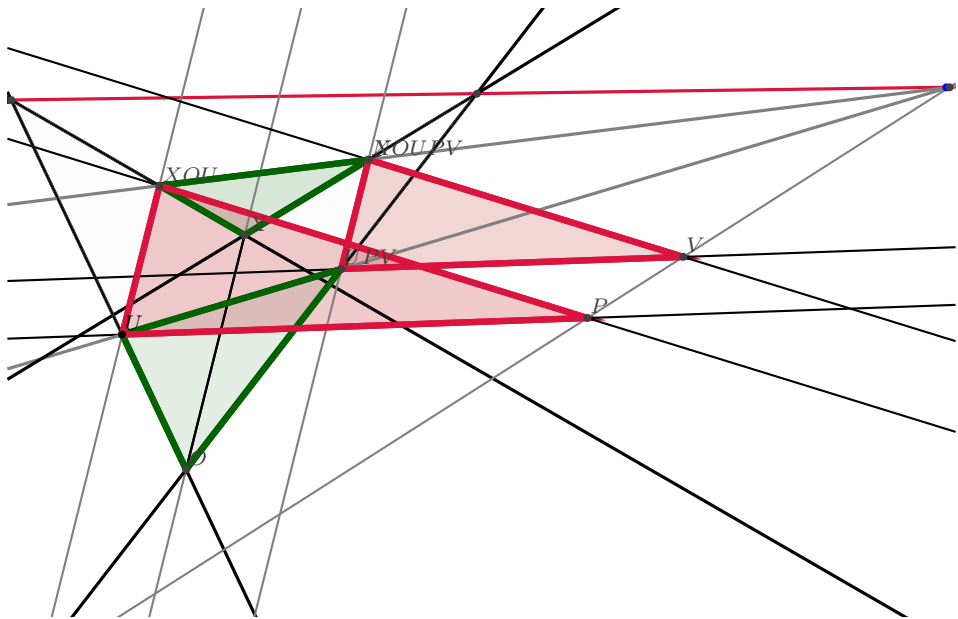
\begin{center}
\newrgbcolor{dcrutc}{0.86 0.08 0.24}
\newrgbcolor{uququq}{0.25 0.25 0.25}
\newrgbcolor{zzttqq}{0.6 0.2 0}
\newrgbcolor{qqwuqq}{0 0.39 0}
\psset{xunit=0.5cm,yunit=0.5cm,algebraic=true,dotstyle=o,dotsize=3pt 0,linewidth=0.8pt,arrowsize=3pt 2,arrowinset=0.25}
\begin{pspicture*}(-4.3,-9.8)(20.64,6.3)
\pspolygon[linecolor=zzttqq,fillcolor=zzttqq,fillstyle=solid,opacity=0.1](-0.3,1.61)(-1.28,-2.34)(10.96,-1.9)
\pspolygon[linecolor=zzttqq,fillcolor=zzttqq,fillstyle=solid,opacity=0.1](4.48,-0.6)(5.2,2.3)(13.48,-0.28)
\pspolygon[linecolor=white,fillcolor=white,fillstyle=solid,opacity=0.1](1.94,0.3)(0.4,-5.92)(8.05,4.05)
\pspolygon[linecolor=white,fillcolor=white,fillstyle=solid,opacity=0.1](1.94,0.3)(0.4,-5.92)(-4.2,3.88)
\pspolygon[fillcolor=black,fillstyle=solid,opacity=0.01](1.94,0.3)(0.4,-5.92)(-4.2,3.88)
\pspolygon[fillcolor=black,fillstyle=solid,opacity=0.01](1.94,0.3)(8.05,4.05)(0.4,-5.92)
\pspolygon[linewidth=2.4pt,linecolor=qqwuqq,fillcolor=qqwuqq,fillstyle=solid,opacity=0.16](1.94,0.3)(5.2,2.3)(-0.3,1.61)
\pspolygon[linewidth=2.4pt,linecolor=qqwuqq,fillcolor=qqwuqq,fillstyle=solid,opacity=0.1](0.4,-5.92)(-1.28,-2.34)(4.48,-0.6)
\pspolygon[linewidth=2.4pt,linecolor=dcrutc,fillcolor=dcrutc,fillstyle=solid,opacity=0.16](-1.28,-2.34)(-0.3,1.61)(10.96,-1.9)
\pspolygon[linewidth=2.4pt,linecolor=dcrutc,fillcolor=dcrutc,fillstyle=solid,opacity=0.1](4.48,-0.6)(5.2,2.3)(13.48,-0.28)
\psplot[linewidth=1.2pt,linecolor=dcrutc]{-4.3}{20.64}{(--97.3--0.34*x)/24.7}
\psplot[linecolor=gray]{-4.3}{20.64}{(--11.6-6.22*x)/-1.54}
\psplot[linewidth=1.2pt]{-4.3}{20.64}{(-8.51-3.58*x)/1.68}
\psplot[linewidth=1.2pt]{-4.3}{20.64}{(-8.79--3.58*x)/-6.14}
\psplot[linecolor=gray]{-4.3}{20.64}{(-4.36-6.22*x)/-1.54}
\psplot[linecolor=gray]{-4.3}{20.64}{(-22.54--1.62*x)/2.52}
\psplot{-4.3}{20.64}{(-28.08--0.44*x)/12.24}
\psplot{-4.3}{20.64}{(--17.05-3.51*x)/11.26}
\psplot{-4.3}{20.64}{(--44.14-3.51*x)/11.26}
\psplot{-4.3}{20.64}{(-9.36--0.44*x)/12.24}
\psplot[linewidth=1.2pt,linecolor=gray]{-4.3}{20.64}{(--34.21--2.61*x)/20.78}
\psplot[linewidth=1.2pt,linecolor=gray]{-4.3}{20.64}{(-42.52--6.56*x)/21.76}
\psplot[linecolor=gray]{-4.3}{20.64}{(--2036.5-439.96*x)/-108.93}
\psplot[linewidth=1.2pt]{-4.3}{20.64}{(--26.28-5.32*x)/-4.08}
\psplot[linewidth=1.2pt]{-4.3}{20.64}{(-5.44--3.75*x)/6.11}
\psline[linecolor=zzttqq](-0.3,1.61)(-1.28,-2.34)
\psline[linecolor=zzttqq](-1.28,-2.34)(10.96,-1.9)
\psline[linecolor=zzttqq](10.96,-1.9)(-0.3,1.61)
\psline[linecolor=zzttqq](4.48,-0.6)(5.2,2.3)
\psline[linecolor=zzttqq](5.2,2.3)(13.48,-0.28)
\psline[linecolor=zzttqq](13.48,-0.28)(4.48,-0.6)
\psline[linecolor=white](1.94,0.3)(0.4,-5.92)
\psline[linecolor=white](0.4,-5.92)(8.05,4.05)
\psline[linecolor=white](8.05,4.05)(1.94,0.3)
\psline[linecolor=white](1.94,0.3)(0.4,-5.92)
\psline[linecolor=white](0.4,-5.92)(-4.2,3.88)
\psline[linecolor=white](-4.2,3.88)(1.94,0.3)
\psline(1.94,0.3)(0.4,-5.92)
\psline(0.4,-5.92)(-4.2,3.88)
\psline(-4.2,3.88)(1.94,0.3)
\psline(1.94,0.3)(8.05,4.05)
\psline(8.05,4.05)(0.4,-5.92)
\psline(0.4,-5.92)(1.94,0.3)
\psline[linewidth=2.4pt,linecolor=qqwuqq](1.94,0.3)(5.2,2.3)
\psline[linewidth=2.4pt,linecolor=qqwuqq](5.2,2.3)(-0.3,1.61)
\psline[linewidth=2.4pt,linecolor=qqwuqq](-0.3,1.61)(1.94,0.3)
\psline[linewidth=2.4pt,linecolor=qqwuqq](0.4,-5.92)(-1.28,-2.34)
\psline[linewidth=2.4pt,linecolor=qqwuqq](-1.28,-2.34)(4.48,-0.6)
\psline[linewidth=2.4pt,linecolor=qqwuqq](4.48,-0.6)(0.4,-5.92)
\psline[linewidth=2.4pt,linecolor=dcrutc](-1.28,-2.34)(-0.3,1.61)
\psline[linewidth=2.4pt,linecolor=dcrutc](-0.3,1.61)(10.96,-1.9)
\psline[linewidth=2.4pt,linecolor=dcrutc](10.96,-1.9)(-1.28,-2.34)
\psline[linewidth=2.4pt,linecolor=dcrutc](4.48,-0.6)(5.2,2.3)
\psline[linewidth=2.4pt,linecolor=dcrutc](5.2,2.3)(13.48,-0.28)
\psline[linewidth=2.4pt,linecolor=dcrutc](13.48,-0.28)(4.48,-0.6)
\begin{scriptsize}
\psdots[dotstyle=*,linecolor=blue](-4.3,3.88)
\psdots[dotstyle=*,linecolor=blue](20.4,4.22)
\psdots[dotstyle=*,linecolor=uququq](1.94,0.3)
\rput[bl](2.02,0.42){\uququq{$X$}}
\psdots[dotstyle=*,linecolor=uququq](0.4,-5.92)
\rput[bl](0.48,-5.8){\uququq{$O$}}
\psdots[dotstyle=*](-1.28,-2.34)
\rput[bl](-1.2,-2.22){$U$}
\psdots[dotstyle=*,linecolor=uququq](-4.2,3.88)
\psdots[dotstyle=*,linecolor=uququq](-0.3,1.61)
\rput[bl](-0.22,1.72){\uququq{$XOU$}}
\psdots[dotstyle=*,linecolor=uququq](10.96,-1.9)
\rput[bl](11.04,-1.78){\uququq{$P$}}
\psdots[dotstyle=*,linecolor=uququq](13.48,-0.28)
\rput[bl](13.56,-0.16){\uququq{$V$}}
\psdots[dotstyle=*,linecolor=uququq](20.48,4.22)
\psdots[dotstyle=*,linecolor=uququq](5.2,2.3)
\rput[bl](5.28,2.42){\uququq{$XOUPV$}}
\psdots[dotstyle=*,linecolor=uququq](4.48,-0.6)
\rput[bl](4.56,-0.48){\uququq{$UPV$}}
\psdots[dotstyle=*,linecolor=uququq](8.05,4.05)
\psdots[dotstyle=*,linecolor=uququq](5.2,2.3)
\rput[bl](5.28,2.42){\uququq{$N$}}
\end{scriptsize}
\end{pspicture*}
\end{center}

Observe that the two green and the two red triangles are in a Desargues' configuration,
respectively (more precisely: the two dual affine versions of the Desargues' configuration
appear here). Use this observation to prove the following theorem:
{\em The general associative law (\ref{eqn:torsor1}) holds in $U=\cP \setminus a$ if, and only if,
Desargues' theorem holds in $\cP$.} 
For the reader familiar with the projective version of Desargues' theorem, we reproduce 
here the projective version of the preceding drawing.
We  have marked the Desargues point $R$ of the ``red triangles'' (whose Desargues line is
$a$) and the Desargues point $G$ of the ``green triangles'' (whose Desargues line is $b$). 


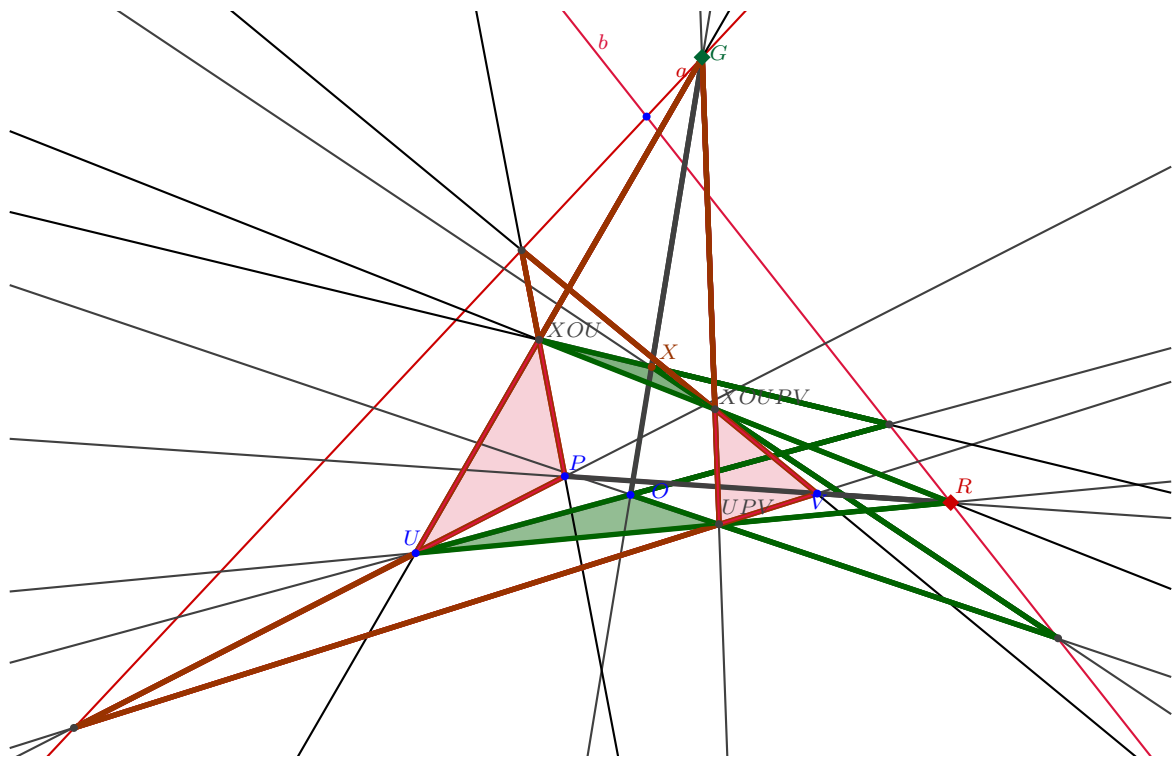
\begin{center}
\newrgbcolor{ccqqqq}{0.8 0 0}
\newrgbcolor{dcrutc}{0.86 0.08 0.24}
\newrgbcolor{zzttqq}{0.6 0.2 0}
\newrgbcolor{uququq}{0.25 0.25 0.25}
\newrgbcolor{qqwwtt}{0 0.4 0.2}
\newrgbcolor{qqwuqq}{0 0.39 0}
\psset{xunit=0.6cm,yunit=0.6cm,algebraic=true,dotstyle=o,dotsize=3pt 0,linewidth=0.8pt,arrowsize=3pt 2,arrowinset=0.25}
\begin{pspicture*}(-7.51,-10.45)(17.95,5.98)
\pspolygon[linecolor=qqwuqq,fillcolor=qqwuqq,fillstyle=solid,opacity=0.1](6.1,-4.7)(1.39,-5.99)(8.04,-5.35)
\pspolygon[linecolor=qqwuqq,fillcolor=qqwuqq,fillstyle=solid,opacity=0.1](6.56,-1.87)(4.09,-1.27)(7.95,-2.81)
\pspolygon[linecolor=dcrutc,fillcolor=dcrutc,fillstyle=solid,opacity=0.1](1.39,-5.99)(4.66,-4.29)(4.09,-1.27)
\pspolygon[linecolor=dcrutc,fillcolor=dcrutc,fillstyle=solid,opacity=0.1](7.95,-2.81)(8.04,-5.35)(10.19,-4.67)
\pspolygon[linecolor=dcrutc,fillcolor=dcrutc,fillstyle=solid,opacity=0.1](4.09,-1.27)(1.39,-5.99)(4.66,-4.29)
\pspolygon[linecolor=dcrutc,fillcolor=dcrutc,fillstyle=solid,opacity=0.1](7.95,-2.81)(8.04,-5.35)(10.19,-4.67)
\pspolygon[linecolor=qqwuqq,fillcolor=qqwuqq,fillstyle=solid,opacity=0.36](1.39,-5.99)(6.1,-4.7)(8.04,-5.35)
\pspolygon[linecolor=qqwuqq,fillcolor=qqwuqq,fillstyle=solid,opacity=0.1](4.09,-1.27)(6.56,-1.87)(7.95,-2.81)
\pspolygon[linecolor=qqwuqq,fillcolor=qqwuqq,fillstyle=solid,opacity=0.38](4.09,-1.27)(6.56,-1.87)(7.95,-2.81)
\psplot[linecolor=ccqqqq]{-7.51}{17.95}{(-27.8--9.11*x)/8.46}
\psplot[linecolor=dcrutc]{-7.51}{17.95}{(--44.78-4.81*x)/3.76}
\psplot[linecolor=uququq]{-7.51}{17.95}{(--30.02-1.29*x)/-4.72}
\psplot[linecolor=uququq]{-7.51}{17.95}{(-19.36--2.82*x)/0.46}
\psplot[linecolor=uququq]{-7.51}{17.95}{(-21.96-0.38*x)/5.52}
\psplot[linecolor=uququq]{-7.51}{17.95}{(-21.96--1.69*x)/3.28}
\psplot[linecolor=uququq]{-7.51}{17.95}{(-128.83--5.18*x)/16.29}
\psplot[linecolor=uququq]{-7.51}{17.95}{(--71.78-1.12*x)/-11.73}
\psplot[linecolor=uququq]{-7.51}{17.95}{(-5.08-0.66*x)/1.94}
\psplot[linecolor=uququq]{-7.51}{17.95}{(--80.95-10.32*x)/0.37}
\psplot[linecolor=uququq]{-7.51}{17.95}{(--22.64-6*x)/8.91}
\psline[linecolor=qqwuqq](6.1,-4.7)(1.39,-5.99)
\psline[linecolor=qqwuqq](1.39,-5.99)(8.04,-5.35)
\psline[linecolor=qqwuqq](8.04,-5.35)(6.1,-4.7)
\psplot{-7.51}{17.95}{(--52.77-10.95*x)/-6.28}
\psplot{-7.51}{17.95}{(-1.45-1.27*x)/5.21}
\psplot{-7.51}{17.95}{(--11.62-3.02*x)/0.57}
\psplot{-7.51}{17.95}{(--24.55-5.38*x)/6.47}
\psplot{-7.51}{17.95}{(-3.2--3.59*x)/-9.03}
\psline[linecolor=qqwuqq](6.56,-1.87)(4.09,-1.27)
\psline[linecolor=qqwuqq](4.09,-1.27)(7.95,-2.81)
\psline[linecolor=qqwuqq](7.95,-2.81)(6.56,-1.87)
\psline[linecolor=dcrutc](1.39,-5.99)(4.66,-4.29)
\psline[linecolor=dcrutc](4.66,-4.29)(4.09,-1.27)
\psline[linecolor=dcrutc](4.09,-1.27)(1.39,-5.99)
\psline[linecolor=dcrutc](7.95,-2.81)(8.04,-5.35)
\psline[linecolor=dcrutc](8.04,-5.35)(10.19,-4.67)
\psline[linecolor=dcrutc](10.19,-4.67)(7.95,-2.81)
\psline[linewidth=1.2pt,linecolor=qqwuqq](7.67,4.97)(1.39,-5.99)
\psline[linewidth=1.2pt,linecolor=qqwuqq](7.67,4.97)(4.09,-1.27)
\psline[linewidth=1.2pt,linecolor=qqwuqq](7.67,4.97)(7.95,-2.81)
\psline[linewidth=1.2pt,linecolor=qqwuqq](7.67,4.97)(8.04,-5.35)
\psline[linewidth=1.2pt,linecolor=qqwuqq](7.67,4.97)(6.56,-1.87)
\psline[linewidth=1.2pt,linecolor=qqwuqq](7.67,4.97)(6.1,-4.7)
\psline[linewidth=1.6pt,linecolor=qqwuqq](7.67,4.97)(8.04,-5.35)
\psline[linewidth=1.6pt,linecolor=qqwuqq](7.67,4.97)(6.1,-4.7)
\psline[linewidth=1.6pt,linecolor=qqwuqq](7.67,4.97)(1.39,-5.99)
\psline[linewidth=1.6pt,linecolor=dcrutc](1.39,-5.99)(13.12,-4.87)
\psline[linewidth=1.6pt,linecolor=dcrutc](8.04,-5.35)(13.12,-4.87)
\psline[linewidth=1.6pt,linecolor=dcrutc](4.09,-1.27)(13.12,-4.87)
\psline[linewidth=1.6pt,linecolor=dcrutc](7.95,-2.81)(13.12,-4.87)
\psline[linewidth=1.6pt,linecolor=dcrutc](13.12,-4.87)(4.66,-4.29)
\psline[linewidth=1.6pt,linecolor=dcrutc](10.19,-4.67)(13.12,-4.87)
\psline[linewidth=2pt,linecolor=uququq](7.67,4.97)(8.04,-5.35)
\psline[linewidth=2pt,linecolor=uququq](7.67,4.97)(7.95,-2.81)
\psline[linewidth=2pt,linecolor=uququq](7.67,4.97)(6.56,-1.87)
\psline[linewidth=2pt,linecolor=uququq](7.67,4.97)(6.1,-4.7)
\psline[linewidth=2pt,linecolor=uququq](7.67,4.97)(4.09,-1.27)
\psline[linewidth=2pt,linecolor=uququq](7.67,4.97)(1.39,-5.99)
\psline[linewidth=2pt,linecolor=uququq](1.39,-5.99)(11.77,-3.14)
\psline[linewidth=2pt,linecolor=uququq](6.1,-4.7)(11.77,-3.14)
\psline[linewidth=2pt,linecolor=uququq](11.77,-3.14)(6.56,-1.87)
\psline[linewidth=2pt,linecolor=uququq](11.77,-3.14)(4.09,-1.27)
\psline[linewidth=2pt,linecolor=uququq](6.56,-1.87)(15.47,-7.87)
\psline[linewidth=2pt,linecolor=uququq](7.95,-2.81)(15.47,-7.87)
\psline[linewidth=2pt,linecolor=uququq](6.1,-4.7)(15.47,-7.87)
\psline[linewidth=2pt,linecolor=qqwuqq](8.04,-5.35)(15.47,-7.87)
\psline[linewidth=2pt,linecolor=qqwuqq](6.1,-4.7)(15.47,-7.87)
\psline[linewidth=2pt,linecolor=qqwuqq](7.95,-2.81)(15.47,-7.87)
\psline[linewidth=2pt,linecolor=qqwuqq](6.56,-1.87)(15.47,-7.87)
\psline[linewidth=2pt,linecolor=qqwuqq](1.39,-5.99)(11.77,-3.14)
\psline[linewidth=2pt,linecolor=qqwuqq](6.1,-4.7)(11.77,-3.14)
\psline[linewidth=2pt,linecolor=qqwuqq](4.09,-1.27)(11.77,-3.14)
\psline[linewidth=2pt,linecolor=qqwuqq](6.56,-1.87)(11.77,-3.14)
\psline[linewidth=2pt,linecolor=dcrutc](3.71,0.71)(4.66,-4.29)
\psline[linewidth=2pt,linecolor=dcrutc](3.71,0.71)(4.09,-1.27)
\psline[linewidth=2pt,linecolor=dcrutc](7.67,4.97)(4.09,-1.27)
\psline[linewidth=2pt,linecolor=dcrutc](7.67,4.97)(1.39,-5.99)
\psline[linewidth=2pt,linecolor=dcrutc](-6.1,-9.85)(4.66,-4.29)
\psline[linewidth=2pt,linecolor=dcrutc](-6.1,-9.85)(1.39,-5.99)
\psline[linewidth=2pt,linecolor=dcrutc](-6.1,-9.85)(10.19,-4.67)
\psline[linewidth=2pt,linecolor=dcrutc](-6.1,-9.85)(8.04,-5.35)
\psline[linewidth=2pt,linecolor=dcrutc](3.71,0.71)(7.95,-2.81)
\psline[linewidth=2pt,linecolor=dcrutc](3.71,0.71)(10.19,-4.67)
\psline[linewidth=2pt,linecolor=dcrutc](7.67,4.97)(7.95,-2.81)
\psline[linewidth=2pt,linecolor=dcrutc](7.67,4.97)(8.04,-5.35)
\psline[linewidth=2pt,linecolor=uququq](4.66,-4.29)(13.12,-4.87)
\psline[linewidth=2pt,linecolor=uququq](10.19,-4.67)(13.12,-4.87)
\psline[linewidth=2pt,linecolor=qqwuqq](4.09,-1.27)(13.12,-4.87)
\psline[linewidth=2pt,linecolor=qqwuqq](7.95,-2.81)(13.12,-4.87)
\psline[linewidth=2pt,linecolor=qqwuqq](8.04,-5.35)(13.12,-4.87)
\psline[linewidth=2pt,linecolor=qqwuqq](13.12,-4.87)(1.39,-5.99)
\psline[linewidth=2pt,linecolor=zzttqq](-6.1,-9.85)(8.04,-5.35)
\psline[linewidth=2pt,linecolor=zzttqq](-6.1,-9.85)(10.19,-4.67)
\psline[linewidth=2pt,linecolor=zzttqq](-6.1,-9.85)(1.39,-5.99)
\psline[linewidth=2pt,linecolor=zzttqq](-6.1,-9.85)(4.66,-4.29)
\psline[linewidth=2pt,linecolor=zzttqq](3.71,0.71)(4.09,-1.27)
\psline[linewidth=2pt,linecolor=zzttqq](3.71,0.71)(4.66,-4.29)
\psline[linewidth=2pt,linecolor=zzttqq](3.71,0.71)(7.95,-2.81)
\psline[linewidth=2pt,linecolor=zzttqq](3.71,0.71)(10.19,-4.67)
\psline[linewidth=2pt,linecolor=zzttqq](7.67,4.97)(4.09,-1.27)
\psline[linewidth=2pt,linecolor=zzttqq](7.67,4.97)(1.39,-5.99)
\psline[linewidth=2pt,linecolor=zzttqq](7.67,4.97)(7.95,-2.81)
\psline[linewidth=2pt,linecolor=zzttqq](7.67,4.97)(8.04,-5.35)
\psline[linecolor=dcrutc](4.09,-1.27)(1.39,-5.99)
\psline[linecolor=dcrutc](1.39,-5.99)(4.66,-4.29)
\psline[linecolor=dcrutc](4.66,-4.29)(4.09,-1.27)
\psline[linecolor=dcrutc](7.95,-2.81)(8.04,-5.35)
\psline[linecolor=dcrutc](8.04,-5.35)(10.19,-4.67)
\psline[linecolor=dcrutc](10.19,-4.67)(7.95,-2.81)
\psline[linecolor=qqwuqq](1.39,-5.99)(6.1,-4.7)
\psline[linecolor=qqwuqq](6.1,-4.7)(8.04,-5.35)
\psline[linecolor=qqwuqq](8.04,-5.35)(1.39,-5.99)
\psline[linecolor=qqwuqq](4.09,-1.27)(6.56,-1.87)
\psline[linecolor=qqwuqq](6.56,-1.87)(7.95,-2.81)
\psline[linecolor=qqwuqq](7.95,-2.81)(4.09,-1.27)
\psline[linecolor=qqwuqq](4.09,-1.27)(6.56,-1.87)
\psline[linecolor=qqwuqq](6.56,-1.87)(7.95,-2.81)
\psline[linecolor=qqwuqq](7.95,-2.81)(4.09,-1.27)
\begin{scriptsize}
\psdots[dotstyle=*,linecolor=blue](6.45,3.66)
\rput[bl](7.09,4.53){\ccqqqq{$a$}}
\rput[bl](5.39,5.15){\dcrutc{$b$}}
\psdots[dotstyle=*,linecolor=blue](6.1,-4.7)
\rput[bl](6.56,-4.78){\blue{$O$}}
\psdots[dotstyle=*,linecolor=blue](1.39,-5.99)
\rput[bl](1.12,-5.8){\blue{$U$}}
\psdots[dotstyle=*,linecolor=zzttqq](6.56,-1.87)
\rput[bl](6.72,-1.69){\zzttqq{$X$}}
\psdots[dotstyle=*,linecolor=blue](4.66,-4.29)
\rput[bl](4.74,-4.16){\blue{$P$}}
\psdots[dotstyle=*,linecolor=blue](10.19,-4.67)
\rput[bl](10.03,-5.02){\blue{$V$}}
\psdots[dotsize=5pt 0,dotstyle=square*,dotangle=45,linecolor=ccqqqq](13.12,-4.87)
\rput[bl](13.21,-4.65){\ccqqqq{$R$}}
\psdots[dotstyle=*,linecolor=darkgray](-6.1,-9.85)
\psdots[dotstyle=*,linecolor=darkgray](8.04,-5.35)
\rput[bl](8.09,-5.12){\darkgray{$UPV$}}
\psdots[dotstyle=*,linecolor=darkgray](15.47,-7.87)
\psdots[dotsize=5pt 0,dotstyle=square*,dotangle=45,linecolor=qqwwtt](7.67,4.97)
\rput[bl](7.84,4.9){\qqwwtt{$G$}}
\psdots[dotstyle=*,linecolor=darkgray](7.95,-2.81)
\rput[bl](8.03,-2.69){\darkgray{$XOUPV$}}
\psdots[dotstyle=*,linecolor=darkgray](11.77,-3.14)
\psdots[dotstyle=*,linecolor=darkgray](4.09,-1.27)
\rput[bl](4.25,-1.22){\darkgray{$XOU$}}
\psdots[dotstyle=*,linecolor=darkgray](3.71,0.71)
\end{scriptsize}
\end{pspicture*}
\end{center}

Summing up, a projective plane is Desarguesian if, and only if,  all laws $(xyz)_{ab}$ 
are associative.
As mentioned in the introduction, the question arises then what can be said for general
projective planes, or even for more general lattices.


\bigskip

Wolfgang Bertram\\
Institut \'{E}lie Cartan de Lorraine (IECL) \\
Universit\'{e} de Lorraine, CNRS, INRIA, \\
Boulevard des Aiguillettes, B.P. 239, \\
F-54506 Vand\oe{}uvre-l\`{e}s-Nancy, France\\
{\tt wolfgang.bertram@univ-lorraine.fr}

\end{document}